\documentclass[twoside,11pt,reqno]{amsart}
\usepackage{amsmath,amssymb,amscd,mathrsfs,epic,empheq}
\usepackage{latexsym,amsthm,amscd, enumerate}
\usepackage[all]{xy}
\usepackage{stmaryrd}

\makeatletter

\@addtoreset{equation}{section}

\newtheorem{Proposition}{Proposition}[section]
\newtheorem{Lemma}[Proposition]{Lemma}
\newtheorem{Theorem}[Proposition]{Theorem}
\newtheorem{Corollary}[Proposition]{Corollary}

\newtheorem{Remark}[Proposition]{Remark}

\def\phantomsubsection#1{\vspace{2mm}\noindent{\bf #1.}}

\newcommand{\cV}{\mathcal{V}}
\newcommand{\nV}{V}
\newcommand{\cL}{\mathcal{L}}
\newcommand{\cD}{D}

\newcommand{\fin}{{f\:\!\!in}}

\def\op{{\operatorname{op}}}

\newcommand{\mg}{\mathfrak{g}}

\newcommand{\cP}{\mathcal{P}}

\newcommand{\sF}{\mathscr{F}\!}
\newcommand{\sFp}{\mathscr{F}'}
\newcommand{\cO}{\mathcal{O}}

\newcommand{\cF}{\mathcal F}
\newcommand{\cE}{\mathcal E}

\newcommand{\nT}{T}
\newcommand{\nP}{P}

\newcommand{\mZ}{\mathbb{Z}}

\newbox\squ  % box character for ends of proofs
\setbox\squ=\hbox{\vrule width.6pt
   \vbox{\hrule height.6pt width.4em\kern1ex
         \hrule height.6pt}%
   \vrule width.6pt}

\def\dd{\operatorname{dim}_{p,q}}

\def\zi{\xi}
\def\deg{\operatorname{ht}}
\def\defect{\operatorname{def}}

\def\op{\operatorname{op}}

\def\down{\vee}
\def\up{\wedge}
\def\Rep#1{\operatorname{rep}(#1)}
\def\id{\operatorname{id}}

\def\Id{\operatorname{Id}}
\def\C{{\mathbb F}}

\def\Z{{\mathbb Z}}

\def\0{{\bar 0}}
\def\1{{\bar 1}}
\def\pr{{\operatorname{pr}}}
\def\hom{{\operatorname{Hom}}}

\def\End{{\operatorname{End}}}

\def\soc{{\operatorname{soc}\:}}

\def\bid{\hbox{\boldmath{$1$}}}
\def\eps{{\varepsilon}}
\def\phi{{\varphi}}

\def\la{{\lambda}}
\def\La{{\Lambda}}
\def\Ga{{\Gamma}}
\def\ga{{\gamma}}

\def\De{{\Delta}}

\def\bh{\text{\boldmath$h$}}
\def\bi{\text{\boldmath$i$}}
\def\bj{\text{\boldmath$j$}}

\def\bGa{\text{\boldmath$\Ga$}}

\def\bt{\text{\boldmath$t$}}
\def\bs{\text{\boldmath$s$}}
\def\br{\text{\boldmath$r$}}
\def\bu{\text{\boldmath$u$}}

\def\bde{\text{\boldmath$\delta$}}
\def\bga{\text{\boldmath$\ga$}}

\def\bsigma{\text{\boldmath$\sigma$}}
\def\btau{\text{\boldmath$\tau$}}

\begin{document}

\title[Khovanov's diagram algebra IV]{\boldmath Highest weight categories
arising from Khovanov's diagram algebra IV: the general linear supergroup}
\author{Jonathan Brundan and Catharina Stroppel}

\address{Department of Mathematics, University of Oregon, Eugene, OR 97403, USA}
\email{brundan@uoregon.edu}
\address{Department of Mathematics, University of Bonn, 53115 Bonn, Germany}
\email{stroppel@math.uni-bonn.de}

\thanks{First author supported in part by NSF grant no. DMS-0654147.}

\begin{abstract}
We prove that blocks of
the general linear
supergroup are Morita
equivalent to
a limiting version of Khovanov's diagram algebra.
We deduce that blocks of the general linear supergroup are Koszul.
\end{abstract}
\maketitle

% TOC seems not so useful as the paper doesn't have very many
% sections; could be removed

\tableofcontents

\section{Introduction}

This is the culmination of a series of four articles
studying various generalisations of Khovanov's diagram algebra from \cite{K2}. The goal is to relate
the limiting version $H^\infty_r$ of this algebra constructed in \cite{BS1}
to blocks of the general linear supergroup $GL(m|n)$.
More precisely,
working always
over a fixed algebraically closed field $\C$ of characteristic zero,
we show that
any block of $GL(m|n)$ of atypicality $r$ is Morita equivalent to the
algebra $H^\infty_r$.
We refer the reader to the introduction of \cite{BS1}
for a detailed account of our approach to the definition of Khovanov's diagram
algebra and the construction of its limiting version; see also \cite{SICM} which discusses
further the connections to link homology.

To formulate our main result in detail, fix $m, n \geq 0$ and let
$G$ denote the algebraic supergroup $GL(m|n)$ over $\C$. 
Using scheme-theoretic language, $G$ can be regarded as a functor
from the category of commutative superalgebras over $\C$
to the category of groups,
mapping a commutative superalgebra $A = A_{\0}\oplus A_{\1}$
to the group $G(A)$ of all  invertible $(m+n) \times (m+n)$
matrices of the form
\begin{equation}\label{supermat}
g=
\left(
\begin{array}{l|l}
a&b\\\hline
c&d
\end{array}
\right)
\end{equation}
where $a$ (resp. $d$) is an $m \times m$ (resp. $n \times n$)
matrix with entries in
$A_{\0}$, and $b$ (resp. $c$) is an $m\times n$ (resp. $n \times m$)
matrix with entries in
$A_{\1}$.

We are interested here in finite dimensional
representations of $G$, which can be viewed equivalently as
{integrable} supermodules over its
Lie superalgebra $\mathfrak{g} \cong \mathfrak{gl}(m|n,\C)$; the
condition for integrability is the same as for $\mathfrak{g}_{\0}
\cong
\mathfrak{gl}(m,\C) \oplus \mathfrak{gl}(n,\C)$.\label{nv}
For example, we have the natural $G$-module
$\nV$ of column vectors, with
standard basis $v_1,\dots,v_m,v_{m+1},\dots,v_{m+n}$
and $\Z_2$-grading defined by putting $v_r$ in degree $\bar r := \0$ if $1 \leq r \leq m$, $\bar r := \1$ if $m+1 \leq r \leq m+n$.
Let $B$ and $T$ be the standard choices of Borel subgroup\label{borels}
and maximal torus: for each commutative superalgebra $A$, the
groups
$B(A)$ and $T(A)$ consist of all
matrices $g \in G(A)$ that are upper triangular and diagonal, respectively.
Let $\eps_1,\dots,\eps_{m+n}$ be the usual basis for the
character group $X(T)$ of $T$, i.e. $\eps_r$ picks
out the $r$th diagonal entry of a diagonal matrix.
Equip $X(T)$ with a symmetric bilinear form $(.,.)$
such that $(\eps_r,\eps_s) = (-1)^{\bar r} \delta_{r,s}$,
and set
\begin{equation}
\rho :=
\sum_{r=1}^m (1-r) \eps_r + \sum_{s=1}^n (m-s) \eps_{m+s}.
\end{equation}
 Let
\begin{equation}
\label{XT}
X^+(T)
:= \left\{
\la \in X(T)\:\bigg|\:
\begin{array}{c}
\:\:\,(\la+\rho,\eps_1) > \cdots > (\la+\rho,\eps_m),\\
(\la+\rho,\eps_{m+1}) < \cdots < (\la+\rho,\eps_{m+n})
\end{array}
\right\}
\end{equation}
denote the set of {\em dominant weights}.

We allow
only {\em even} morphisms between $G$-modules, so that the
category of all finite dimensional
$G$-modules is obviously an abelian category.
Any $G$-module $M$ decomposes as
$M = M_+ \oplus M_-$,
where
$M_+$ (resp.\ $M_-$) is the $G$-submodule of $M$ spanned by
the degree
$\bar\lambda$
(resp.\ the degree $(\bar\lambda+\bar1)$)
component of the $\lambda$-weight space of $M$ for all $\lambda \in X(T)$;
here $\bar\lambda :=
(\lambda,\eps_{m+1}+\cdots+\eps_{m+n}) \pmod{2}$.
It follows that the category of all finite dimensional
$G$-modules
decomposes as $\sF \oplus \Pi \sF$, where
$\sF = \sF(m|n)$\label{fcat}
(resp.\ $\Pi \sF = \Pi \sF(m|n)$) is the full subcategory consisting
of all $M$ such that $M = M_+$ (resp.\ $M = M_-$).
Moreover $\sF$ and $\Pi \sF$ are obviously equivalent.
In view of this decomposition, we will focus just on $\sF$ from now on.
Note further that
$\sF$ is closed under tensor product, and it contains
both the natural module $\nV$ and its dual $\nV^*$.

By \cite[Theorem 4.47]{B}, the category $\sF$ is
a highest weight category with
weight poset $(X^+(T), \leq)$, where $\leq$ is the
{\em Bruhat ordering} defined combinatorially in
the next paragraph.
We fix representatives
$\{\cL(\la)\:|\:\la \in X^+(T)\}$
for the isomorphism classes of
irreducible modules in $\sF$ so that
$\cL(\la)$ is an
irreducible object in $\sF$
generated by a one-dimensional $B$-submodule of
weight $\la$.
We also denote the
{standard}
and {projective indecomposable modules}
in the highest weight category $\sF$
by
$\{\cV(\la)\:|\:\la \in X^+(T)\}$ and
$\{\cP(\la)\:|\:\la \in X^+(T)\}$, respectively.
So $\cP(\la) \twoheadrightarrow \cV(\la) \twoheadrightarrow \cL(\la)$.
In this setting, the standard module $\cV(\la)$ is often referred to as a {\em Kac module}
after \cite{Kac2}.

Now we turn our attention to the diagram algebra side.\label{dwts}
Let $\La = \La(m|n)$ denote the set of all weights in the
diagrammatic sense of \cite[$\S$2]{BS1} drawn on a number line with
vertices indexed by $\Z$, such that a total of $m$ vertices are
labelled $\times$ or $\down$, a total of $n$ vertices
are labelled $\circ$ or $\down$, and all of the (infinitely many)
remaining  vertices are labelled $\up$.
From now on, we identify the set $X^+(T)$ introduced above with the set
$\La$ via the following {\em weight dictionary}.
Given $\la \in X^+(T)$, we define
\iffalse
\begin{align}
I_{\scriptstyle\times}(\la) &:= \{(\la+\rho,\eps_1),\dots,(\la+\rho,\eps_m)\},\\
I_\circ(\la) &:= \{(\la+\rho,\eps_{m+1}),\dots,(\la+\rho,\eps_{m+n})\}.
\end{align}
Then we identify $\la$ with the element of $\La$
whose $i$th vertex is labelled
\begin{equation}\label{wtdict}
\left\{
\begin{array}{ll}
{\scriptstyle\up}&\text{if $i$ does not belong to either $I_{\scriptstyle\times}(\la)$ or $I_\circ(\la)$,}\\
{\scriptstyle\times}&\text{if $i$ belongs to $I_{\scriptstyle\times}(\la)$ but not to
$I_\circ(\la)$,}\\
\circ&\text{if $i$ belongs to $I_\circ(\la)$ but not to $I_{\scriptstyle\times}(\la)$,}\\
{\scriptstyle\down}&\text{if $i$ belongs to both $I_{\scriptstyle\times}(\la)$ and $I_\circ(\la)$.}
\end{array}\right.
\end{equation}
\fi
\begin{align}
I_{\down}(\la) &:= \{(\la+\rho,\eps_1),\dots,(\la+\rho,\eps_m)\},\\
I_\up(\la) &:= \Z \setminus \{(\la+\rho,\eps_{m+1}),\dots,(\la+\rho,\eps_{m+n})\}.
\end{align}
Then we identify $\la$ with the element of $\La$
whose $i$th vertex is labelled
\begin{equation}\label{wtdict}
\left\{
\begin{array}{ll}
\circ&\text{if $i$ does not belong to either $I_\down(\la)$ or $I_\up(\la)$,}\\
{\scriptstyle\down} &\text{if $i$ belongs to $I_\down(\la)$ but not to
$I_\up(\la)$,}\\
{\scriptstyle\up} &\text{if $i$ belongs to $I_\up(\la)$ but not to
$I_\down(\la)$,}\\
{\scriptstyle\times}&\text{if $i$ belongs to both $I_\down(\la)$ and $I_\up(\la)$.}
\end{array}\right.
\end{equation}
For example, the zero weight (which parametrises the trivial
$G$-module) is identified with the diagram
\begin{equation*}\label{groundstate}
\begin{array}{ll}
\hspace{100mm}
&\text{if $m \geq n$,}\\\\
&\text{if $m \leq n$.}\\
\end{array}
\begin{picture}(0,30)
\put(-333,12.5){$\cdots$}
\put(-90,12.5){$\cdots$}
\put(-270,14.2){$\overbrace{\phantom{hellow worl}}^{n}$}
\put(-201,14.2){$\overbrace{\phantom{hellow worl}}^{m-n}$}
\put(-313,15.2){\line(1,0){214}}
\put(-288.7,10.6){$\scriptstyle\up$}
\put(-311.7,10.6){$\scriptstyle\up$}
\put(-268.7,15.3){$\scriptstyle\down$}
\put(-245.7,15.3){$\scriptstyle\down$}
\put(-222.7,15.3){$\scriptstyle\down$}
\put(-200.2,13.3){$\scriptstyle\times$}
\put(-177.2,13.3){$\scriptstyle\times$}
\put(-154.2,13.3){$\scriptstyle\times$}
\put(-130.7,10.6){$\scriptstyle\up$}
\put(-107.7,10.6){$\scriptstyle\up$}
\end{picture}
\begin{picture}(0,0)
\put(-333,-12.6){$\cdots$}
\put(-90,-12.6){$\cdots$}
\put(-270,-17){$\underbrace{\phantom{hellow worl}}_{m}$}
\put(-201,-17){$\underbrace{\phantom{hellow worl}}_{n-m}$}
\put(-313,-10){\line(1,0){214}}
\put(-288.7,-14.6){$\scriptstyle\up$}
\put(-311.7,-14.6){$\scriptstyle\up$}
\put(-268.7,-9.9){$\scriptstyle\down$}
\put(-245.7,-9.9){$\scriptstyle\down$}
\put(-222.7,-9.9){$\scriptstyle\down$}
\put(-200,-12.6){$\circ$}
\put(-177,-12.6){$\circ$}
\put(-154,-12.6){$\circ$}
\put(-130.7,-14.6){$\scriptstyle\up$}
\put(-107.7,-14.6){$\scriptstyle\up$}
\end{picture}
\end{equation*}
\vspace{4mm}

\noindent
where the leftmost $\down$ is on vertex $(1-m)$.
In these diagrammatic terms, the Bruhat ordering on $X^+(T)$
mentioned earlier is
the same as the Bruhat ordering on $\La$ from \cite[$\S$2]{BS1},
that is, the partial order $\leq$ on diagrams
generated by the basic operation of swapping a $\down$ and an $\up$
so that $\down$'s move to the right.

Let $\sim$ be the
equivalence relation on $\La$
generated by permuting $\down$'s and $\up$'s.
Following the language of \cite{BS1} again,
the $\sim$-equivalence classes of weights from $\La$
are called {\em blocks}. The {\em defect} $\defect(\Ga)$
of each block $\Ga \in \La / \sim$ is simply equal to the number
of vertices labelled $\down$ in any weight $\la \in\Ga$;
this is the same thing as the usual notion of {\em atypicality} in the
representation theory of $GL(m|n)$
as in e.g. \cite[(1.1)]{Serg}.

Let $K = K(m|n)$ denote the direct sum
of the diagram algebras
$K_\Ga$ associated to all the blocks $\Ga \in \La / \sim$
as defined in \cite[$\S$4]{BS1}.
As a vector space, $K$ has a basis
\begin{equation}\label{hbase}
\left\{(a \la b)\:|\:\text{for all oriented circle diagrams
$a \la b$ with $\la \in \La$}\right\},
\end{equation}
and its multiplication is defined by an explicit combinatorial procedure
in terms of such diagrams as in \cite[$\S$6]{BS1}; see the discussion
at the end of this introduction for some examples illustrating the precise
meaning of all this.
As explained in \cite[$\S$5]{BS1}, to each $\la \in \La$
there is associated
an idempotent
$e_\la \in K$. The left ideal
$P(\la) := K e_\la$ is a projective
indecomposable module with irreducible head denoted $L(\la)$.
The modules $\{L(\la)\:|\:\la \in \La\}$ are all one dimensional and give a complete set of
irreducible $K$-modules. Finally let $V(\la)$ be the
standard module corresponding to $\la$, which was
referred to as a {\em cell module}
in \cite[$\S$5]{BS1}.

The main result of the paper is the following.

\begin{Theorem}\label{main}
\iffalse
Let $\nP := \bigoplus_{\la \in \La(m|n)} \cP(\la)$. Let
$\operatorname{End}^{f.d.}_G(P)$ denote the algebra of all $G$-module
endomorphisms of $P$ with finite dimensional image.
Then there is an algebra isomorphism
$$
i:K(m|n) \stackrel{\sim}{\rightarrow}\End^{f.d.}_G(\nP)^{\op}
$$
mapping $e_\la$ to the projection of $\nP$ onto
the summand $\cP(\la)$.
Hence, viewing $\nP$ as a right $K(m|n)$-module via this isomorphism,
\fi
There is an equivalence of categories $\mathbb E$ from
 $\sF(m|n)$ to the category of finite
dimensional left $K(m|n)$-modules,
such that $\mathbb E \cL(\la) \cong L(\la)$,
$\mathbb E \cV(\la) \cong V(\la)$
and
$\mathbb E \cP(\la) \cong P(\la)$ for each $\la \in \La(m|n)$.
\end{Theorem}

Our proof of Theorem~\ref{main} involves showing that
$K$ is isomorphic to the locally finite endomorphism
algebra
$\End_G^{fin}(P)^{\op}$
of a canonical minimal projective generator $P \cong \bigoplus_{\la\in\La}\cP(\la)$
for $\sF$; see Lemmas~\ref{l1}--\ref{l2} below.
To construct $P$, we first consider the weight
\begin{equation}\label{laab}
\la_{p,q} := \sum_{r=1}^m p \eps_r - \sum_{s=1}^n (q+m) \eps_{m+s}
\end{equation}
for integers $p \leq q$.
This is represented diagrammatically by
\begin{equation}
\begin{picture}(-320,25)
\put(-222,14){$_p$}
\put(-153.3,14){$_q$}
\put(-294,14){$_{p-m}$}
\put(-90,14){$_{q+n}$}
\put(-333,-1.5){$\cdots$}
\put(-22,-1.5){$\cdots$}
\put(-270,-5.8){$\underbrace{\phantom{hellow worl}}_{m}$}
\put(-133,-5.8){$\underbrace{\phantom{hellow worl}}_{n}$}
\put(-313,1.2){\line(1,0){283}}
\put(-288.7,-3.3){$\scriptstyle\up$}
\put(-311.7,-3.3){$\scriptstyle\up$}
\put(-153.7,-3.3){$\scriptstyle\up$}
\put(-176.7,-3.3){$\scriptstyle\up$}
\put(-199.7,-3.3){$\scriptstyle\up$}
\put(-269.2,-0.6){$\scriptstyle\times$}
\put(-246.2,-0.6){$\scriptstyle\times$}
\put(-223.2,-0.6){$\scriptstyle\times$}
\put(-131.4,-1.4){$\circ$}
\put(-108.4,-1.4){$\circ$}
\put(-85.4,-1.4){$\circ$}
\put(-61.7,-3.3){$\scriptstyle\up$}
\put(-38.7,-3.3){$\scriptstyle\up$}
\end{picture}\label{tuesday}
\end{equation}

\vspace{5mm}

\noindent
where the rightmost $\times$ is on vertex $p$ and the
rightmost $\circ$ is on vertex $(q+n)$.
The $G$-module $\cV(\la_{p,q})$ is
projective,
hence the ``tensor space''
$\cV(\la_{p,q}) \otimes \nV^{\otimes d}$ is projective for any $d \geq 0$.
Moreover, any $\cP(\la)$ appears as a summand of
$\cV(\la_{p,q}) \otimes \nV^{\otimes d}$ for suitable $p, q$ and $d$.
The key step in our approach is to compute the endomorphism algebra
of
$\cV(\la_{p,q}) \otimes \nV^{\otimes d}$ for $d \geq 0$.
For $d \leq \min(m,n)$, we show that it is a certain degenerate cyclotomic Hecke algebra of level two,
giving a new ``super'' version of the level two Schur-Weyl duality from \cite{BKschur}.
Then we invoke results from \cite{BS3} which show that the basic
algebra that is
Morita equivalent to
this
cyclotomic Hecke algebra is a
generalised Khovanov algebra; this equivalence relies in particular on the connection between
cyclotomic Hecke algebras and Khovanov-Lauda-Rouquier algebras in type $A$
from \cite{BKinv}.
Finally we let $p, q$ and $d$ vary, taking a suitable direct limit
to derive our main result.

We briefly collect here some applications of Theorem~\ref{main}.

\phantomsubsection{Blocks of the same atypicality are equivalent}
The algebras $K_\Ga$
for all $\Ga \in \La(m|n) / \sim$ are the
blocks of the algebra $K(m|n)$.
Hence by Theorem~\ref{main}
they are the
basic algebras representing the individual blocks of the category
$\sF(m|n)$.
In the diagrammatic setting, it is obvious
for $\Ga \in \La(m|n) / \sim$ and $\Ga' \in \La(m'|n') / \sim$
(for possibly different $m'$ and $n'$)
that the algebras $K_\Ga$ and $K_{\Ga'}$ are isomorphic if
and only if $\Ga$ and $\Ga'$ have the same defect. Thus we recover
a result of Serganova from \cite{Serg3}: the blocks of
$GL(m|n)$ for all $m, n$ depend up to equivalence only on the
degree of atypicality of the block.

\phantomsubsection{Gradings on blocks and Koszulity}
Each of the algebras $K_\Ga$ carries a canonical positive grading with respect to which
it is a (locally unital) Koszul algebra; see \cite[Corollary 5.13]{BS2}.
So Theorem~\ref{main} implies that blocks of $GL(m|n)$ are Koszul.
The appearence of such hidden Koszul gradings in representation theory goes back to the
classic paper of Beilinson, Ginzburg and Soergel \cite{BGS}
on blocks of category $\cO$ for a semisimple Lie algebra.
In that work, the grading is of geometric origin,
whereas in our situation we establish the Koszulity in a purely algebraic way.

\phantomsubsection{Rigidity of Kac modules}
Another consequence of Theorem~\ref{main}, combined with \cite[Corollary 6.7]{BS2} on the diagram algebra side,
is that all the Kac modules $\cV(\la)$ are {\em rigid}, i.e. their
radical and socle filtrations coincide.
See \cite[Theorem 5.2]{BS1} for the explicit combinatorial description of the layers.

\phantomsubsection{Kostant modules and BGG resolution}
In \cite{BS2} we studied in detail the {\em Kostant
  modules} for the generalised Khovanov
algebras,  i.e. the irreducible modules whose Kazhdan-Lusztig polynomials
are multiplicity-free. In particular in \cite[Lemma 7.2]{BS2} we
classified the highest weights of these modules
via a pattern avoidance condition.  Combining this with Theorem~\ref{main},
we obtain the following
classification of all Kostant modules for $GL(m|n)$: they are the irreducible modules parametrised
by the weights
in which no two vertices labelled $\down$
have a vertex labelled $\up$ between them.
By \cite[Theorem 7.3]{BS2},
Kostant modules possess a BGG
resolution by multiplicity-free direct sums of standard modules.
All irreducible polynomial representations of $GL(m|n)$ satisfy the
combinatorial criterion to be Kostant
modules, so this gives another proof of
the main result of \cite{CKL}.

\phantomsubsection{Endomorphism algebras of PIMs}
For any $\la \in \La$,
Theorem~\ref{main} implies that
the endomorphism algebra $\End_G(\cP(\la))^{\op}$ of
the projective indecomposable module $\cP(\la)$
is isomorphic to the algebra $e_\la K e_\la$.
By the definition of multiplication in $K$, this algebra is isomorphic to
$\C[x_1,\dots,x_r] / (x_1^2,\dots,x_r^2)$ where $r$ is the defect (atypicality) of
the block containing $\la$, answering a question raised recently by
several authors; see \cite[(4.2)]{BKN} and \cite[Conjecture 4.3.3]{Dr}.
(It should also be possible to give a proof of the commutativity of these endomorphism algebras
using some deformation theory like in \cite[$\S$2.8]{Scomp}, invoking the fact that the multiplicities $(\cP(\la):\cV(\mu))$ are at most one
by Theorem \ref{form1} below; see \cite[Theorem 7.1]{Squiv} for a similar situation.)

\phantomsubsection{Super duality}
When combined with the results from \cite{BS3},
our results can be used to prove the
``Super Duality Conjecture'' as formulated in \cite{CWZ}.
A direct algebraic proof of this
conjecture, and its substantial generalisation from \cite{CW},
has recently been found by Cheng and Lam \cite{CL}.
All of these results suggest some more direct geometric connection between the representation theory of $GL(m|n)$ and the
category of perverse sheaves on Grassmannians may exist.

\vspace{2mm}

To conclude this introduction, we recall in more detail the definition
of the algebra $K$ following \cite[$\S$6]{BS1}.
We assume $m=n=r$ and
focus just on the principal block of $G= GL(r|r)$, which is the
basic example of a block of atypicality $r$.
The dominant weights in this block are the weights
\begin{equation*}
\la_{i_1,\dots,i_r} := \sum_{s=1}^r (i_s+s-1) (\eps_s-\eps_{2r+1-s}) \in X^+(T)
\end{equation*}
parametrised by
sequences $i_1 > \cdots > i_r$ of integers.
According to the weight dictionary (\ref{wtdict}), the diagram
for $\la_{i_1,\dots,i_r}$
has label $\down$ at the vertices indexed by $i_1,\dots,i_r$, and label
$\up$ at all remaining vertices.
The corresponding block of the algebra $K$
is exactly the algebra denoted $K^\infty_r$ in the introduction of
\cite{BS1}. Theorem~\ref{main} (or rather the more precise
Lemmas~\ref{l1}--\ref{l2} below)
asserts in this situation
that
\begin{equation}\label{pb}
K^\infty_r \cong \End^{fin}_{G}\left(\bigoplus_{i_1 > \cdots > i_r}
  \!\!\cP(\la_{i_1,\dots,i_r})\right)^{\op}.
%= \!\!\!\bigoplus_{\substack{j_1 > \cdots > j_r \\ k_1 > \cdots > k_r}}
%\!\!\!\hom_G(\cP(\la_{j_1,\dots,j_r}), \cP(\la_{k_1,\dots,k_r})).
\end{equation}

Our explicit basis of $K^\infty_r$ is given by the {\em oriented circle diagrams} from (\ref{hbase}). These are
obtained by
taking the diagram of some
$\la_{i_1,\dots,i_r}$, then
gluing $r$ cups
and infinitely many rays to the bottom
and
$r$ caps and infinitely many rays to the top of the diagram so that
\begin{itemize}
\item
every vertex meets exactly one cup or ray below and exactly one cap or ray
above the number line;
\item each cup and each cap is incident with one vertex labelled
$\up$ and one vertex labelled $\down$;
\item
each ray
is incident with a
vertex labelled $\up$ and extends from there vertically up or down to infinity;
\item
no crossings of cups, caps and rays are allowed.
\end{itemize}
Under the isomorphism (\ref{pb}), such a
diagram
represents a homomorphism
$\cP(\la_{j_1,\dots,j_r}) \rightarrow \cP(\la_{k_1,\dots,k_r})$ where
$j_1 > \cdots > j_r$ (resp.\ $k_1 > \cdots > k_r$)
index the leftmost vertices of the
cups (resp.\ caps) in the diagram.
For example, here are
two oriented circle diagrams
corresponding to basis vectors in $K^\infty_2$
(where $\cdots$ indicates infinitely many
pairs of vertical rays labelled $\up$):
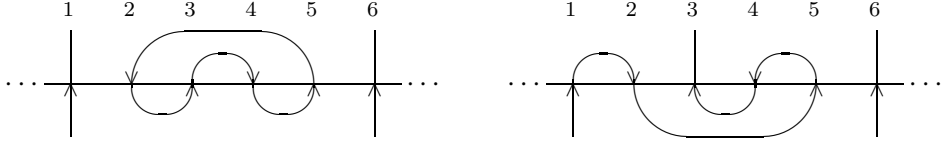
\begin{equation*}
\begin{picture}(310,55)
\put(-10,20){\line(1,0){135}}
\put(180,20){\line(1,0){135}}
\put(-25,17.2){$\cdots$}
\put(127,17.2){$\cdots$}
\put(165,17.2){$\cdots$}
\put(317,17.2){$\cdots$}

\put(-3,47.2){$_1$}
\put(20,47.2){$_2$}
\put(43,47.2){$_3$}
\put(66,47.2){$_4$}
\put(89,47.2){$_5$}
\put(112,47.2){$_6$}

\put(187,47.2){$_1$}
\put(210,47.2){$_2$}
\put(233,47.2){$_3$}
\put(256,47.2){$_4$}
\put(279,47.2){$_5$}
\put(302,47.2){$_6$}

\put(-2.9,15.5){$\scriptstyle\up$}
\put(20.1,20.1){$\scriptstyle\down$}
\put(43.1,15.5){$\scriptstyle\up$}
\put(66.1,20.1){$\scriptstyle\down$}
\put(89.1,15.5){$\scriptstyle\up$}
\put(112.1,15.5){$\scriptstyle\up$}

\put(57.5,20){\oval(23,23)[t]}
\put(57.5,20){\oval(69,40)[t]}
\put(80.5,20){\oval(23,23)[b]}
\put(34.5,20){\oval(23,23)[b]}
\put(115,0){\line(0,1){40}}
\put(0,0){\line(0,1){40}}

\put(187.1,15.5){$\scriptstyle\up$}
\put(210.1,20.1){$\scriptstyle\down$}
\put(233.1,15.5){$\scriptstyle\up$}
\put(256.1,20.1){$\scriptstyle\down$}
\put(279.1,15.5){$\scriptstyle\up$}
\put(302.1,15.5){$\scriptstyle\up$}

\put(247.5,20){\oval(23,23)[b]}
\put(247.5,20){\oval(69,40)[b]}
\put(190,0){\line(0,1){20}}
\put(236,20){\line(0,1){20}}
\put(270.5,20){\oval(23,23)[t]}
\put(201.5,20){\oval(23,23)[t]}
\put(305,0){\line(0,1){40}}
\end{picture}
\end{equation*}
The first diagram here represents a homomorphism
$\cP(\la_{4,2}) \rightarrow \cP(\la_{3,2})$ and the second one
represents a homomorphism
$\cP(\la_{3,2}) \rightarrow \cP(\la_{4,1})$.
Multiplying these two basis vectors together
as described in the next paragraph,
bearing in mind the ``$\op$'' in (\ref{pb}) which means for once that we are
writing maps on the right,
one gets the basis vector
$$
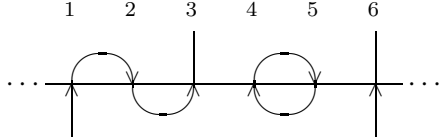
\begin{picture}(120,50)
\put(-10,20){\line(1,0){135}}
\put(-25,17.2){$\cdots$}
\put(127,17.2){$\cdots$}

\put(-3,47.2){$_1$}
\put(20,47.2){$_2$}
\put(43,47.2){$_3$}
\put(66,47.2){$_4$}
\put(89,47.2){$_5$}
\put(112,47.2){$_6$}

\put(-2.9,15.5){$\scriptstyle\up$}
\put(20.1,20.1){$\scriptstyle\down$}
\put(43.1,15.5){$\scriptstyle\up$}
\put(66.1,15.5){$\scriptstyle\up$}
\put(89.1,20.1){$\scriptstyle\down$}
\put(112.1,15.5){$\scriptstyle\up$}

\put(80.5,20){\oval(23,23)[b]}
\put(34.5,20){\oval(23,23)[b]}
\put(0,0){\line(0,1){20}}

\put(46,20){\line(0,1){20}}
\put(80.5,20){\oval(23,23)[t]}
\put(11.5,20){\oval(23,23)[t]}
\put(115,0){\line(0,1){40}}
\end{picture}
$$
which is some
homomorphism $\cP(\la_{4,2}) \rightarrow \cP(\la_{4,1})$.

We now sketch briefly the combinatorial procedure for multiplying
basis vectors.
Given two basis vectors, their product is necessarily zero unless the caps 
at the top of the first diagram are in exactly the same positions as the cups
at the bottom of the second.
Assuming that is the case, we glue the first diagram
underneath the second and join matching pairs of rays.
Then we perform a sequence of {\em generalised surgery procedures} to smooth
out all cup-cap pairs in the symmetric middle section of the resulting
composite diagram,
obtaining zero or more new diagrams in which
the middle section only involves vertical line segments. Finally we collapse
these middle sections to obtain a sum of basis vectors, which is the
desired  product.
Each generalised surgery procedure in this algorithm either
involves two components in the diagram merging into one or one
component splitting into two.
The rules for relabelling the new component(s) produced when this
operation is performed are summarized as follows:
\begin{align*}
1 \otimes 1 \mapsto 1,&
&&1 \otimes x \mapsto x,
&&x \otimes x \mapsto 0,\\
1 \otimes y \mapsto y,&
&&x \otimes y \mapsto 0,
&&y \otimes y \mapsto 0,\\
1 \mapsto 1 \otimes x + x \otimes 1,&
&&x \mapsto x \otimes x,
&&y \mapsto x \otimes y,
\end{align*}
where $1$ represents an anti-clockwise circle, $x$ represents a clockwise
circle, and $y$ represents a line.
This is a little cryptic;
we refer the reader to \cite{BS1} for a fuller account (and
explanation of
the connection to Khovanov's original construction via a certain TQFT).
Let us at least
apply this algorithm to the example from the previous
paragraph: two surgeries are needed,
the first of which involves an anti-clockwise circle and a line
merging together ($1 \otimes y
\mapsto y$)  and the second of which
involves a line splitting into a
clockwise circle and a line ($y \mapsto x \otimes y$):
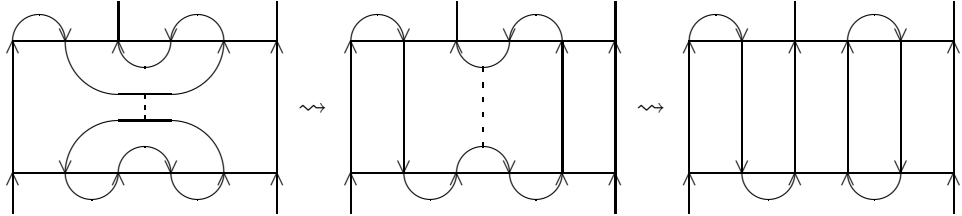
$$
\begin{picture}(138,80)
\dashline{2}(60,35)(60,45)
\put(118,37){$\rightsquigarrow$}
\put(246,37){$\rightsquigarrow$}
\put(10,15){\line(1,0){100}}
\put(10,65){\line(1,0){100}}
\put(7.1,10.5){$\scriptstyle\up$}
\put(27.1,15.1){$\scriptstyle\down$}
\put(47.1,10.5){$\scriptstyle\up$}
\put(67.1,15.1){$\scriptstyle\down$}
\put(87.1,10.5){$\scriptstyle\up$}
\put(107.1,10.5){$\scriptstyle\up$}
\put(60,15){\oval(20,20)[t]}
\put(60,15){\oval(60,40)[t]}
\put(80,15){\oval(20,20)[b]}
\put(40,15){\oval(20,20)[b]}
\put(110,0){\line(0,1){80}}
\put(10,0){\line(0,1){65}}

\put(7.1,60.5){$\scriptstyle\up$}
\put(27.1,65.1){$\scriptstyle\down$}
\put(47.1,60.5){$\scriptstyle\up$}
\put(67.1,65.1){$\scriptstyle\down$}
\put(87.1,60.5){$\scriptstyle\up$}
\put(107.1,60.5){$\scriptstyle\up$}

\put(60,65){\oval(20,20)[b]}
\put(60,65){\oval(60,40)[b]}
\put(50,65){\line(0,1){15}}
\put(80,65){\oval(20,20)[t]}
\put(20,65){\oval(20,20)[t]}
\end{picture}
\begin{picture}(128,80)
\put(0,15){\line(1,0){100}}
\put(0,65){\line(1,0){100}}
\put(-2.9,10.5){$\scriptstyle\up$}
\put(17.1,15.1){$\scriptstyle\down$}
\put(37.1,10.5){$\scriptstyle\up$}
\put(57.1,15.1){$\scriptstyle\down$}
\put(77.1,10.5){$\scriptstyle\up$}
\put(97.1,10.5){$\scriptstyle\up$}
\put(50,15){\oval(20,20)[t]}
\put(70,15){\oval(20,20)[b]}
\put(30,15){\oval(20,20)[b]}
\put(100,0){\line(0,1){80}}
\put(0,0){\line(0,1){65}}

\put(-2.9,60.5){$\scriptstyle\up$}
\put(17.1,65.1){$\scriptstyle\down$}
\put(37.1,60.5){$\scriptstyle\up$}
\put(57.1,65.1){$\scriptstyle\down$}
\put(77.1,60.5){$\scriptstyle\up$}
\put(97.1,60.5){$\scriptstyle\up$}

\put(50,65){\oval(20,20)[b]}
\put(40,65){\line(0,1){15}}
\put(70,65){\oval(20,20)[t]}
\put(10,65){\oval(20,20)[t]}
\put(20,65){\line(0,-1){50}}
\put(80,65){\line(0,-1){50}}
\dashline{2}(50,25)(50,55)
\end{picture}
\begin{picture}(110,80)
\put(20,65){\line(0,-1){50}}
\put(80,65){\line(0,-1){50}}
\put(40,65){\line(0,-1){50}}
\put(60,65){\line(0,-1){50}}
\put(0,15){\line(1,0){100}}
\put(0,65){\line(1,0){100}}
\put(-2.9,10.5){$\scriptstyle\up$}
\put(17.1,15.1){$\scriptstyle\down$}
\put(37.1,10.5){$\scriptstyle\up$}
\put(57.1,10.5){$\scriptstyle\up$}
\put(77.1,15.1){$\scriptstyle\down$}
\put(97.1,10.5){$\scriptstyle\up$}
\put(70,15){\oval(20,20)[b]}
\put(30,15){\oval(20,20)[b]}
\put(100,0){\line(0,1){80}}
\put(0,0){\line(0,1){65}}

\put(-2.9,60.5){$\scriptstyle\up$}
\put(17.1,65.1){$\scriptstyle\down$}
\put(37.1,60.5){$\scriptstyle\up$}
\put(57.1,60.5){$\scriptstyle\up$}
\put(77.1,65.1){$\scriptstyle\down$}
\put(97.1,60.5){$\scriptstyle\up$}

\put(40,65){\line(0,1){15}}
\put(70,65){\oval(20,20)[t]}
\put(10,65){\oval(20,20)[t]}
\end{picture}
$$
Contracting the middle section of the diagram on the right hand side here gives the final product recorded
already at the end of the previous paragraph.

The case $r=1$ in the above discussion (the principal block of $GL(1|1)$) is
easy to derive from scratch,
but still this is quite instructive.
So now the irreducible modules
are indexed simply by the weights $\{\la_i\:|\:i \in \Z\}$.
It is well known that $\cP(\la_i)$ has irreducible socle and head
isomorphic 
to $\cL(\la_i)$, with $\operatorname{rad} \cP(\la_i) / \soc \cP(\la_i) \cong
\cL(\la_{i-1}) \oplus \cL(\la_{i+1})$.
Hence in this case the locally finite endomorphism algebra
from (\ref{pb})
has basis $\{e_i, c_i, a_i, b_i\:|\:i \in \Z\}$,
where $e_i$ is the projection onto $\cP(\la_i)$,
$a_i:\cP(\la_{i}) \rightarrow \cP(\la_{i+1})$ and $b_i:\cP(\la_{i+1}) \rightarrow \cP(\la_{i})$ are non-zero homomorphisms
chosen so that
$b_i \circ a_i = a_{i-1} \circ b_{i-1}$,
and $c_i := b_i \circ a_i$ sends the head of $\cP(\la_i)$ onto its socle.
This corresponds to our diagram basis for $K^\infty_1$ 
so that
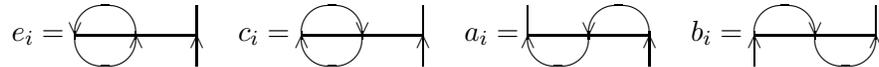
\begin{center}
\begin{picture}(72,28)
\put(-35,12){$e_i=$}
\put(9.2,9.4){$\scriptstyle\up$}
\put(32.2,9.4){$\scriptstyle\up$}
\put(-13.8,14.1){$\scriptstyle\down$}
\put(.5,14){\oval(23,23)[b]}
\put(.5,14){\oval(23,23)[t]}
\put(35,2.5){\line(0,1){23}}
\put(35,14){\line(-1,0){46}}
\end{picture}
\begin{picture}(62,28)
\put(-25,12){$c_i=$}
\put(19.2,14.1){$\scriptstyle\down$}
\put(-3.8,9.4){$\scriptstyle\up$}
\put(42.2,9.4){$\scriptstyle\up$}
\put(10.5,14){\oval(23,23)[b]}
\put(10.5,14){\oval(23,23)[t]}
\put(45,2.5){\line(0,1){23}}
\put(45,14){\line(-1,0){46}}
\end{picture}
\begin{picture}(62,28)
\put(-5,12){$a_i=$}
\put(39.2,14.1){$\scriptstyle\down$}
\put(62.2,9.4){$\scriptstyle\up$}
\put(16.2,9.4){$\scriptstyle\up$}
\put(53.5,14){\oval(23,23)[t]}
\put(30.5,14){\oval(23,23)[b]}
\put(19,14){\line(0,1){11.5}}
\put(65,14){\line(0,-1){11.5}}
\put(65,14){\line(-1,0){46}}
\end{picture}
\begin{picture}(62,28)
\put(15,12){$b_i=$}
\put(59.2,14.1){$\scriptstyle\down$}
\put(82.2,9.4){$\scriptstyle\up$}
\put(36.2,9.4){$\scriptstyle\up$}
\put(73.5,14){\oval(23,23)[b]}
\put(50.5,14){\oval(23,23)[t]}
\put(39,14){\line(0,-1){11.5}}
\put(85,14){\line(0,1){11.5}}
\put(85,14){\line(-1,0){46}}
\end{picture}
\end{center}
where we display only vertices $i, i+1, i+2$ 
and there are infinitely many pairs of vertical rays labelled
$\up$ at all other vertices.
In fact, $K^\infty_1$ is simply the
path algebra of the infinite quiver
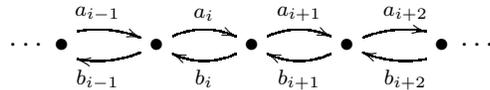
\begin{displaymath}
\xymatrix{
\cdots\:\bullet
\ar@/^/[r]^{a_{i-1}}&\ar@/^/[l]^{b_{i-1}}\bullet
\ar@/^/[r]^{a_{i}}&\ar@/^/[l]^{b_{i}}\bullet
\ar@/^/[r]^{a_{i+1}}&\ar@/^/[l]^{b_{i+1}}\bullet
\ar@/^/[r]^{a_{i+2}}&\ar@/^/[l]^{b_{i+2}}\bullet
\:\cdots
}
\end{displaymath}
modulo the relations $a_i b_i = b_{i-1} a_{i-1}$ 
and $a_{i} a_{i+1} = 0=b_{i+1} b_{i}$ for all
$i\in\mZ$.
It is clear from the quiver description that 
$K^\infty_1$ is naturally graded by path length; this is actually a Koszul grading. For general $r$ the canonical Koszul grading on
$K^\infty_r$ is defined by declaring that a basis vector is of degree
equal to the total number of clockwise cups and caps in the oriented
circle diagram.

\vspace{2mm}
\noindent
{\em Acknowledgements.}
This article was written up during stays by both authors at
the Isaac Newton Institute in Spring 2009.
We thank the INI staff and the Algebraic Lie Theory programme
organisers for the opportunity.

\section{Combinatorics of Grothendieck groups}

In this preliminary section, we compare the combinatorics
underlying the representation theory of $GL(m|n)$ with that of
the diagram algebra
$K(m|n)$.
Our exposition is largely
independent of \cite{B}, indeed, we will reprove the relevant results from there as we go.
On the other hand, we do assume that the reader is
familiar with the general theory of diagram algebras developed in \cite{BS1, BS2}.
Later in the article we will also need to appeal to various results from \cite{BS3}.

\phantomsubsection{\boldmath Representation theory of $K(m|n)$}
Fix once and for all integers $m,n \geq 0$.
Let $K = K(m|n)$ and $\La = \La(m|n)$ be as in the introduction.
The elements $\{e_\la\:|\:\la\in \La\}$ form a system of (in general infinitely many)
mutually orthogonal idempotents in $K$
such that
\begin{equation}\label{locun}
K = \bigoplus_{\la,\mu\in \La} e_\la K e_\mu.
\end{equation}
So the algebra $K$ is {\em locally unital}, but it
is not unital (except in the trivial case $m=n=0$).
By a $K$-module we always mean a locally unital module; for
a left $K$-module $M$ this means that $M$
decomposes as $$
M = \bigoplus_{\la \in \La}
e_\la M.
$$
The irreducible $K$-modules $\{L(\la)\:|\:\la \in \La\}$ defined in
the introduction
are all one dimensional, so $K$ is a basic algebra.

Let $\Rep{K}$ denote the category of
finite dimensional left $K$-modules.
The Grothendieck group $[\Rep{K}]$
of this category is the free $\Z$-module
on basis
$\{[L(\la)]\:|\:\la \in \La\}$.
The standard modules
$\{V(\la)\:|\:\la \in \La\}$
and the projective indecomposable
modules
$\{P(\la)\:|\:\la \in \La\}$
from \cite[$\S$5]{BS1}
are
finite dimensional, so it makes sense to
consider their classes $[V(\la)]$ and $[P(\la)]$
in $[\Rep{K}]$.
Finally, we use the notation $\mu \supset \la$
(resp.\ $\mu \subset \la$) from \cite[$\S$2]{BS1}
to indicate that the composite diagram
$\mu \overline{\la}$ (resp.\ $\underline{\mu} \la$)
is oriented in the obvious sense.

\begin{Theorem}\label{form1}
We have
in $[\Rep{K}]$ that
\begin{equation*}
[P(\la)] = \sum_{\mu \supset \la}
[V(\mu)],\qquad
[V(\la)] = \sum_{\mu \subset \la}
[L(\mu)]
\end{equation*}
for each $\la \in \La$.
\end{Theorem}

\begin{proof}
This follows from \cite[Theorem 5.1]{BS1} and \cite[Theorem 5.2]{BS1}.
\end{proof}

As $\mu
\supset \la$
(resp. $\mu \subset \la$) implies that $\mu \geq \la$
(resp. $\mu \leq \la$) in the Bruhat ordering,
we deduce from Theorem~\ref{form1} that the
classes $\{[P(\la)]\}$ and $\{[V(\la)]\}$
are linearly independent in $[\Rep{K}]$.
However they do not span $[\Rep{K}]$ as the
chains in the Bruhat order are infinite.

\begin{Remark}\label{grep}\rm
The algebra $K$ possesses a natural $\Z$-grading
defined by declaring that each basis vector
$(a \la b)$ from (\ref{hbase}) is of degree equal to
the number
of clockwise cups and caps
in the diagram $a \la b$.
This means that one can consider the
{\em graded} representation theory of $K$.
The various modules
$L(\la), V(\la)$ and $P(\la)$ also possess canonical
gradings, as is discussed in detail in \cite[$\S$5]{BS1}.
\end{Remark}

\phantomsubsection{Special projective functors: the diagram side}
As in \cite[(2.5)]{BS3},
let us represent a block $\Ga\in\La / \sim$
by means of its {\em block diagram}, that is, the diagram obtained by taking
any $\la \in \Ga$ and replacing all the $\up$'s and $\down$'s
labelling its vertices by the symbol $\bullet$.
Because $m$ and $n$ are fixed, the block $\Ga$ can be recovered
uniquely from its block diagram.
Recall also the notion of the {\em defect} of a weight $\la \in \La$
from \cite[$\S$2]{BS1}. In this setting, this simply means the number of
vertices labelled $\down$ in $\la$, and the defect of $\la$ is the
same thing as the defect $\defect(\Ga)$ of the unique block $\Ga \in
\La / \sim$ containing $\la$.

Given a block $\Ga$,
we say that $i \in \Z$ is {\em $\Ga$-admissible} if
the $i$th and $(i+1)$th vertices of the block diagram of $\Ga$
match
the top number line
of a unique one of the following pictures, and $\defect(\Ga)$ is as indicated:
\begin{equation}
% scriptstyles at (1,x) x = -68,-45,-22,1,24,47,70,93,...
% \times is -0.2
% \circ is +1
% \up is +0.3 and down 2.6
% \down is +0.3 and up 2.1
\begin{picture}(75,33)
\put(159.5,12){$\Ga$}
\put(155,1){$t_i(\Ga)$}
\put(152,-11.4){$\Ga-\alpha_i$}
\put(-106,17){\vector(0,-1){28}}\put(-125,-.5){$F_i$}
\put(-86,28){$_\text{$\defect(\Ga) \geq 1$}$}
\put(-26,28){$_\text{$\defect(\Ga) \geq 0$}$}
\put(34,28){$_\text{$\defect(\Ga) \geq 0$}$}
\put(94,28){$_\text{$\defect(\Ga) \geq 0$}$}
%\put(-81,28){$\text{$_i$\quad\:\,$_{i+1}$}$}
%\put(-21,28){$\text{$_i$\quad\:\,$_{i+1}$}$}
%\put(36,28){$\text{$_i$\quad\:\,$_{i+1}$}$}
%\put(99.5,28){$\text{$_i$\quad\:\,$_{i+1}$}$}
\put(-82,-24){$\text{``cup''}$}
\put(-21,-24){$\text{``cap''}$}
\put(21,-24){$\text{``right-shift''}$}
\put(89,-24){$\text{``left-shift''}$}

\put(-84,15){\line(1,0){33}}
\put(-84,-9){\line(1,0){33}}
\put(-67.5,15){\oval(23,23)[b]}
\put(-59.3,-10.9){$\scriptstyle\times$}
\put(-81.8,-11.6){$\circ$}
\put(-80.9,13.1){{$\scriptstyle\bullet$}}
\put(-57.9,13.1){{$\scriptstyle\bullet$}}

\put(-24,15){\line(1,0){33}}
\put(-24,-9){\line(1,0){33}}
\put(-7.5,-9){\oval(23,23)[t]}
\put(1.2,12.4){$\circ$}
\put(-22.2,13.1){$\scriptstyle\times$}
\put(-21.1,-10.9){{$\scriptstyle\bullet$}}
\put(1.9,-10.9){{$\scriptstyle\bullet$}}

\put(36,15){\line(1,0){33}}
\put(36,-9){\line(1,0){33}}
\put(61.2,12.4){$\circ$}
\put(38.2,-11.6){$\circ$}
\put(64,-8.2){\line(-1,1){22.9}}
\put(61.9,-10.9){{$\scriptstyle\bullet$}}
\put(38.9,13.1){{$\scriptstyle\bullet$}}

\put(96,15){\line(1,0){33}}
\put(96,-9){\line(1,0){33}}
\put(120.7,-10.9){$\scriptstyle\times$}
\put(97.7,13.1){$\scriptstyle\times$}
\put(101.1,-8.2){\line(1,1){22.9}}
\put(98.9,-10.9){{$\scriptstyle\bullet$}}
\put(121.9,13.1){{$\scriptstyle\bullet$}}
\end{picture}\label{CKLR}
\end{equation}
\vspace{4mm}

\noindent
Assuming $i$ is $\Ga$-admissible, we let $(\Ga - \alpha_i)$ denote the block
obtained from $\Ga$
by relabelling the $i$th and $(i+1)$th vertices of its block diagram
according to the bottom
number line of the appropriate picture.
Also define a $(\Ga-\alpha_i)\Ga$-matching $t_i(\Ga)$
in the sense of \cite[$\S$2]{BS2}
so that the strip between the $i$th and $(i+1)$th vertices of $t_i(\Ga)$
is as in the picture, and there are only vertical ``identity'' line segments
elsewhere.

For blocks $\Ga,\De \in \La / \sim$ and a $\Ga\De$-matching $t$,
recall the geometric bimodule $K^t_{\Ga\De}$\label{cmat}
from \cite[$\S$3]{BS2}. By definition this is a
$(K_\Ga,K_\De)$-bimodule. We can view it as a
$(K,K)$-bimodule by extending the actions of $K_\Ga$ and $K_\De$ to all
of $K$ so that the other blocks act as zero.
The functor $K^t_{\Ga\De}\otimes_K ?$ is an
endofunctor of $\Rep{K}$ called a {\em projective
  functor}. Writing $t^*$ for the mirror image of $t$ in a horizontal
axis, the functor $K^{t^*}_{\De\Ga} \otimes_K ?$ gives another projective
functor which is biadjoint to
$K^t_{\Ga\De}\otimes_K ?$
by \cite[Corollary 4.9]{BS2}.

For any $i \in \Z$, introduce the $(K,K)$-bimodules
\begin{equation}\label{bims}
\widetilde{F}_i := \bigoplus_{\Ga} K^{t_i(\Ga)}_{(\Ga-\alpha_i)\Ga},
\qquad
\widetilde{E}_i := \bigoplus_{\Ga} K^{t_i(\Ga)^*}_{\Ga(\Ga-\alpha_i)},
\end{equation}
where the direct sums are over all
$\Ga \in \La / \sim$ such that $i$ is $\Ga$-admissible.
The {\em special projective functors} are the
endofunctors $F_i := \widetilde{F}_i \otimes_K ?$ and $E_i :=
\widetilde{E}_i \otimes_K ?$ of $\Rep{K}$ defined by tensoring
with these bimodules.
The discussion in the previous paragraph implies that the functors $F_i$ and $E_i$
are biadjoint, hence
they are both exact and map projectives to projectives.

For $\la \in \La$, let $I_\times(\la) := I_\down(\la) \cap I_\up(\la)$
(resp.\ $I_\circ(\la) := \Z \setminus (I_\down(\la) \cup I_\up(\la))$) denote the set of
integers indexing the vertices labelled $\times$
(resp.\ $\circ$)
in $\la$; cf. (\ref{wtdict}).
Introduce the notion of the
{\em height} of $\la$:
\begin{align}\label{absdeg}
\deg(\la) &:= \sum_{i \in I_\times(\la)} i - \sum_{i \in I_\circ(\la)} i.
\end{align}
Note all weights belonging to the same block have the same
height.

\begin{Lemma}\label{degrees}
For $\la \in \La$ and $i \in \Z$,
all composition factors of $F_i L(\la)$
(resp.\ $E_i L(\la)$)
are of the form $L(\mu)$
with $\deg(\mu) = \deg(\la)+1$ (resp.\ $\deg(\la)-1$).
\end{Lemma}

\begin{proof}
This follows by inspecting (\ref{CKLR}).
\end{proof}

\begin{Lemma}\label{ga}
Let $\la \in \La$ and $i \in \Z$.
For symbols $x,y \in \{\circ,\up,\down,\times\}$
we write $\la_{xy}$ for the diagram obtained from $\la$
by relabelling the $i$th and $(i+1)$th vertices
by $x$ and $y$, respectively.
\begin{itemize}
\item[\rm(i)]
If $\la = \la_{{\scriptscriptstyle\down}\circ}$ then
$F_i P(\la) \cong P(\la_{\circ{\scriptscriptstyle\down}})$,
$F_i V(\la) \cong V(\la_{\circ{\scriptscriptstyle\down}})$,
$F_i L(\la) \cong L(\la_{\circ{\scriptscriptstyle\down}})$.
\item[\rm(ii)]
If $\la = \la_{{\scriptscriptstyle\up}\circ}$ then
$F_i P(\la) \cong P(\la_{\circ{\scriptscriptstyle\up}})$,
$F_i V(\la) \cong V(\la_{\circ{\scriptscriptstyle\up}})$,
$F_i L(\la) \cong L(\la_{\circ{\scriptscriptstyle\up}})$.
\item[\rm(iii)]
If $\la = \la_{\scriptscriptstyle\times\down}$ then
$F_i P(\la) \cong P(\la_{\scriptscriptstyle\down\times})$,
$\!F_i V(\la) \cong V(\la_{\scriptscriptstyle\down\times})$,
$\!F_i L(\la) \cong L(\la_{\scriptscriptstyle\down\times})$.
\item[\rm(iv)]
If $\la = \la_{\scriptscriptstyle\times\up}$ then
$F_i P(\la) \cong P(\la_{\scriptscriptstyle\up\times})$,
$\!F_i V(\la) \cong V(\la_{\scriptscriptstyle\up\times})$,
$\!F_i L(\la) \cong L(\la_{\scriptscriptstyle\up\times})$.
\item[\rm(v)]
If $\la = \la_{{\scriptscriptstyle\times}\circ}$ then:
\begin{itemize}
\item[(a)]
$F_i P(\la) \cong P(\la_{\scriptscriptstyle\down\up})$;
\item[(b)]
there is a short exact sequence
$$
0 \rightarrow V(\la_{\scriptscriptstyle\up\down}) \rightarrow F_i V(\la)
\rightarrow V(\la_{\scriptscriptstyle\down\up}) \rightarrow 0;
$$
\item[(c)]
$F_i L(\la)$ has irreducible
socle and head both isomorphic to $L(\la_{\scriptscriptstyle\down\up})$, and all other composition
factors are of the form $L(\mu)$ for $\mu \in\La$
such that
$\mu = \mu_{\scriptscriptstyle\down\down}$,
$\mu = \mu_{\scriptscriptstyle\up\up}$ or
$\mu = \mu_{\scriptscriptstyle\up\down}$.
\end{itemize}
\item[\rm(vi)]
If $\la = \la_{{\scriptscriptstyle\down\up}}$ then
$F_i P(\la) \cong
P(\la_{\circ\scriptscriptstyle{\times}}) \oplus
P(\la_{\circ\scriptscriptstyle{\times}})$,
$F_i V(\la) \cong V(\la_{\circ{\scriptscriptstyle\times}})$
and
$F_i L(\la) \cong L(\la_{\circ{\scriptscriptstyle\times}})$.
\item[\rm(vii)]
If $\la = \la_{{\scriptscriptstyle\up\down}}$ then
$F_i V(\la) \cong V(\la_{\circ{\scriptscriptstyle\times}})$
and $F_i L(\la) = \{0\}$.
\item[\rm(viii)]
If
$\la = \la_{{\scriptscriptstyle\down\down}}$
then
$F_i V(\la) = F_i L(\la) = \{0\}$.
\item[\rm(ix)]
If
$\la = \la_{{\scriptscriptstyle\up\up}}$
then
$F_i V(\la) = F_i L(\la) = \{0\}$.
\item[\rm(x)]
For all other $\la$ we have that
$F_i P(\la) = F_i V(\la) = F_i L(\la) = \{0\}$.
\end{itemize}
For the dual statement about $E_i$,
interchange all occurrences of $\circ$ and $\times$.
\end{Lemma}

\begin{proof}
Apply \cite[Theorems 4.2]{BS2}, \cite[Theorem 4.5]{BS2} and \cite[Theorem 4.11]{BS2}, exactly as was done in \cite[Lemma 3.4]{BS3}.
\end{proof}

\begin{Remark}\rm\label{gradings}
Using Lemma~\ref{ga}, one can check that the
endomorphisms of $[\Rep{K}]$ induced by the functors $F_i$ and
$E_i$ for all $i \in \Z$ satisfy the Serre relations
defining the Lie algebra
$\mathfrak{sl}_{\infty}$.
Indeed, letting $V_\infty$ denote the natural
$\mathfrak{sl}_{\infty}$-module of column vectors,
the category $\Rep{K}$ can be interpreted in a precise sense as a
categorification of a certain completion of the
 $\mathfrak{sl}_{\infty}$-module
$\bigwedge^m V_\infty
\otimes \bigwedge^n V_\infty^*$;
see also \cite{B} where this point of view is taken on the supergroup side.
Using the graded representation theory mentioned in Remark~\ref{grep},
i.e. replacing $\Rep{K}$ with the category of finite dimensional
{\em graded} $K$-modules, one gets a categorification of the
$q$-analogue of this module over the
quantised enveloping algebra $U_q(\mathfrak{sl}_\infty)$; the action of $q$ comes from shifting the grading on a module up by one. We are not going to pursue this connection further here, but refer the reader to
\cite[Theorem 3.5]{BS3} where an analogous ``graded categorification theorem'' is discussed in detail.
\end{Remark}

\phantomsubsection{The crystal graph}
Define the {\em crystal graph} to be the directed coloured graph
with vertex set equal to $\La$ and a directed edge
$\mu \stackrel{i}{\rightarrow} \la$ of colour $i \in \Z$
if $L(\la)$ is a quotient of $F_i L(\mu)$.
It is clear from Lemma~\ref{ga}
that $\mu \stackrel{i}{\rightarrow} \la$ if and only if
the $i$th and $(i+1)$th vertices of $\la$ and $\mu$ are labelled according
to one of the six cases in the following table, and all other vertices
of $\la$ and $\mu$
are labelled in the same way:
\begin{equation}\label{cg}
\begin{array}{c|c|c|c|c|c|c}
\mu&\quad{\scriptstyle\down}\:\circ\quad&\quad{\scriptstyle\up}\:\circ\quad&\quad{\scriptstyle\times}\:{\scriptstyle\down}\quad&\quad{\scriptstyle\times}\:{\scriptstyle\up}\quad&\quad{\scriptstyle\times}\:\circ\quad&\quad{\scriptstyle\down}\:{\scriptstyle\up}\quad\\
\hline
\la&\circ\:{\scriptstyle\down}&\circ\:{\scriptstyle\up}&{\scriptstyle\down}\:{\scriptstyle\times}&{\scriptstyle\up}\:{\scriptstyle\times}&{\scriptstyle\down}\:{\scriptstyle\up}&\circ\:{\scriptstyle\times}
\end{array}
\end{equation}
Comparing this explicit description with \cite[$\S$3-d]{B}, it follows that
our crystal graph is isomorphic to Kashiwara's crystal graph
associated to the $\mathfrak{sl}_{\infty}$-module mentioned in Remark~\ref{gradings}, which hopefully
explains our choice of terminology.

Suppose we are given integers $p \leq q$.
Define the following intervals
\begin{align}\label{intervals}
I_{p,q} &:= \{p-m+1,p-m+2,\dots,q+n-1\},\\
I_{p,q}^+ &:= \{p-m+1,p-m+2,\dots,q+n-1,q+n\}.
\end{align}
(The reader may find it helpful at this point to note which vertices of the weight
$\la_{p,q}$ from (\ref{tuesday}) are indexed by the set $I_{p,q}^+$.)
Then introduce the following subsets of $\La$:
\begin{align}\label{lamn1}
\La_{p,q} &:= \{\la \in \La\:|\:
\text{the $i$th vertex of $\la$
is labelled $\up$ for all $i \notin I_{p,q}^+$}\},\\
\La_{p,q}^{\circ} &:= \Big\{\la \in \La_{p,q}\:\Big|\:
\begin{array}{l}
\text{amongst vertices $j,\dots,q+n$ of $\la$, the number}\\
\text{of $\up$'s is $\geq$ the number of $\down$'s, for all $j \in I_{p,q}^+$}\\
\end{array}
\Big\}.\label{lamn2}
\end{align}
Note that the weight $\la_{p,q}$ from (\ref{tuesday})
belongs to
$\La_{p,q}^\circ$. It
is the unique weight
in $\La_{p,q}$ of minimal
height.

\begin{Lemma}\label{cgt}
Given $\la \in \La$,
choose $p \leq q$ such that
$\la \in \La_{p,q}^\circ$ (which is always possible as there are infinitely many
$\up$'s and finitely many $\down$'s).
Then there are integers $i_1,\dots,i_d \in I_{p,q}$, where $d = \deg(\la) - \deg(\la_{p,q})$,
such that $\la_{p,q}
\stackrel{i_1}{\rightarrow} \cdots \stackrel{i_d}{\rightarrow}
\la$ is a path in the crystal graph.
Moreover we have that
$$
F_{i_d} \cdots F_{i_1} V(\la_{p,q}) \cong
P(\la)^{\oplus 2^r},
$$
where $r$ is the number of
edges in the given path
of the form ${\scriptstyle \down\up}
\rightarrow \circ{\scriptstyle\times}$.
\end{Lemma}

\begin{proof}
For the first statement, we proceed by induction on $\deg(\la)$.
If $\deg(\la) = \deg(\la_{p,q})$, then $\la = \la_{p,q}$
and the conclusion is trivial. Now assume that
$\deg(\la) > \deg(\la_{p,q})$.
As $\la \in \La_{p,q}^\circ$ and $\la \neq \la_{p,q}$,
it is possible to find $i \in I_{p,q}$ such that the $i$th
and $(i+1)$th vertices of $\la$ are labelled
$\circ \down, \circ \up, \circ \times,
\down\times, \up\times$ or $\down\up$.
Inspecting (\ref{cg}), there is a unique weight $\mu \in \La_{p,q}^\circ$
with $\mu \stackrel{i}{\rightarrow} \la$ in the crystal
graph. Noting $\deg(\mu) = \deg(\la)-1$, we are now done by induction.
To deduce the second statement,
we apply Lemma~\ref{ga} to get easily that
$F_{i_d} \cdots F_{i_1} P(\la_{p,q})
\cong P(\la)^{\oplus 2^r}$.
Finally $P(\la_{p,q})
\cong V(\la_{p,q})$
as $\la_{p,q}$ is of defect zero, by \cite[Theorem 5.1]{BS1}.
\end{proof}

\phantomsubsection{\boldmath Representation theory of $GL(m|n)$}
Now we turn to discussing the representation theory of $G = GL(m|n)$.
In the introduction, we defined already the abelian category
$\sF = \sF(m|n)$ and the irreducible
modules
$\{\cL(\la)\}$,\label{reps}
the standard modules $\{\cV(\la)\}$ and the projective indecomposable
modules
$\{\cP(\la)\}$, all of which are
parametrised by the set $X^+(T)$ of dominant weights.
We are using an unusual font here
(and a few other places later on) to avoid confusion with the analogous
$K$-modules
$\{L(\la)\}$, $\{V(\la)\}$ and $\{P(\la)\}$.
Recall in particular that the $\Z_2$-grading on $\cL(\la)$
is defined so that its $\la$-weight space is
concentrated in degree
$\bar\la
:= (\la,\eps_{m+1}+\dots+\eps_{m+n}) \pmod{2}$.
Bearing in mind that we consider only even morphisms,
the modules
$$
\{\cL(\la)\:|\:\la \in X^+(T)\}
\cup
\{\Pi \cL(\la)\:|\:\la \in X^+(T)\}
$$
give a complete set of pairwise non-isomorphic irreducible
$G$-modules, where $\Pi$ denotes the change of parity functor.

The standard module $\cV(\la)$ is usually called a {\em Kac module} in this setting
after \cite{Kac2}, and can be constructed explicitly as follows.
Let $P$ be the parabolic subgroup of $G$
such that $P(A)$
consists
of all invertible
matrices of the form (\ref{supermat}) with $c=0$,
for each commutative superalgebra $A$.
Given $\la \in X^+(T)$, we let $E(\la)$
denote the usual
finite dimensional irreducible module of highest weight $\la$ for the
underlying even subgroup $G_{\0} \cong GL(m) \times GL(n)$,
viewing $E(\la)$ as a supermodule with $\Z_2$-grading concentrated
in degree $\bar\la$.
We can regard $E(\la)$ also as a $P$-module
by inflating through the obvious homomorphism $P \twoheadrightarrow G_{\0}$.
Then we have that
\begin{equation}\label{kacmod}
\cV(\la) = U(\mathfrak{g}) \otimes_{U(\mathfrak{p})}
E(\la),
\end{equation}
where $\mg$ and $\mathfrak{p}$ denote the Lie superalgebras of
$G$ and $P$, respectively.
This construction makes sense because the induced module on the right hand side of (\ref{kacmod}) is an integrable $\mathfrak{g}$-supermodule, i.e. it lifts in a unique way to a $G$-module.

The module $\cL(\la)$ is isomorphic to the unique irreducible quotient of
$\cV(\la)$.
Also $\cP(\la)$ is the projective cover of $\cL(\la)$ in the category
$\sF$. It has a
{\em standard flag}, that is, a filtration whose sections are standard
modules. The multiplicity $(\cP(\la):\cV(\mu))$ of
$\cV(\mu)$ as a section of any such standard flag is given by
the {\em BGG reciprocity formula}
\begin{equation}\label{bggrec}
(\cP(\la):\cV(\mu)) = [\cV(\mu):\cL(\la)],
\end{equation}
as follows from
\cite{Z} or the discussion in \cite[Example 7.5]{Btilt}.

\phantomsubsection{Special projective functors: the supergroup side}
Recall the {\em weight dictionary} from (\ref{wtdict}) by means of which
we identify the set $X^+(T)$ with the set $\La = \La(m|n)$.
Under this identification, the usual
notion of the {degree of atypicality} of a weight $\la \in X^+(T)$
corresponds to the notion of {\em defect} of $\la \in \La$.
Given $\la, \mu \in \La$,
the irreducible $G$-modules $\cL(\la)$ and $\cL(\mu)$ have the
same central character if and only if
$\la \sim \mu$ in the diagrammatic sense;
this can be deduced from \cite[Corollary 1.9]{Serg}.
Hence the category $\sF$ decomposes as
\begin{equation}\label{blockclass}
\sF = \bigoplus_{\Ga \in \La / \sim} \sF_\Ga,
\end{equation}
where $\sF_\Ga$ is the full subcategory consisting of the
modules all of whose composition factors are of the form
$\cL(\la)$ for $\la \in \Ga$.
We let $\pr_\Ga:\sF \rightarrow \sF$ be the exact
functor defined by projection onto $\sF_\Ga$ along (\ref{blockclass}).

Recall that $\nV$ denotes the natural $G$-module
and $\nV^*$ is its dual.
Following \cite[(4.21)--(4.22)]{B}, we define the
{\em special projective functors}
$\cF_i$ and $\cE_i$ for each $i \in \Z$
to be the following endofunctors of $\sF$:
\begin{equation}\label{spf}
\cF_i := \bigoplus_{\Ga} \pr_{\Ga-\alpha_i} \circ (? \otimes\nV)
\circ \pr_\Ga,
\quad
\cE_i := \bigoplus_{\Ga} \pr_{\Ga} \circ (? \otimes\nV^*)
\circ \pr_{\Ga-\alpha_i},
\end{equation}
where the direct sums are over all $\Ga \in \La / \sim$
such that $i$ is $\Ga$-admissible
(as in (\ref{bims})).
The functors $? \otimes \nV$ and $? \otimes \nV^*$ are biadjoint,
hence so are
$\cF_i$ and $\cE_i$.
In particular, all these functors are exact
and send projectives to projectives.
For later use, let us
fix once and for all a choice of
an adjunction making $(\cF_i, \cE_i)$ into an adjoint pair for each $i \in \Z$.

\begin{Lemma}\label{onstdmods}
The following hold for any $\la \in X^+(T)$:
\begin{itemize}
\item[(i)] $\cV(\la) \otimes \nV$ has a filtration with
sections $\cV(\la + \eps_r)$
for all $r=1,\dots,m+n$ such that $\la+\eps_r \in X^+(T)$,
arranged in order from bottom to top.
\item[(ii)]
$\cV(\la) \otimes \nV^*$ has a filtration with
sections $\cV(\la - \eps_r)$
for all $r=1,\dots,m+n$ such that $\la-\eps_r \in X^+(T)$,
arranged in order from top to bottom.
\end{itemize}
\end{Lemma}

\iffalse
\begin{Lemma}\label{onstdmods}
Given $\la \in X^+(T)$, let $v_\la$ be a non-zero highest weight vector
in $\cV(\la)$ of weight $\la$.
\begin{itemize}
\item[(i)] List the integers
$\{1 \leq r \leq m+n\:|\:\la+\eps_r \in X^+(T)\}$
as $r_1 < \cdots < r_k$.
Then there is a filtration
$$
0 = M_0 < M_1 < \cdots < M_k = \cV(\la) \otimes V
$$
such that $M_j / M_{j-1} \cong \cV(\la+\eps_{r_j})$ for each
$j=1,\dots, k$. The submodule $M_j$ here is
generated by the vectors $v_\la \otimes v_{r_1},\dots,v_\la \otimes
v_{r_j}$, and the image of $v_\la \otimes
v_{r_j}$ in $M_j / M_{j-1}$ is a non-zero
highest weight vector.
\item[(ii)]
List the integers
$\{1 \leq s \leq m+n\:|\:\la-\eps_s \in X^+(T)\}$
as $s_1 > \cdots > s_l$.
Then there is a filtration
$$
0 = N_0 < N_1 < \cdots < N_k = \cV(\la) \otimes V^*
$$
such that $N_j / N_{j-1} \cong \cV(\la+\eps_{r_j})$ for each
$j=1,\dots, k$. Denoting the basis of $V^*$ dual to the standard basis
of $V$ by $f_1,\dots,f_{m+n}$,
the submodule $N_j$ here is generated by the vectors $v_\la \otimes f_{r_1},\dots,v_\la \otimes
f_{r_j}$, and the image of $v_\la \otimes
f_{s_j}$ in $N_j / N_{j-1}$ is a non-zero
highest weight vector.
\end{itemize}
\end{Lemma}
\fi

\begin{proof}
This follows
from the definition (\ref{kacmod}) and the tensor identity.
\end{proof}

\begin{Corollary}\label{onstdmodsc}
The following hold for any $\la \in X^+(T)$ and $i \in \Z$:
\begin{itemize}
\item[(i)]
$\cF_i \cV(\la)$ has a filtration with
sections $\cV(\la + \eps_r)$
for all $r=1,\dots,m+n$ such that $\la+\eps_r \in X^+(T)$
and $(\la+\rho,\eps_r) = i + (1-(-1)^{\bar r})/2$,
arranged in order from bottom to top.
\item[(ii)]
$\cE_i \cV(\la)$ has a filtration with
sections $\cV(\la - \eps_r)$
for all $r=1,\dots,m+n$ such that $\la-\eps_r \in X^+(T)$
and $(\la+\rho,\eps_r) = i+(1+(-1)^{\bar r})/2$,
arranged in order from top to bottom.
\end{itemize}
\end{Corollary}

\begin{proof}
For (i),
apply $\pr_{\Ga-\alpha_i}$ to
the statement of Lemma~\ref{onstdmods}(i),
where $\Ga$ is the block containing $\la$ (and do a little work to
translate the combinatorics).
The proof of (ii) is similar.
\end{proof}

\begin{Corollary}\label{fullt}
We have that
$? \otimes \nV = \bigoplus_{i \in \Z} \cF_i$
and
$? \otimes \nV^* = \bigoplus_{i \in \Z} \cE_i$.
\end{Corollary}

The next lemma gives an alternative definition of the functors
$\cF_i$ and $\cE_i$ which will be needed in the next section;
cf. \cite[Proposition 5.2]{CW}. Let
\begin{equation}\label{omega}
\Omega := \sum_{r,s=1}^{m+n} (-1)^{\bar s} e_{r,s} \otimes e_{s,r}
\in \mathfrak{g} \otimes \mathfrak{g},
\end{equation}
where $e_{r,s}$ denotes the $rs$-matrix unit.
This corresponds to the
supertrace form on $\mathfrak{g}$,
so left multiplication by $\Omega$ (interpreted with the usual superalgebra
sign conventions) defines a $G$-module endomorphism of $M \otimes N$
for any $G$-modules $M$ and $N$.

\begin{Lemma}\label{eig}
For any $G$-module $M$,
we have that $\cF_i M$
(resp.\ $\cE_i M$)
is the generalised $i$-eigenspace (resp. the generalised $-(m-n+i)$-eigenspace)
of the operator $\Omega$ acting on $M \otimes \nV$
(resp.\ $M \otimes \nV^*$).
\end{Lemma}

\begin{proof}
We just explain for $\cF_i$.
Let $c := \sum_{r,s=1}^{m+n} (-1)^{\bar s} e_{r,s} e_{s,r} \in U(\mathfrak{g})$ be the Casimir element.
It acts on $\cV(\la)$ by multiplication by the scalar
$$
c_\la := \left(\la+2\rho + (m-n-1) \delta, \la\right)
$$
where $\delta = \eps_1+\cdots+\eps_m-\eps_{m+1}-\cdots-\eps_{m+n}$.
Also,
we have that $\Omega = (\Delta(c) - c \otimes 1 - 1 \otimes c) / 2$
where $\Delta$ is the comultiplication of $U(\mathfrak{g})$.
Now to prove the lemma,
it suffices to verify it for the special case
$M = \cV(\la)$.
Using the
observations just made, we see that
multiplication by $\Omega$ preserves the filtration from
Lemma~\ref{onstdmods}(i), and
the induced action of $\Omega$
on the
section $\cV(\la+\eps_r)$ is
by multiplication by the scalar
$$
(c_{\la+\eps_r} - c_\la - m+n) / 2 =
(\la+\rho,\eps_r) + (1-(-1)^{\bar i})/2.
$$
The result follows on comparing with Corollary~\ref{onstdmodsc}(i).
\end{proof}

The next two lemmas are
the key to understanding the
representation theory of $GL(m|n)$ from a combinatorial point of view.

\begin{Lemma}\label{jsf}
Let $i \in \Z$ and $\la
\in \La$ be a weight such that the
$i$th and $(i+1)$th vertices of $\la$ are labelled
$\up$ and $\down$, respectively.
Let $\mu$ be the weight obtained from $\la$
by interchanging the labels on these two vertices.
Then $\cL(\mu)$ is
a composition factor of
$\cV(\la)$.
\end{Lemma}

\begin{proof}
This is a reformulation of
\cite[Theorem 5.5]{Serg}.
It can be proved directly
by an explicit calculation with
certain lowering operators in $U(\mathfrak{g})$ as in
\cite[Lemma 4.8]{BS3}.
\end{proof}

\begin{Lemma}\label{ga2}
Exactly the same statement as Lemma~\ref{ga} holds in the category
$\sF$,
replacing $L(\la), V(\la), P(\la), F_i$ and $E_i$ by $\cL(\la),
\cV(\la)$, $\cP(\la)$, $\cF_i$ and $\cE_i$, respectively.
\end{Lemma}

\begin{proof}
The statements involving $\cV(\la)$ follow
from Corollary~\ref{onstdmodsc}.
The remaining parts then follow by mimicking the arguments used to prove
\cite[Lemma 4.9]{BS3}, using Lemma~\ref{jsf} in place of
\cite[Lemma 4.8]{BS3}.
 \end{proof}

\begin{Corollary}\label{cgt2}
Given $\la \in \La$,
pick $p, q$, $d$, $i_1,\dots,i_d$ and $r$ as in Lemma~\ref{cgt}.
Then we have that
$\cF_{i_d} \cdots \cF_{i_1} \cV(\la_{p,q}) \cong  \cP(\la)^{\bigoplus 2^r}$.
\end{Corollary}

\begin{proof}
We note as $\la_{p,q}$ is of defect zero that it is the
only weight in its block.
Using also (\ref{bggrec}), this implies that $\cP(\la_{p,q}) =
\cV(\la_{p,q})$.
Given this, the corollary
follows from Lemma~\ref{ga2} in exactly the same way
that Lemma~\ref{cgt} was deduced from Lemma~\ref{ga}.
\end{proof}

\phantomsubsection{Identification of Grothendieck groups}
Consider the Grothendieck group $[\sF]$ of
$\sF$. It is the
free $\Z$-module on basis $\{[\cL(\la)]\:|\:\la \in \La\}$.
The exact functors $\cF_i$ and $\cE_i$ (resp.\ $F_i$ and $E_i$) induce endomorphisms of the
Grothendieck group $[\sF]$ (resp.\ $[\Rep{K}]$),
which we denote by the same notation.
The last part of the
following theorem recovers the main result of \cite{B}.

\begin{Theorem}\label{cog}
Define a $\Z$-module isomorphism
$\iota:[\sF] \stackrel{\sim}{\rightarrow} [\Rep{K}]$
by declaring that
$\iota([\cL(\la)]) = [L(\la)]$ for each $\la \in \La$.
\begin{itemize}
\item[(i)]
We have that
 $\iota([\cV(\la)]) = [V(\la)]$ and
$\iota([\cP(\la)])= [P(\la)]$
for each $\la \in \La$.
\item[(ii)]
For each $i \in \Z$, we have that
$F_i \circ \iota = \iota \circ \cF_i$
and $E_i \circ \iota = \iota \circ \cE_i$ as linear maps
from $[\sF]$ to $[\Rep{K}]$.
\item[(iii)]
We have in $[\sF]$ that
\begin{equation*}
\displaystyle[\cP(\la)] = \sum_{\mu \supset \la} [\cV(\mu)],
\qquad
\displaystyle[\cV(\la)] = \sum_{\mu \subset \la}
[\cL(\mu)]
\end{equation*}
for each $\la \in \La$.
\end{itemize}
\end{Theorem}

\begin{proof}
Given $\la \in \La$,
let $p, q$, $d$, $r$ and $i_1,\dots,i_d$ be as in Lemma~\ref{cgt}.
By Lemma~\ref{cgt} and Theorem~\ref{form1}, we know already that
\begin{equation}\label{sec}
[P(\la)] = \frac{1}{2^r}
\cdot
F_{i_d} \cdots F_{i_1} [V(\la_{p,q})] = \sum_{\mu \supset \la}
[V(\mu)],
\end{equation}
all equalities written in $[\Rep{K}]$.
In view of Lemma~\ref{ga2},
the action of
$F_i$
on the classes of
standard modules in $[\Rep{K}]$
is described by exactly the same matrix
as the action of $\cF_i$
on the classes of standard modules in $[\sF]$.
So we deduce from the second equality in (\ref{sec})
that
$$
\frac{1}{2^r} \cdot \cF_{i_d} \cdots \cF_{i_1} [\cV(\la_{p,q})]
=
\sum_{\mu \supset \la} [\cV(\mu)],
$$
equality in $[\sF]$.
By Corollary~\ref{cgt2} this also equals $[\cP(\la)]$, proving the
first formula in (iii). The second formula
in (iii) follows
from the first and (\ref{bggrec}).

Then (i) is immediate from the definition of $\iota$ and
the coincidence of the formulae in (iii) and Theorem~\ref{form1}.

Finally to deduce (ii), we have already noted that
$\iota (F_i [V(\la)]) = \cF_i [\cV(\la)]$
for every $\la$.
It follows easily from this that
$\iota (F_i [P(\la)]) = \cF_i [\cP(\la)]$
for every $\la$. Using also
 the adjointness of $F_i$ and $E_i$
(resp.\ $\cF_i$ and $\cE_i$) we deduce that
\begin{align*}
[E_i L(\mu):L(\la)] &= \dim \hom_{K}(P(\la), E_i L(\mu))\\ &=
\dim \hom_{K}(F_i P(\la), L(\mu))=
\dim \hom_{G}(\cF_i \cP(\la), \cL(\mu))\\
&=\dim \hom_{G}(\cP(\la), \cE_i \cL(\mu))
= [\cE_i \cL(\mu):\cL(\la)]
\end{align*}
for every $\la,\mu \in \La$.
This is enough to show that
$\iota( E_i[L(\mu)]) = \cE_i [\cL(\mu)]$
for every $\mu$, which implies (ii) for $E_i$ and $\cE_i$. The argument for
$F_i$ and $\cF_i$ is similar.
\end{proof}

\phantomsubsection{Highest weight structure and duality}
At this point, we can also deduce the following result, which recovers
\cite[Theorem 4.47]{B}.

\begin{Theorem}\label{hwc}
The category $\sF$ is a highest weight category
in the sense of \cite{CPS} with weight poset
$(\La, \leq)$. The modules
$\{\cL(\la)\}, \{\cV(\la)\}$ and $\{\cP(\la)\}$
give its irreducible, standard and projective indecomposable modules,
respectively.
\end{Theorem}

\begin{proof}
We already noted just before (\ref{bggrec})
that $\cP(\la)$ has
a standard flag with $\cV(\la)$ at the top.
Moreover by
Theorem~\ref{cog}(iii) all the other
sections of this flag are all of the form
$\cV(\mu)$
with $\mu > \la$ in the Bruhat order.
The theorem follows from this, (\ref{bggrec})
and the definition of highest weight category.
\end{proof}

The costandard modules in the highest weight category $\sF$
can be constructed explicitly as the duals $\cV(\la)^{\circledast}$
of the standard modules with respect to a natural duality
$\circledast$. This duality maps a $G$-module $M$
to the linear dual $M^*$ with the action of $G$ defined
using the {\em supertranspose} anti-automorphism $g \mapsto
g^{st}$, where
$$
g^{st} =
\left(
\begin{array}{r|r}
a^t&-c^t\\\hline
b^t&d^t
\end{array}
\right)
$$
for $g$ of the form (\ref{supermat}).
Note $\circledast$ fixes irreducible modules, i.e.
$\cL(\la)^\circledast \cong \cL(\la)$ for each $\la \in \La$.

\section{Cyclotomic Hecke algebras and level two Schur-Weyl duality}

Fix integers $p \leq q$
and let $\la_{p,q}$ be the weight of defect zero from (\ref{laab}).
The standard module
$\cV(\la_{p,q})$ is projective.
As the functor $? \otimes \nV$ sends projectives to projectives,
the $G$-module
$\cV(\la_{p,q}) \otimes \nV^{\otimes d}$
is again projective
for any $d \geq 0$.
We want to describe its endomorphism algebra.

\phantomsubsection{Action of the degenerate affine Hecke algebra}
We begin by constructing an explicit
basis for $\cV(\la_{p,q}) \otimes \nV^{\otimes d}$.
Recalling (\ref{kacmod}), we have that
\begin{equation}\label{lab}
\cV(\la_{p,q}) = U(\mathfrak{g})\otimes_{U(\mathfrak{p})} E(\la_{p,q}).
\end{equation}
Let $\det_m$ (resp.\ $\det_n$) denote the
one-dimensional
$G_{\bar 0}$-module defined
by taking the determinant of $GL(m)$
(resp.\ $GL(n)$), with $\Z_2$-grading concentrated in degree $\0$.
Then the module $E(\la_{p,q})$ in (\ref{lab})
 is the inflation to $P$ of the module
$\Pi^{n(q+m)} (\det_m^{p} \otimes \det_n^{-(q+m)})$, so it
is also one dimensional.
Hence, fixing
a non-zero highest weight vector $v_{p,q} \in \cV(\la_{p,q})$, the
induced module
$\cV(\la_{p,q})$ is of dimension $2^{mn}$ with basis
\begin{equation}\label{labbasis}
\left\{
\prod_{r=m+1}^{m+n} \prod_{s=1}^m
e_{r,s}^{\tau_{r,s}}\cdot v_{p,q}\:\bigg|\:0 \leq \tau_{r,s} \leq 1\right\},
\end{equation}
where the products here are taken in any fixed order
(changing the order only changes the vectors by $\pm 1$).
Recall also that $v_1,\dots,v_{m+n}$ is the standard basis
for the natural module $\nV$, from which we get the obvious
monomial basis
\begin{equation}\label{mon}
\{v_{i_1} \otimes \cdots \otimes v_{i_d}\:|\:1 \leq i_1,\dots,i_d \leq m+n\}
\end{equation}
for $\nV^{\otimes d}$. Tensoring (\ref{labbasis})
and (\ref{mon}), we get the desired
basis for $\cV(\la_{p,q}) \otimes \nV^{\otimes d}$.

Now let $H_d$ be the {\em degenerate affine Hecke algebra}
from \cite{D2}.
This is the associative algebra
equal as a vector space to
$\C[x_1,\dots,x_d] \otimes \C S_d$,
the tensor product
of a polynomial algebra and the group algebra of the symmetric group $S_d$.
Multiplication
is defined so that
 $\C[x_1,\dots,x_d] \equiv \C[x_1,\dots,x_d] \otimes 1$ and $\C S_d
\equiv 1 \otimes \C S_d$ are subalgebras of $H_d$,
and also
\begin{equation*}
s_r x_{s} = x_{s} s_r \:\:\text{if $s \neq r,r+1$},
\qquad
s_r x_{r+1} = x_{r} s_r + 1,
\end{equation*}
where $s_r$
denotes the $r$th basic transposition $(r\:\,r+1)$.

By \cite[Proposition 5.1]{CW},
there is a right action
of $H_d$ on $\cV(\la_{p,q})\otimes \nV^{\otimes d}$
by $G$-module endomorphisms.
The transposition $s_r$ acts as the
``super'' flip
$$
(v \otimes v_{i_1} \otimes \cdots \otimes v_{i_r} \otimes v_{i_{r+1}}\otimes\cdots \otimes v_{i_d}) s_r
= (-1)^{\bar i_r \bar i_{r+1}}
v \otimes v_{i_1} \otimes \cdots \otimes v_{i_{r+1}} \otimes v_{i_{r}}\otimes\cdots \otimes v_{i_d}.
$$
This is the same as the endomorphism
defined by left multiplication by the
element $\Omega$ from (\ref{omega}) so that the first and second tensors in
$\Omega$ hit the $(r+1)$th and $(r+2)$th tensor positions in
$\cV(\la_{p,q}) \otimes \nV^{\otimes d}$, respectively.
The polynomial generator $x_s$ acts by left multiplication
by $\Omega$ so that the first tensor in $\Omega$ is spread across
tensor positions $1,\dots,s$ using the comultiplication of $U(\mathfrak{g})$
and the second tensor in $\Omega$ hits the $(s+1)$th tensor position
in
$\cV(\la_{p,q}) \otimes \nV^{\otimes d}$. The following lemma gives an explicit
formula for the action of $x_s$ in a special case.

\begin{Lemma}\label{litlem}
For $1 \leq i_1,\dots,i_d \leq m+n$ and $1 \leq s \leq d$, we have that
\begin{multline*}
(v_{p,q} \otimes v_{i_1} \otimes\cdots\otimes v_{i_d}) x_s
=
p v_{p,q} \otimes v_{i_1}\otimes\cdots\otimes v_{i_d}\\
+ \sum_{r=1}^{s-1} (-1)^{\bar i_r \bar i_s + \sum_{r < t < s} (\bar i_r + \bar i_s) \bar i_t}
v_{p,q} \otimes v_{i_1} \otimes \cdots \otimes v_{i_s} \otimes \cdots\otimes
v_{i_r} \otimes \cdots \otimes v_{i_d}
\end{multline*}
if $1 \leq i_s \leq m$, and
\begin{multline*}
(v_{p,q} \otimes v_{i_1} \otimes\cdots\otimes v_{i_d}) x_s
=
(q+m) v_{p,q} \otimes v_{i_1}\otimes\cdots\otimes v_{i_d}\\
+ \sum_{r=1}^{s-1} (-1)^{\bar i_r \bar i_s + \sum_{r < t < s} (\bar i_r + \bar i_s) \bar i_t}
v_{p,q} \otimes v_{i_1} \otimes \cdots \otimes v_{i_s} \otimes \cdots\otimes
v_{i_r} \otimes \cdots \otimes v_{i_d}\\
+\sum_{j=1}^m
(-1)^{n(q+m)+\bar i_1+\cdots+\bar i_{s-1}} (e_{i_s,j} \cdot v_{p,q})
\otimes v_{i_1}\otimes\cdots\otimes  v_j
\otimes\cdots\otimes v_{i_d}
\end{multline*}
if $m+1 \leq i_s \leq m+n$. (In the first two summations we have
interchanged $v_{i_r}$ and $v_{i_s}$, while in the last one
we have replaced $v_{i_s}$ by $v_j$.)
\end{Lemma}

\begin{proof}
Note for any $1 \leq i,j \leq m+n$ that
$$
e_{i,j}\cdot v_{p,q}
=
\left\{
\begin{array}{ll}
pv_{p,q}&\text{if $1 \leq i = j \leq m$,}\\
-(q+m)v_{p,q}&\text{if $m+1 \leq i = j \leq m+n$,}\\
e_{i,j}\cdot v_{p,q}&\text{if $m+1 \leq i \leq m+n$ and $1 \leq j \leq m$,}\\
0&\text{otherwise.}
\end{array}\right.
$$
Using this,
the lemma is a routine calculation (taking care with superalgebra signs).
\end{proof}

\begin{Corollary}\label{kernel}
The element $(x_1-p)(x_1-q) \in H_d$ acts as zero
on $\cV(\la_{p,q}) \otimes \nV^{\otimes d}$.
\end{Corollary}

\begin{proof}
It suffices to check this in the special case that $d=1$.
In that case,
Lemma~\ref{litlem} shows that
$$
(v_{p,q} \otimes v_i) x_1
=
\left\{\begin{array}{ll}
p v_{p,q}\otimes v_i
&\text{if $1 \leq i \leq m$,}\\
q v_{p,q}\otimes v_i\\
\quad+ \sum_{j=1}^m (-1)^{n(q+m)} e_{i,j} (v_{p,q} \otimes v_{j})&\text{if $m+1 \leq i \leq m+n$.}
\end{array}\right.
$$
It follows easily that
$(x_1-p)(x_1-q)$ acts
as zero on the vector $v_{p,q} \otimes v_i$ for every $1 \leq i \leq m+n$.
These vectors generate $\cV(\la_{p,q}) \otimes \nV$
as a $G$-module so we deduce that
$(x_1-p)(x_1-q)$ acts
as zero the whole module.
\end{proof}

\begin{Corollary}\label{indep}
If $d \leq \min(m,n)$
then the endomorphisms of $\cV(\la_{p,q}) \otimes \nV^{\otimes d}$
defined by right multiplication by
$\{x_1^{\sigma_1} \cdots x_d^{\sigma_d} w \:|\:
0 \leq \sigma_1,\dots,\sigma_d \leq 1, w \in S_d\}$
are linearly independent.
\end{Corollary}

\begin{proof}
Any vector $v \in \cV(\la_{p,q}) \otimes \nV^{\otimes d}$
can be written as $v = \sum_{i \in I} b_i \otimes c_i$
where $\{b_i\:|\:i \in I\}$ is the basis from (\ref{labbasis}) and
the $c_i$'s are unique vectors in $\nV^{\otimes d}$.
We refer to $c_i$ as the {\em $b_i$-component} of $v$.
Exploiting the assumption on $d$, we can pick distinct integers $m+1 \leq i_1,\dots,i_d \leq m+n$
and $1 \leq j_1,\dots,j_d \leq m$.
Take $0 \leq \sigma_1,\dots,\sigma_d \leq 1$
and consider the vector
$$
(v_{p,q}\otimes v_{i_1}\otimes\cdots\otimes v_{i_d})
x_1^{\sigma_1} \cdots x_d^{\sigma_d}.
$$
For $0 \leq \tau_1,\dots,\tau_d \leq 1$,
Lemma~\ref{litlem} implies that
the $e_{i_1,j_1}^{\tau_1}\cdots e_{i_d,j_d}^{\tau_d}\cdot v_{p,q}$-component
of
$(v_{p,q}\otimes v_{i_1}\otimes\cdots\otimes v_{i_d})
x_1^{\sigma_1} \cdots x_d^{\sigma_d}$
is zero either if $\tau_1+\cdots+\tau_d > \sigma_1+\cdots+\sigma_d$,
or if $\tau_1+\cdots+\tau_d = \sigma_1+\cdots+\sigma_d$
but $\tau_r \neq \sigma_r$ for some $r$.
Moreover, if $\tau_r = \sigma_r$ for all $r$, then
the $e_{i_1,j_1}^{\tau_1}\cdots e_{i_d,j_d}^{\tau_d}\cdot v_{p,q}$-component
of
$(v_{p,q}\otimes v_{i_1}\otimes\cdots\otimes v_{i_d})
x_1^{\sigma_1} \cdots x_d^{\sigma_d}$
is equal to $\pm v_{k_1}\otimes\cdots\otimes v_{k_d}$
where $k_r = i_r$ if $\sigma_r = 0$ and
$k_r = j_r$ if $\sigma_r = 1$.
This is enough to
show that the vectors
$(v_{p,q} \otimes v_{i_1}\otimes\cdots\otimes v_{i_d})x_1^{\sigma_1}\cdots x_d^{\sigma_d} w$
for all $0 \leq \sigma_1,\dots,\sigma_d \leq 1$ and $w \in S_d$
are linearly independent, and the corollary follows.
\end{proof}

In view of Corollary~\ref{kernel}, the right action of $H_d$
on $\cV(\la_{p,q}) \otimes \nV^{\otimes d}$
induces an action of the quotient algebra
\begin{equation}\label{cyclo}
H_d^{p,q} := H_d / \langle (x_1-p)(x_1-q) \rangle.
\end{equation}
This algebra is a particular example of
a {\em degenerate cyclotomic Hecke algebra} of level two.
It is well known (e.g. see \cite[Lemma 3.5]{BKschur}) that
$\dim H_d^{p,q} = 2^d d!$.

\begin{Corollary}\label{faith}
If $d \leq \min(m,n)$
the action of $H_d^{p,q}$ on
$\cV(\la_{p,q}) \otimes \nV^{\otimes d}$ is faithful.
\end{Corollary}

\begin{proof}
This follows on comparing the dimension of $H_d^{p,q}$ with the number of linearly
independent endomorphisms constructed in Corollary~\ref{indep}.
\end{proof}

Since the action of $H^{p,q}_d$ on $\cV(\la_{p,q})\otimes V^{\otimes d}$
is by $G$-module endomorphisms, it induces an algebra homomorphism
\begin{equation}\label{Phi}
\Phi:H_d^{p,q} \rightarrow
\End_G(\cV(\la_{p,q}) \otimes V^{\otimes d})^{\op}.
\end{equation}
The main goal in the remainder of the section is to show that this homorphism is {\em surjective}.

\phantomsubsection{\boldmath Weight idempotents and
the space $\nT_d^{p,q}$}
For a tuple $\bi = (i_1,\dots,i_d) \in \Z^d$, there is
an idempotent $e(\bi) \in H_d^{p,q}$ determined uniquely by the
property that
multiplication by $e(\bi)$ projects any $H_d^{p,q}$-module
onto its {\em $\bi$-weight space}, that is, the simultaneous
generalised eigenspace for the commuting operators
$x_1,\dots,x_d$ and eigenvalues $i_1,\dots,i_d$, respectively.
All but finitely many of the $e(\bi)$'s are zero, and the non-zero ones give
a system of mutually orthogonal idempotents in $H_d^{p,q}$ summing to
$1$; see e.g. \cite[$\S$3.1]{BKinv}.

The action of the idempotent $e(\bi)$ on the module
$\cV(\la_{p,q}) \otimes
\nV^{\otimes d}$ can be interpreted as follows.
In view of Corollary~\ref{fullt}, we have that
\begin{equation}\label{wtdec}
\cV(\la_{p,q}) \otimes \nV^{\otimes d} = \bigoplus_{\bi \in \Z^d}
\cF_{\bi} \cV(\la_{p,q})
\end{equation}
where $\cF_{\bi}$ denotes the composite
$\cF_{i_d} \circ \cdots \circ \cF_{i_1}$ of the functors
from (\ref{spf}). By Lemma~\ref{eig} and the definition
of the actions of $x_1,\dots,x_d$, the summand $\cF_{\bi} \cV(\la_{p,q})$
in this decomposition is precisely the $\bi$-weight space of
$\cV(\la_{p,q}) \otimes \nV^{\otimes d}$.
Hence the weight idempotent $e(\bi)$ acts on
$\cV(\la_{p,q}) \otimes \nV^{\otimes d}$ as the projection onto the summand
$\cF_{\bi} \cV(\la_{p,q})$
along the decomposition (\ref{wtdec}).

Recalling the interval $I_{p,q}$ from (\ref{intervals}),
we are usually from now on going to restrict our attention to the
summand
\begin{equation}\label{littlet}
\nT_d^{p,q} := \bigoplus_{\bi \in (I_{p,q})^d} \cF_{\bi} \cV(\la_{p,q})
\end{equation}
of $\cV(\la_{p,q}) \otimes \nV^{\otimes d}$.
By the discussion in the previous paragraph, we have equivalently that
$\nT^{p,q}_d = (\cV(\la_{p,q}) \otimes \nV^{\otimes d}) 1^{p,q}_d$
where
\begin{equation}\label{theidemp}
1^{p,q}_d := \sum_{\bi \in (I_{p,q})^d} e(\bi) \in H_d^{p,q}.
\end{equation}
As a consequence of the fact that
any symmetric polynomial in $x_1,\dots,x_d$ is central
in $H_d$, the idempotent $1^{p,q}_d$ is
{central} in $H_d^{p,q}$. The space $\nT^{p,q}_d$ is naturally a right module
over $1^{p,q}_d H_d^{p,q}$,
which is a sum of blocks of $H_d^{p,q}$.
Hence the map $\Phi$ from (\ref{Phi}) induces an algebra homomorphism
\begin{equation}\label{rep}
1^{p,q}_d H_d^{p,q} \rightarrow
\End_G(\nT_d^{p,q})^{\op}.
\end{equation}
As a refinement of the surjectivity of $\Phi$ proved below, we will also see later in the section
that the induced map (\ref{rep}) is an {\em isomorphism}. Note from (\ref{ident}) onwards we
will denote the algebra $1^{p,q}_d H^{p,q}_d$ instead by $R^{p,q}_d$.

\phantomsubsection{Stretched diagrams}
In this subsection, we develop some combinatorial tools which will be used initially to compute the dimension
of the various endomorphism algebras that we are interested in.
We say that a tuple $\bi\in \Z^d$
is {\em $(p,q)$-admissible} if
$i_r$ is $\Ga_{r-1}$-admissible for each $r=1,\dots,d$,
where $\Ga_0,\dots,\Ga_d$ are defined recursively
from
$\Ga_0 := \{\la_{p,q}\}$ and
$\Ga_r := \Ga_{r-1} - \alpha_{i_r}$,
notation as in (\ref{CKLR}). We refer to the sequence
$\bGa := \Ga_d \cdots \Ga_1 \Ga_0$ of blocks
here
as the
{\em associated block sequence}.
The composite matching $\bt = t_d \cdots t_1$
defined by setting $t_r := t_{i_r}(\Ga_{r-1})$ for each $r$
is the {\em associated composite matching}.
Both of these things make sense only if $\bi \in \Z^d$ is
$(p,q)$-admissible.

\begin{Lemma}\label{easyzero}
If $\bi \in \Z^d$ is not $(p,q)$-admissible then
$\cF_{\bi} \cV(\la_{p,q})$ is zero.
\end{Lemma}

\begin{proof}
This follows from the definitions and (\ref{spf}).
\end{proof}

By a {\em stretched cap diagram} $\bt = t_d \cdots t_1$ of height $d$, we mean
the associated composite matching
for some $(p,q)$-admissible
sequence $\bi \in \Z^d$.
We can uniquely recover the sequence $\bi$, hence also the
associated block sequence $\bGa$, from
the stretched cap diagram
$\bt$.
Here is an example of
a stretched cap diagram of height $5$, taking
$m=2,n=1$ and $q-p=1$; we
draw only the strip containing the vertices
indexed by $I_{p,q}^+$,
as the picture outside of this strip consists only of vertical lines,
and also label the horizontal number lines by the associated block
sequence
$\bGa = \Ga_5 \cdots \Ga_0$.
$$
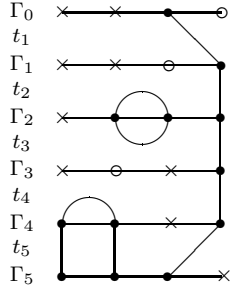
\begin{picture}(40,110)
\put(-20,102){\text{$_{\Ga_0}$}}
\put(-19,92.5){\text{$_{t_1}$}}
\put(-20,82){\text{$_{\Ga_1}$}}
\put(-19,72.5){\text{$_{t_2}$}}
\put(-20,62){\text{$_{\Ga_2}$}}
\put(-19,52.5){\text{$_{t_3}$}}
\put(-20,42){\text{$_{\Ga_3}$}}
\put(-19,32.5){\text{$_{t_4}$}}
\put(-20,22){\text{$_{\Ga_4}$}}
\put(-19,12.5){\text{$_{t_5}$}}
\put(-20,2){\text{$_{\Ga_5}$}}
\put(0,2){\line(1,0){60}}
\put(0,22){\line(1,0){60}}
\put(0,42){\line(1,0){60}}
\put(0,62){\line(1,0){60}}
\put(0,82){\line(1,0){60}}
\put(0,102){\line(1,0){60}}
\put(-3.3,100.2){$\scriptstyle\times$}
\put(16.7,100.2){$\scriptstyle\times$}
\put(37.8,100){$\scriptstyle\bullet$}
\put(57.5,99.3){$\circ$}
\put(40,102){\line(1,-1){20}}
\put(-3.3,80.2){$\scriptstyle\times$}
\put(16.7,80.2){$\scriptstyle\times$}
\put(57.8,80){$\scriptstyle\bullet$}
\put(37.5,79.3){$\circ$}
\put(30,62){\oval(20,20)[t]}
\put(59.5,82){\line(0,-1){20}}
\put(17.8,60){$\scriptstyle\bullet$}
\put(37.8,60){$\scriptstyle\bullet$}
\put(-3.3,60.2){$\scriptstyle\times$}
\put(57.5,60){$\scriptstyle\bullet$}
\put(30,62){\oval(20,20)[b]}
\put(59.5,62){\line(0,-1){20}}
\put(37.7,40.2){$\scriptstyle\times$}
\put(17.5,39.3){$\circ$}
\put(-3.3,40.2){$\scriptstyle\times$}
\put(57.5,40){$\scriptstyle\bullet$}
\put(60,22){\line(-1,-1){20}}
\put(37.7,20.2){$\scriptstyle\times$}
\put(57.5,20){$\scriptstyle\bullet$}
\put(10,22){\oval(20,20)[t]}
\put(19.5,22){\line(0,-1){20}}
\put(59.5,42){\line(0,-1){20}}
\put(-0.5,22){\line(0,-1){20}}
\put(57.7,.2){$\scriptstyle\times$}
\put(37.5,0){$\scriptstyle\bullet$}
\put(-2.5,0){$\scriptstyle\bullet$}
\put(17.5,0){$\scriptstyle\bullet$}
\put(-2.5,20){$\scriptstyle\bullet$}
\put(17.5,20){$\scriptstyle\bullet$}
\end{picture}
$$
By a {\em generalised cap} in a stretched cap diagram we
 mean a component
that meets the bottom number line
at two different vertices.
An {\em oriented stretched cap diagram}
is a consistently oriented diagram of the form
$$
\bt[\bga] = \ga_d t_d \ga_{d-1}\cdots\ga_1 t_1 \ga_0
$$
where $\bga = \ga_d \cdots \ga_0$ is a sequence of
weights chosen from the associated block sequence
$\bGa = \Ga_d \cdots \Ga_0$, i.e.
$\ga_r \in \Ga_r$ for each $r=0,\dots,d$.
In other words, we decorate the number lines of $\bt$ by weights
from the appropriate blocks, in such a way that the resulting
diagram is consistently oriented. (For a precise definition of the term oriented we refer to \cite[$\S$2]{BS1}).

\begin{Theorem}\label{Dec}
There are $G$-module isomorphisms
\begin{align*}
\cV(\la_{p,q}) \otimes \nV^{\otimes d}
&\cong \bigoplus_{\la \in \La,\,
\deg(\la) = \deg(\la_{p,q})+d}
\cP(\la)^{\oplus \dd(\la)},\\
\nT_d^{p,q} &\cong
\!\!\!\bigoplus_{\la \in \La_{p,q},\,\deg(\la) = \deg(\la_{p,q})+d}\!\!\!
\cP(\la)^{\oplus \dd(\la)},
\end{align*}
where
$\dd(\la)$
is the number of
oriented stretched cap diagrams $\bt[\bga]$ of height $d$
such that $\ga_0 = \la_{p,q}$,
$\ga_d = \la$, and all generalised caps are anti-clockwise.
\end{Theorem}

\begin{proof}
For the first isomorphism,
in view of Theorem~\ref{cog} and Corollary~\ref{fullt}, it suffices to prove the analogous statement
on the diagram algebra side, namely, that
\begin{align}\label{init}
\bigoplus_{\bi \in \Z^d} F_{\bi} V(\la_{p,q})
&\cong \bigoplus_{\la \in \La,\,\deg(\la) = \deg(\la_{p,q})+d}
P(\la)^{\oplus \dd(\la)}
\end{align}
as $K$-modules.
Remembering that $V(\la_{p,q}) = P(\la_{p,q})$, this follows as an application
of \cite[Theorem~4.2]{BS2},
first
using \cite[Theorem~3.5]{BS2} and \cite[Theorem~3.6]{BS2} to write
the composite projective
functor $F_{\bi} = F_{i_d} \circ \cdots \circ F_{i_1}$
in terms of
indecomposable projective functors.

The proof of the second isomorphism is similar, taking only
$\bi \in (I_{p,q})^d$ in (\ref{init}). It is helpful to note that if $\la \in
\La_{p,q}$ and $\bt[\bga]$ is one of the oriented stretched cap
diagrams counted by $\dd(\la)$ then $\bt[\bga]$ is trivial outside
the strip containing the vertices indexed by $I_{p,q}^+$, i.e. it
consists only of straight lines oriented $\up$ outside that region. This follows by
considering (\ref{CKLR}).
\end{proof}

\begin{Corollary}\label{summands}
The modules
$\{\cP(\la)\:|\:\la \in \La_{p,q}^\circ,
\deg(\la) = \deg(\la_{p,q})+d\}$
give a complete set of representatives for the isomorphism
classes of indecomposable direct summands of $\nT_d^{p,q}$.
\end{Corollary}

\begin{proof}
Suppose we are given
$\la \in \La_{p,q}^{\circ}$ with $\deg(\la) =
\deg(\la_{p,q})+d$.
Applying Corollary~\ref{cgt2}, there is
a sequence $\bi = (i_1,\dots,i_d) \in (I_{p,q})^d$ such that
$\cP(\la)$ is a summand of
$\cF_\bi \cV(\la_{p,q})$. Hence $\cP(\la)$ is a summand of
$\nT_d^{p,q}$.
Conversely,
applying Theorem~\ref{Dec}, we take
$\la \in \La_{p,q}$ with
$\deg(\la) = \deg(\la_{p,q})+d$ and
$\dd(\la) \neq 0$, and must show that
$\la \in \La_{p,q}^\circ$.
There exists an oriented stretched cap diagram
$\bt[\bga]$ of height $d$ with $\ga_0 = \la_{p,q}$ and
$\ga_d = \la$, all of whose generalised caps are anti-clockwise.
Every vertex labelled $\down$ in $\la$ must be at the
left end of one of these anti-clockwise generalised caps,
the right end of which gives a vertex
labelled $\up$ indexed by an integer $\leq q+n$.
Recalling the definition (\ref{lamn2}),
these observations prove that $\la \in \La_{p,q}^\circ$.
\end{proof}

\begin{Corollary}\label{zer}
$\nT_d^{p,q} = \{0\}$
for $d > (m+n)(q-p)+2mn$.
\end{Corollary}

\begin{proof}
The set $\La_{p,q}$ has a unique element
$\mu_{p,q}$ of maximal height, namely, the weight
$$
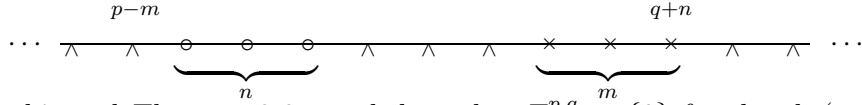
\begin{picture}(-350,40)
\put(-294,28){$_{p-m}$}
\put(-90,28){$_{q+n}$}
\put(-333,12.5){$\cdots$}
\put(-22,12.5){$\cdots$}
\put(-270,8.2){$\underbrace{\phantom{hellow worl}}_{n}$}
\put(-133,8.2){$\underbrace{\phantom{hellow worl}}_{m}$}
\put(-313,15.2){\line(1,0){283}}
\put(-288.7,10.7){$\scriptstyle\up$}
\put(-311.7,10.7){$\scriptstyle\up$}
\put(-153.7,10.7){$\scriptstyle\up$}
\put(-176.7,10.7){$\scriptstyle\up$}
\put(-199.7,10.7){$\scriptstyle\up$}
\put(-268.2,12.4){$\circ$}
\put(-245.2,12.4){$\circ$}
\put(-222.2,12.4){$\circ$}
\put(-131.4,13.4){$\scriptstyle\times$}
\put(-108.4,13.4){$\scriptstyle\times$}
\put(-85.4,13.4){$\scriptstyle\times$}
\put(-61.7,10.7){$\scriptstyle\up$}
\put(-38.7,10.7){$\scriptstyle\up$}
\end{picture}
$$
Using this and Theorem~\ref{Dec},
we deduce that $\nT_d^{p,q} = \{0\}$
for $d > \deg(\mu_{p,q}) - \deg(\la_{p,q}) =
(m+n)(q-p) +2mn$.
\end{proof}

The mirror image of the oriented stretched cap diagram
$\bu[\bde]$ in a horizontal axis is denoted
$\bu^*[\bde^*]$. We call it an {\em oriented stretched cup diagram}.
Then an {\em oriented stretched circle diagram} of height $d$
means a composite diagram
of the form\label{tprb}
$$
\bu^*[\bde^*] \wr \bt[\bga]
= \delta_0 u_1^* \delta_1 \cdots \delta_{d-1} u_d^* \ga_d
t_d \ga_{d-1} \cdots \ga_1 t_1 \ga_0
$$
where $\bt[\bga]$ and $\bu[\bde]$ are
oriented stretched cap diagrams of height $d$ with $\ga_d = \delta_d$;
see \cite[(6.17)]{BS3} for an example.

\begin{Theorem}\label{dimformt}
The dimension of the algebra
$\End_G(\nT_d^{p,q})^{\op}$
is equal to the number of oriented stretched circle diagrams
$\bu^*[\bde^*] \wr \bt [\bga]$
of height $d$ such that
$\ga_0 = \delta_0 = \la_{p,q}$
and $\ga_d = \delta_d \in \La_{p,q}$.
\end{Theorem}

\begin{proof}
Applying Theorem~\ref{Dec}, we see that
the dimension of the endomorphism algebra is equal to
$$
\sum_{\la,\mu \in \La_{p,q},\,\deg(\la) = \deg(\mu) = \deg(\la_{p,q})+d}
\dd(\la) \cdot \dd(\mu) \cdot \dim \hom_G(\cP(\la), \cP(\mu)).
$$
Also in view of Theorem~\ref{cog},
$\dim \hom_G(\cP(\la), \cP(\mu)) = [\cP(\mu):\cL(\la)]$
is equal to the analogous dimension
$\dim \hom_{K}(P(\la), P(\mu)) = [P(\mu):L(\la)]$
on the diagram algebra side, which is described
explicitly by \cite[(5.9)]{BS1}.
We deduce that
$\dim \hom_G(\cP(\la), \cP(\mu))$ is equal to
the number of weights $\nu$
such that $\la \sim \nu \sim \mu$
and the circle diagram $\underline{\la} \nu \overline{\mu}$
is consistently oriented.
The theorem follows easily on combining this
with the combinatorial definitions
of $\dd(\la)$ and $\dd(\mu)$
from Theorem~\ref{Dec}.
\end{proof}

\phantomsubsection{\boldmath The algebra $R_d^{p,q}$ and the
isomorphism theorem}
Now we need to recall some of the main results of \cite{BS3} which
give an alternative diagrammatic
description of the algebra
$1_d^{p,q} H^{p,q}_d$.
This will allow us to see to start with that this
algebra has the same dimension as the
endomorphism algebra from
Theorem~\ref{dimformt}.
For the reader wishing to understand already the relationship between
the diagram algebra
$R_d^{p,q}$ defined in the next paragraph and the algebra $K(m|n)$ from the introduction,
we point to Lemma~\ref{itsall}, (\ref{splodge}) and Corollary~\ref{sef} in
the next section.
On the other hand Corollary~\ref{endoc}
established later in this section
explains the connection between
$R_d^{p,q}$ and the representation theory of $G$.

Let $R_d^{p,q}$ be the associative, unital algebra
with basis
$$
\left\{|\bu^*[\bde^*] \wr \bt[\bga]|\:\:\bigg|\:
\begin{array}{l}
\text{for all oriented stretched circle diagrams
$\bu^*[\bde^*]\wr \bt[\bga]$} \\
\text{of height $d$ with $\ga_0 = \delta_0 = \la_{p,q}$
and $\ga_d = \delta_d \in \La_{p,q}$}
\end{array}
\right\}.
$$\label{tpr}
The multiplication is defined
by an explicit algorithm
described in detail in \cite{BS3}.
Briefly, to multiply two basis vectors
$|\bs^*[\btau^*]\wr \br[\bsigma]|$ and
$|\bu^*[\bde^*]\wr \bt[\bga]|$, the product is zero
unless $\br = \bu$ and all mirror image pairs of internal circles
in $\br[\bsigma]$ and $\bu^*[\bde^*]$ are oriented so that one is clockwise,
the other anti-clockwise. Assuming these conditions hold, the product
is computed by putting $\bs^*[\btau^*] \wr \br[\bsigma]$
underneath
$\bu^*[\bde^*]\wr \bt[\bga]$, erasing all internal circles and number lines in
$\br[\bsigma]$ and $\bu^*[\bde^*]$, then iterating the generalised surgery procedure
to smooth out the symmetric middle
section of the diagram.

\begin{Lemma}\label{eqd}
The algebras
$R_d^{p,q}$ and $\End_G(\nT_d^{p,q})^{\op}$ have the same dimension.
In particular, $R_d^{p,q}$ is the zero algebra for $d > (m+n)(q-p)+2mn$.
\end{Lemma}

\begin{proof}
The number of elements in the diagram basis for $R_d^{p,q}$
is the same as the dimension of the algebra
$\End_G(\nT_d^{p,q})^{\op}$
thanks to Theorem~\ref{dimformt}.
The last statement follows from Corollary~\ref{zer}.
\end{proof}

As a consequence of \cite[Corollary 8.6]{BS3},
we can identify $R_d^{p,q}$ with
a certain {\em cyclotomic Khovanov-Lauda-Rouquier algebra}
in the sense of \cite{KhL, Rou}.
To make this identification explicit, we need to
define some special elements
\begin{equation}\label{klgens}
\{e(\bi)\:|\:\bi \in (I_{p,q})^d\} \cup \{y_1,\dots,y_d\}\cup
\{\psi_1,\dots,\psi_{d-1}\}
\end{equation}
in $R_d^{p,q}$. For $\bi \in (I_{p,q})^d$, we let $e(\bi) \in R_d^{p,q}$
be the idempotent defined as follows.
If $\bi$ is not $(p,q)$-admissible then $e(\bi) := 0$.
If it is admissible, let $\bt = t_d\cdots t_1$ be the
associated composite matching
and $\bGa = \Ga_d \cdots \Ga_0$ be the associated block sequence.
Then
\begin{equation}\label{idmsps}
e(\bi) := \sum_{\bde, \bga} |\bt^*[\bde^*] \wr \bt[\bga]|
\end{equation}
where the sum is over all sequences
$\bga = \ga_d \cdots \ga_0$ and
$\bde = \delta_d\cdots\delta_0$ of weights
with each $\ga_r, \delta_r \in \Ga_r$
chosen so that every circle of
$\bt^*[\bde^*] \wr \bt[\bga]$
crossing the middle number
line is anti-clockwise, and all remaining mirror image pairs of circles are
oriented so that one is clockwise, the other anti-clockwise.
The elements $\{e(\bi)\:|\:\bi \in (I_{p,q})^d\}$ give a system of mutually
orthogonal idempotents whose sum is the identity in $R_d^{p,q}$.

Next we define the elements $y_1,\dots,y_d$.
Let $\bar y_1,\dots,\bar y_d$ be the unique elements of $R_d^{p,q}$ such that the product
$|\bu^*[\bde^*]\wr \bt[\bga]| \cdot \bar y_r$ (resp.\ $\bar y_r \cdot |\bu^*[\bde^*]\wr \bt[\bga]|$)
is computed by making a {positive circle move} in the section of $\bu^* \bt$ containing $t_r$
(resp.\ $u_r^*$), as described in detail in \cite[(5.5), (5.11)]{BS3}. Also introduce the signs
\begin{equation}\label{signs}
\sigma_{p,q}^r(\bi) :=
(-1)^{\min(p,i_r)+\min(q,i_r)+m-p-\delta_{i_1,i_r} -\cdots
-\delta_{i_{r-1},i_r}}.
\end{equation}
Then we
define
$y_r := \sum_{\bi \in (I_{p,q})^d} y_r e(\bi)$
where
\begin{equation}
y_re(\bi) :=
\sigma_{p,q}^r(\bi)
\bar y_r e(\bi),\label{yr}
\end{equation}
to get the elements $y_r \in R_d^{p,q}$ for $r=1,\dots,d$.

Finally we define $\psi_1,\dots,\psi_{d-1}$.
Let $\bar \psi_1,\dots,\bar \psi_{d-1}$
be the unique elements of $R_d^{p,q}$ such that the product
$|\bu^*[\bde^*]\wr \bt[\bga]|
\cdot \bar \psi_r$
(resp.\ $\bar \psi_r \cdot |\bu^*[\bde^*]\wr \bt[\bga]|$)
is computed by making a {negative circle move},
a crossing move or a height move
in the section of $\bu^* \bt$
containing
$t_{r+1} t_{r}$
(resp.\ $u_r^* u_{r+1}^*$),
as described in detail in \cite[(5.7),(5.12)]{BS3}.
Then we define $\psi_r := \sum_{\bi \in (I_{p,q})^d} \psi_r e(\bi)$
where
\begin{equation}
\psi_re(\bi) :=
\left\{
\begin{array}{ll}
-\sigma_{p,q}^r(\bi)
\bar \psi_r e(\bi)&\text{if $i_{r+1} = i_r$ or $i_{r+1}=i_{r}+1$,}\\
\bar \psi_r e(\bi)&\text{otherwise,}\label{psir}
\end{array}
\right.
\end{equation}
to get the elements
$\psi_r \in R_d^{p,q}$ for $r=1,\dots,d-1$.

\begin{Theorem}\label{iso1}
The elements (\ref{klgens})
generate $R_d^{p,q}$
subject only to
the defining relations of the Khovanov-Lauda-Rouquier algebra
associated to the linear quiver
$$
\begin{picture}(80,10)
\put(0,4){$_\bullet$}
\put(16,4){$_\bullet$}
\put(32,4){$_\bullet$}
\put(0,1.3){$\longrightarrow$}
\put(16,1.3){$\longrightarrow$}
\put(38,1.4){$\cdots$}
\put(52,4){$_\bullet$}
\put(68,4){$_\bullet$}
\put(52,1.3){$\longrightarrow$}
\end{picture}
$$
with vertices indexed by the set $I_{p,q}$ in order from left to right (see e.g. \cite[(6.8)--(6.16)]{BS3}),
plus the additional cyclotomic relations
$y_1^{\delta_{i_1,p}+\delta_{i_1,q}} e(\bi) = 0$
for $\bi = (i_1,\dots,i_d) \in (I_{p,q})^d$.
\end{Theorem}

\begin{proof}
This is a consequence of \cite[Corollary 8.6]{BS3}.
More precisely, we apply \cite[Corollary 8.6]{BS3},
taking
the index set $I$ there to be the set
$I_{p,q}$, the pair $(o+m,o+n)$ there to be
$(p,q)$, and summing over all $\alpha \in Q_+$
of height $d$.
This implies that the given quotient of the Khovanov-Lauda-Rouquier
algebra is isomorphic to the diagram algebra with basis consisting
of oriented stretched circle diagrams
$|\bu^*[\bde^*] \wr \bt[\bga]|$ just like the ones considered here,
except they are drawn only in the strip containing the vertices
indexed by $I_{p,q}^+$.
The isomorphism in \cite{BS3}
is not quite the same as the map here,
because
the sign in (\ref{signs})
differs from the corresponding sign
chosen in \cite{BS3} by a factor of $(-1)^{m-p}$;
this causes no problems as it amounts to twisting
by an automorphism of
the Khovanov-Lauda-Rouquier algebra.
It remains to observe that all the
oriented stretched circle
diagrams in the statement of the present
theorem are trivial outside the strip $I_{p,q}^+$, consisting only
of straight lines
oriented $\up$ in that region; these have no
effect on the multiplication.
\end{proof}

Now we can formulate the following key {\em isomorphism theorem},
which identifies the algebras
$1_d^{p,q} H^{p,q}_d$ and
$R_d^{p,q}$.

\begin{Theorem}\label{iso2}
There is a unique algebra isomorphism
$$
1_d^{p,q} H^{p,q}_d \stackrel{\sim}{\rightarrow}
R_d^{p,q}
$$
such that $e(\bi) \mapsto e(\bi),
x_r e(\bi) \mapsto (y_r+i_r) e(\bi)$
and $s_r e(\bi) \mapsto (\psi_r q_r(\bi) - p_r(\bi))e(\bi)$
for each $r$ and $\bi \in (I_{p,q})^d$, where
$p_r(\bi), q_r(\bi) \in R_d^{p,q}$ are chosen as
in \cite[$\S$3.3]{BKinv}, e.g. one could take
\begin{align*}
p_r(\bi) &:= \left\{
\begin{array}{ll}
1&\text{if $i_r = i_{r+1}$,}\\
-(i_{r+1}-i_r + y_{r+1}-y_r)^{-1}\hspace{6mm}&\text{if $i_r \neq i_{r+1}$;}
\end{array}\right.\\
q_r(\bi) &:= \left\{
\begin{array}{ll}
1+y_{r+1}-y_r&\text{if $i_r = i_{r+1}$,}\\
(2+y_{r+1}-y_r)(1+y_{r+1}-y_r)^{-2}&\text{if $i_{r+1}=i_r+1$,}\\
1&\text{if $i_{r+1}= i_r-1$,}\\
1+(i_{r+1}-i_r+y_{r+1}-y_r)^{-1}&\text{if $|i_r-i_{r+1}|>1$.}
\end{array}\right.
\end{align*}
(The inverses on the right hand sides of these formulae make sense
because
each $y_{r+1}-y_r$ is nilpotent with nilpotency degree at most two, as is clear
from the diagrammatic definition of the $y_r$'s.)
\end{Theorem}

\begin{proof}
This is a consequence of Theorem~\ref{iso1} combined with
the main theorem of \cite{BKinv};
see also \cite[Theorem 8.5]{BS3}.
\end{proof}

Henceforth, we will use the isomorphism from the above theorem to
{\em identify}
the algebra $1_d^{p,q} H^{p,q}_d$
with $R_d^{p,q}$, so
\begin{equation}\label{ident}
1_d^{p,q} H^{p,q}_d
\equiv R_d^{p,q}.
\end{equation}
We will denote it always by the more compact notation $R_d^{p,q}$.
Thus there are three different ways of viewing $R_d^{p,q}$:
it is a diagram algebra with basis given by oriented stretched
circle diagrams, it is a cyclotomic Khovanov-Lauda-Rouquier algebra,
and it is a sum of blocks of the cyclotomic Hecke algebra $H^{p,q}_d$.

\phantomsubsection{Super version of level two Schur-Weyl duality}
Now we can prove the main results of the section, namely,
that the map $\Phi$ from (\ref{Phi}) is surjective and the induced map (\ref{rep}) is an isomorphism.
In the case $d \leq \min(m,n)$ we have already done most of the work:

\begin{Theorem}\label{step}
If $d \leq \min(m,n)$ then we have that
$\nT_d^{p,q} = \cV(\la_{p,q}) \otimes \nV^{\otimes d}$, $R_d^{p,q} =
H_d^{p,q}$,
and the map
$$
\Phi: H_d^{p,q} \rightarrow
\End_G(\cV(\la_{p,q}) \otimes \nV^{\otimes d})^{\op}
$$
is an algebra isomorphism.
\end{Theorem}

\begin{proof}
Let us first show that $\nT_d^{p,q} = \cV(\la_{p,q})
\otimes \nV^{\otimes d}$.
Observe for $d \leq \min(m,n)$ that
any $(p,q)$-admissible sequence $\bi \in \Z^d$
necessarily lies in $(I_{p,q})^d$. This is clear from (\ref{CKLR})
and the form of the diagram $\la_{p,q}$.
Hence applying Lemma~\ref{easyzero}
we get that $\cF_{\bi} \cV(\la_{p,q}) = \{0\}$
for $\bi \in \Z^d \setminus (I_{p,q})^d$.
So we are done by (\ref{littlet}).

Now consider the map
$\Phi$.
It is injective by Corollary~\ref{faith}.
To show that it is an isomorphism,
we apply Lemma~\ref{eqd},
recalling the identification (\ref{ident}),
to see that
\begin{align*}
\dim \End_G(\cV(\la_{p,q})\otimes \nV^{\otimes d})^{\op}
&= \dim \End_G(\nT_d^{p,q})^{\op}= \dim R_d^{p,q}
\leq \dim H_d^{p,q}.
\end{align*}
Hence our injective map is an isomorphism.
At the same time, we deduce that $\dim R_d^{p,q} = \dim H_d^{p,q}$,
hence $R_d^{p,q}=H_d^{p,q}$.
\end{proof}

It remains to consider the cases with
$d > \min(m,n)$. For that, following a standard argument,
we need to allow $m$ and $n$ to vary.
So we take some other integers $m', n' \geq 0$ and consider
the supergroup $G' = GL(m'|n')$. To avoid any confusion, we decorate
all notation related to $G'$ with a prime, e.g. $\nV'$ is its natural
module, its irreducible modules are the modules denoted $\cL'(\la)$
for $\la \in \La':=\La(m'|n')$, and $\sFp := \sF(m'|n')$. We are going to exploit the following
standard lemma.

\begin{Lemma}\label{ontolem}
Let $F:\sFp \rightarrow \sF$ be an exact functor
and $X \subseteq \La'$ be a subset
with the following properties:
\begin{itemize}
\item[(i)] $F$ commutes with duality, i.e. $F \circ \circledast
\cong \circledast \circ F$;
\item[(ii)] the modules $F \cV'(\la)$
for $\la \in X$ have standard flags;
\item[(iii)]
the map $\hom_{G'}(\cV'(\la), \cV'(\mu)^\circledast)
\rightarrow \hom_{G}(F \cV'(\la), F \cV'(\mu)^\circledast)$
defined by the functor $F$ is surjective
for all $\la,\mu \in X$.
\end{itemize}
Suppose $M, N$
are $G'$-modules with standard flags
all of whose sections are of the form
$\cV'(\la)$ for $\la \in X$.
Then the map
$\hom_{G'}(M,N^\circledast) \rightarrow \hom_G(FM, FN^\circledast)$
defined by the functor
$F$ is surjective.
\end{Lemma}

\begin{proof}
We proceed by induction on the sum of the lengths of the standard flags
of $M$ and $N$, the base case being covered by (iii).
For the induction step, either $M$ or $N$ has a standard flag
of length greater than one. It suffices to consider the case
when the standard flag of $M$ has length greater than $1$, since the other case reduces to that using duality.
Pick a submodule $K$ of $M$
such that both $K$ and $Q := M / K$
are non-zero and possess standard flags. By a general fact about highest weight categories
(see also \cite[Lemma 3.6]{Btilt} for a short
direct proof in this context), the functor $\hom_{G'}(?, N^\circledast)$
is exact on sequences of $G'$-modules possessing a standard flag.
So
applying it to the short exact sequence
$0 \rightarrow K \rightarrow M \rightarrow Q
\rightarrow 0$ we get a short exact sequence as on the top line of the
following diagram:
$$
\begin{CD}
0\: \rightarrow &\hom_{G'}(Q,N^\circledast) &\rightarrow &\hom_{G'}(M,N^\circledast)
&\rightarrow &\hom_{G'}(K,N^\circledast)
&\rightarrow \:0\\
&@VVV@VVV@VVV\\
0\: \rightarrow\: &\hom_{G}(F Q,FN^\circledast) &\:\rightarrow\: &\hom_{G}(FM,FN^\circledast)
&\:\rightarrow\: &\hom_{G}(FK,FN^\circledast)\:
&\rightarrow \:0.
\end{CD}
$$
Similar considerations applying the functor
$\hom_G(?, FN^\circledast)$ to
$0 \rightarrow F K \rightarrow FM \rightarrow F Q \rightarrow 0$
gives the short exact sequence in the bottom row.
The vertical maps making the diagram commute are the maps defined
by the functor $F$.
The left and right vertical arrows are surjective by the induction hypothesis.
Hence the middle vertical arrow is surjective too by the
five lemma.
\end{proof}

Now we consider the situation that
\begin{equation}\label{case1}
G' = GL(m|n+1).
\end{equation}
 Embed $G = GL(m|n)$
into $G'$ in the top left hand
corner in the obvious way.
Also let $S$ be the one-dimensional torus embedded into $G'$
in the bottom right hand corner, so that $S$ centralises the subgroup $G$.
The character group $X(S)$ is generated by $\eps_{m+n+1}$.
Let
$F_+:\sFp \rightarrow \sF$
be the functor mapping $M \in \sFp$ to
the $-(q+m)\eps_{m+n+1}$-weight space
of $\Pi^{q+m} M$ with respect to the torus $S$.

\begin{Lemma}\label{res1}
Let $G' = GL(m|n+1)$ and $G = GL(m|n)$ as in (\ref{case1}).
\begin{itemize}
\item[(i)]
For $\la \in X^+(T')$
with $(\la,\eps_{m+n+1}) = q+m$, we have that
$F_+ \cV'(\la) \cong \cV(\mu)$, where $\mu$
is the restriction of $\la$ to $T < T'$.
\item[(ii)] $F_+ \cV'(\la'_{p,q}) \cong \cV(\la_{p,q})$.
\item[(iii)]
For $\la \in X^+(T)$
with $(\la,\eps_{m+n+1}) < q+m$, we have that
$F_+ \cV'(\la) = \{0\}$.
\end{itemize}
\end{Lemma}

\begin{proof}
The proof of (i) reduces easily using the definition (\ref{kacmod})
to checking that the $-(q+m)\eps_{m+n+1}$-weight space
of $\Pi^{q+m}E'(\la)$ with respect to $S$
is isomorphic to $E(\mu)$ as a $G_{\0}$-module, which is
well known.
Then (ii) is a consequence of (i), noting that
$\la_{p,q}' = \sum_{r=1}^m p \eps_r - \sum_{s=1}^{n+1} (q+m) \eps_{m+s}$.
The proof of (iii) is similar to (i).
\end{proof}

\begin{Lemma}\label{uniq}
Let $G' = GL(m|n+1)$ and $G = GL(m|n)$ as in (\ref{case1}).
There is a unique $G$-module
isomorphism
$$
F_+ (\cV'(\la'_{p,q}) \otimes (\nV')^{\otimes d})
\stackrel{\sim}{\rightarrow} \cV(\la_{p,q}) \otimes \nV^{\otimes d}
$$
such that
$v_{p,q}' \otimes v'_{i_1} \otimes \cdots \otimes v'_{i_d}
\mapsto
v_{p,q} \otimes v_{i_1} \otimes \cdots \otimes v_{i_d}$
for all $1 \leq i_1,\dots,i_d \leq m+n$.
Moreover, this map intertwines the natural actions of $H^{p,q}_d$.
\end{Lemma}

\begin{proof}
The first statement follows using the isomorphism
$F_+ \cV'(\la'_{p,q}) \cong \cV(\la_{p,q})$ from the first part of the previous lemma,
together with the following observations:
\begin{itemize}
\item the weights $\mu$ arising
with non-zero multiplicity in $\cV'(\la'_{p,q})$
all satisfy $(\mu, \eps_{m+n+1}) \leq q+m$;
\item the weights $\mu$ arising with non-zero multiplicity
in $(\nV')^{\otimes d}$ all satisfy $(\mu, \eps_{m+n+1}) \leq 0$;
\item the zero weight space of
$(\nV')^{\otimes d}$ with respect to $S$ is $\nV^{\otimes d}$.
\end{itemize}
The second statement is straightforward.
\end{proof}

\begin{Lemma}\label{red1}
Let $G' = GL(m|n+1)$ and $G = GL(m|n)$ as in (\ref{case1}).
Assume the map
$\Phi':H^{p,q}_d \rightarrow \End_{G'}(\cV'(\la'_{p,q})
\otimes
{(\nV')}^{\otimes d})^{\op}$
is surjective.
Then the map
$\Phi:H^{p,q}_d \rightarrow \End_{G}(\cV(\la_{p,q}) \otimes
\nV^{\otimes d})^{\op}$
is surjective too.
\end{Lemma}

\begin{proof}
We apply Lemma~\ref{ontolem}
to $F := F_+$, taking
$m':=m, n':=n+1$ and
the set of weights $X$ to be
$\{\la \in X^+(T')\:|\:(\la,\eps_{m+n+1}) \leq q+m\}$.
The hypothesis in Lemma~\ref{ontolem}(ii)
follows from Lemma~\ref{res1}, and the other two hypotheses are clear.
Since $\cV'(\la'_{p,q}) \otimes (\nV')^{\otimes d}$ is self-dual
and has a standard flag,
we deduce that the functor $F_+$ defines a surjection
$$
\End_{G'}(\cV'(\la'_{p,q}) \otimes (\nV')^{\otimes d})^{\op}
\twoheadrightarrow
\End_G(F_+ (\cV'(\la'_{p,q}) \otimes (\nV')^{\otimes d}))^{\op}.
$$
Composing with the isomorphism from Lemma~\ref{uniq}
and using also the last part of that
lemma, we deduce that there is a commutative
triangle
$$
\begin{CD}
&\!H^{p,q}_d&\vspace{-2mm}\\
\:\:\qquad\qquad\qquad\qquad\qquad{^{\Phi'}}\!\!\!\swarrow&&\!\!\searrow\!^{\Phi}\\
\End_{G'}(\cV'(\la'_{p,q}) \otimes (\nV')^{\otimes d})^{\op}
@>>>
&\End_G(\cV(\la_{p,q}) \otimes \nV^{\otimes d})^{\op}
\end{CD}
$$
in which the bottom map is surjective.
The map $\Phi'$ is surjective by assumption. So we deduce that
$\Phi$ is surjective too.
\end{proof}

Instead consider the situation that
\begin{equation}\label{case2}
G' = GL(m+1|n)
\end{equation}
and embed $G = GL(m|n)$ into $G'$ into the bottom right hand
corner in the obvious way.
Also let $S$ be the one-dimensional torus embedded into $G'$
in the top left hand
corner, so again $S$ centralises the subgroup $G$.
The character group $X(S)$ is generated by $\eps_{1}$. Let
$F_-:\sFp \rightarrow \sF$
be the functor mapping $M \in \sFp$ to
the $(p-n)\eps_{1}$-weight space
of $M$ with respect to the torus $S$.
The analogues of Lemmas~\ref{res1}--\ref{red1} in the new situation are as follows.

\begin{Lemma}\label{res2}
Let $G' = GL(m+1|n)$ and $G = GL(m|n)$ as in (\ref{case2}).
\begin{itemize}
\item[(i)]
For $\la \in X^+(T')$
with $(\la,\eps_{1}) = p$, we have that
$F_- \cV'(\la) \cong \cV(\mu)$, where $\mu$
is the restriction of
$\la -n \eps_1 + \eps_{m+2}+\cdots+\eps_{m+n+1}$ to $T$.
\item[(ii)] $F_- \cV'(\la'_{p,q}) \cong \cV(\la_{p,q})$.
\item[(iii)]
For $\la \in X^+(T')$,
with $(\la,\eps_{1}) > p$, we have that
$F_- \cV'(\la) = \{0\}$.
\end{itemize}
\end{Lemma}

\begin{proof}
This follows by similar
arguments to the proof of
Lemma~\ref{res1}, but there is an additional subtlety.
The main new point is that if $v_+$ is a non-zero
highest weight vector of weight $\la$
in $\cV'(\la)$ as in (i),
then the vector $e_{1,m+2} \cdots e_{1,m+n+1} v_+$ gives
a highest weight vector for $G$ in $F_- \cV'(\la)$
of weight
$\la - n \eps_1+\eps_{m+2}+\cdots+\eps_{m+n+1}$. This statement
is checked by explicit calculation in $U(\mathfrak{g}')$. It then
follows from the PBW theorem
that this vector generates $F_- \cV'(\la)$ and $F_- \cV'(\la) \cong
\cV(\mu)$
to give (i). For (ii) we note that $\la'_{p,q} = \sum_{r=1}^{m+1} p
\eps_r - \sum_{s=1}^n (q+m+1)\eps_{m+1+s}$.
\end{proof}

\begin{Lemma}\label{non2}
Let $G' = GL(m+1|n)$ and $G = GL(m|n)$ as in (\ref{case2}).
There is a unique $G$-module
isomorphism
$$
F_- (\cV'(\la'_{p,q}) \otimes (\nV')^{\otimes d})
\stackrel{\sim}{\rightarrow} \cV(\la_{p,q}) \otimes \nV^{\otimes d}
$$
such that
$e_{1,m+2} \cdots e_{1,m+n+1} \cdot
v_{p,q}' \otimes v'_{i_1} \otimes \cdots \otimes v'_{i_d}
\mapsto
v_{p,q} \otimes v_{i_1-1} \otimes \cdots \otimes v_{i_d-1}$
for all $2 \leq i_1,\dots,i_d \leq m+n+1$.
Moreover, this map intertwines the natural actions of $H^{p,q}_d$.
\end{Lemma}

\begin{proof}
Similar to the proof of Lemma~\ref{uniq}, using Lemma~\ref{res2}.
\end{proof}

\begin{Lemma}\label{red2}
Let $G' = GL(m+1|n)$ and $G = GL(m|n)$ as in (\ref{case2}).
Assume the map
$\Phi':H^{p,q}_d \rightarrow \End_{G'}(\cV'(\la'_{p,q})
\otimes
{(\nV')}^{\otimes d})^{\op}$
is surjective.
Then the map
$\Phi:H^{p,q}_d \rightarrow \End_{G}(\cV(\la_{p,q}) \otimes
\nV^{\otimes d})^{\op}$
is surjective too.
\end{Lemma}

\begin{proof}
Apply Lemma~\ref{ontolem}
taking $m' :=m+1, n' := n$
and the set $X$ of weights to be
$\{\la \in X^+(T')\:|\:(\la,\eps_1) \geq p\}$,
 arguing in the same way as in the proof of Lemma~\ref{red1}.
\end{proof}

Finally we can assemble the pieces to prove a key result, which
is a super analogue of the Schur-Weyl duality for level two
from \cite{BKschur}.

\begin{Theorem}[Super Schur-Weyl duality]\label{endo}
For any $d \geq 0$, the map
$$
\Phi:H^{p,q}_d \rightarrow \End_G(\cV(\la_{p,q})\otimes \nV^{\otimes d})^{\op}
$$
is surjective.
\end{Theorem}

\begin{proof}
In the case that $d \leq \min(m,n)$, this is immediate from Theorem~\ref{step}.
To prove it in general, pick $m \leq m'$ and $n \leq n'$ so that $d \leq \min(m',n')$.
We already know the surjectivity of the
map $\Phi'$
for $G' = GL(m'|n')$.
Now apply Lemma~\ref{red1} a total of $(n'-n)$ times and Lemma~\ref{red2}
a total of $(m'-m)$ times to
deduce the surjectivity for $G = GL(m|n)$.
\end{proof}

\begin{Corollary}\label{endoc}
Recalling the identification (\ref{ident}),
the map $\Phi$ induces an algebra isomorphism
$R_d^{p,q}
\stackrel{\sim}{\rightarrow} \End_G(\nT_d^{p,q})^{\op}$.
\end{Corollary}

\begin{proof}
Theorem~\ref{endo} shows the induced map
$R_d^{p,q}\rightarrow \End_G(\nT_d^{p,q})^{\op}$
 is surjective. It is an isomorphism by Lemma~\ref{eqd}.
\end{proof}

\phantomsubsection{\boldmath Irreducible representations of $R_d^{p,q}$}
As an application of Corollary~\ref{endoc}, we can recover the known classification
of the irreducible $R_d^{p,q}$-modules.
For $\la \in \La^\circ_{p,q}$ with $\deg(\la) =
\deg(\la_{p,q})+d$, let
\begin{equation}\label{irreeps}
\cD^{p,q}(\la) := \hom_G(\nT_d^{p,q}, \cL(\la)),
\end{equation}
viewed as a left $R^{p,q}_d$-module in the natural way.

\begin{Theorem}\label{cycloclass}
The modules
$\{\cD^{p,q}(\la)\:|\:\la \in \La^\circ_{p,q},
\deg(\la) = \deg(\la_{p,q})+d\}$
give a complete set of pairwise inequivalent irreducible $R^{p,q}_d$-modules.
Moreover, we have that $\dim \cD^{p,q}(\la) = \dd(\la)$, where
$\dd(\la)$
is as defined in Theorem~\ref{Dec}.
\end{Theorem}

\begin{proof}
As $\nT_d^{p,q}$ is a projective module,
Corollary~\ref{endoc} and the usual
theory of functors of the form $\hom(P,?)$
imply that the non-zero modules of the form
$\hom_G(\nT_d^{p,q}, \cL(\la))$ for $\la \in \La$ give a complete set of pairwise
inequivalent irreducible $R^{p,q}_d$-modules.
The non-zero ones are parametrised by the weights $\la \in \La^\circ_{p,q}$
with $\deg(\la) =\deg(\la_{p,q})+d$, thanks to Corollary~\ref{summands}.
Finally the stated dimension formula is a consequence of Theorem~\ref{Dec}.
\end{proof}

\begin{Remark}\rm
For a graded version of the dimension formula for the irreducible $R^{p,q}_d$-modules
derived in Theorem~\ref{cycloclass}, we refer the reader to \cite[Theorem 9.9]{BS3}.
(The identification of the labellings of irreducible representations in the above theorem with the one in \cite{BS3}
can be deduced using the methods of the next subsection.)
\end{Remark}

\phantomsubsection{\boldmath $i$-Restriction and $i$-induction}
To identify the labelling of irreducible
$R^{p,q}_d$-modules from Theorem~\ref{cycloclass} with
other known
parametrisations, it is useful to have available a more intrinsic
characterisation of $\cD^{p,q}(\la)$. We explain one inductive approach to this here in
terms of the well-known {\em $i$-restriction functors}.

Suppose that $i \in I_{p,q}$.
The natural inclusion $H_{d}
\hookrightarrow H_{d+1}$ induces an
embedding $H_{d}^{p,q} \hookrightarrow H_{d+1}^{p,q}$.
Composing before and after with the inclusion
$R_{d}^{p,q} = 1_{d}^{p,q} H_{d}^{p,q}
\hookrightarrow H_{d+1}^{p,q}$ and the projection
$H_{d+1}^{p,q} \twoheadrightarrow 1_{d+1}^{p,q} H_{d+1}^{p,q} = R_{d+1}^{p,q}$,
we get a unital algebra homomorphism
\begin{equation}
\theta_d:R^{p,q}_{d} \rightarrow R^{p,q}_{d+1}.
\end{equation}
Note this map need not be injective, e.g. if $d =(m+n)(p-q)+2mn$
then
the algebra
$R^{p,q}_{d}$ is non-zero but $R^{p,q}_{d+1}$ is zero
by Lemma~\ref{eqd}.
The image of $x_{d+1}$ in $R^{p,q}_{d+1}$ centralises
$\theta_d(R^{p,q}_d)$.
So it makes sense to define the {\em $i$-restriction functor}
\begin{equation}\label{ires}
\cE_i:\Rep{R^{p,q}_{d+1}}
\rightarrow \Rep{R^{p,q}_{d}}
\end{equation}
to be the exact functor mapping an $R^{p,q}_{d+1}$-module $M$ to the
generalised $i$-eigenspace of $x_{d+1}$ on $M$, viewed as an
$R^{p,q}_{d}$-module via
$\theta_d$.

For us, a slightly different formulation of this definition will be
more convenient.
Let
\begin{equation}\label{resid}
1^{p,q}_{d;i}
:= \sum_{\bi \in (I_{p,q})^{d+1},\,i_{d+1} = i} e(\bi)\in
R^{p,q}_{d+1}.
\end{equation}
Multiplication by this idempotent projects any $R^{p,q}_{d+1}$-module $M$ onto
the generalised $i$-eigenspace of $x_{d+1}$, which is a module over
the subring $1^{p,q}_{d;i} R^{p,q}_{d+1} 1^{p,q}_{d;i}$ of $R^{p,q}_{d+1}$.
As $1^{p,q}_{d;i}$ centralises the image of the
homomorphism $\theta_d$, we can define another unital algebra homomorphism
\begin{equation}\label{thetai}
\theta_{d;i}: R^{p,q}_{d} \rightarrow 1^{p,q}_{d;i} R^{p,q}_{d+1} 1^{p,q}_{d;i},
\qquad x \mapsto \theta_d(x) 1^{p,q}_{d;i}.
\end{equation}
Because $1^{p,q}_{d+1} = \sum_{i \in I_{p,q}} 1^{p,q}_{d;i}$, we have that
$\theta_d = \sum_{i \in I_{p,q}} \theta_{d;i}$.
The functor $\cE_i$ from the previous paragraph can be defined
equivalently
as the functor mapping an $R^{p,q}_{d+1}$-module $M$
to the space $1^{p,q}_{d;i} M$ viewed as an $R^{p,q}_{d}$-module
via the homomorphism $\theta_{d;i}$.
So:
\begin{equation}\label{anotherde}
\cE_i M = 1^{p,q}_{d;i} M \cong \hom_{R^{p,q}_{d+1}}(R^{p,q}_{d+1}
1^{p,q}_{d;i}, M),
\end{equation}
where we view $R^{p,q}_{d+1} 1^{p,q}_{d;i}$ as an $(R^{p,q}_{d+1},
R^{p,q}_{d})$-bimodule using the homomorphism $\theta_{d;i}$ to get the
right module structure.
It is clear from (\ref{anotherde}) that the $i$-restriction functor
$\cE_i$ has a left adjoint
\begin{equation}\label{iind}
\cF_i := R^{p,q}_{d+1} 1^{p,q}_{d;i} \otimes_{R^{p,q}_{d}} ?:
\Rep{R^{p,q}_{d}} \rightarrow \Rep{R^{p,q}_{d+1}}.
\end{equation}
We refer to this as the {\em $i$-induction functor}.

\begin{Lemma}\label{rfun}
There is an isomorphism
$$
r:
\cE_i \circ \hom_G(\nT_{d+1}^{p,q},?)
\stackrel{\sim}{\rightarrow}
\hom_G(\nT_{d}^{p,q},?) \circ \cE_i
$$ of functors from $\sF$
to $\Rep{R^{p,q}_{d}}$.
\end{Lemma}

\begin{proof}
Take $M \in \sF$.
Note recalling (\ref{wtdec}) that
$\nT_{d+1}^{p,q} 1^{p,q}_{d;i} = \cF_i \nT_{d}^{p,q}$.
So we can identify
\begin{align*}
\cE_i (\hom_G(\nT_{d+1}^{p,q}, M))
&= 1^{p,q}_{d;i}
\hom_G(\nT_{d+1}^{p,q}, M)\\
&=
\hom_G(\nT_{d+1}^{p,q} 1^{p,q}_{d;i}, M)
=
\hom_G(\cF_i \nT_{d}^{p,q},M).
\end{align*}
Then the adjunction between $\cF_i$ and $\cE_i$
fixed earlier defines a natural isomorphism
$\hom_G(\cF_i \nT_{d}^{p,q},M)
\stackrel{\sim}{\rightarrow}
\hom_G(\nT_{d}^{p,q}, \cE_i M)$.
Naturality gives automatically that this is an $R^{p,q}_{d}$-module homorphism.
So we have defined the desired isomorphism of functors $r$.
\end{proof}

\begin{Corollary}\label{branching}
Take $\la \in \La^\circ_{p,q}$
with $\deg(\la) = \deg(\la_{p,q})+d+1$ for some $d \geq 0$.
ick $\mu \in \La^\circ_{p,q}$ such that
$\mu \stackrel{i}{\rightarrow}\la$ is an edge in the crystal graph
for some $i \in I_{p,q}$.
Then $\cD^{p,q}(\la)$ is the unique irreducible representation of
$R^{p,q}_{d+1}$ with the property that
$\cE_i \cD^{p,q}(\la)$ has a quotient isomorphic to $\cD^{p,q}(\mu)$.
\end{Corollary}

\begin{proof}
By Lemma~\ref{ga2}, we know for $\la$ as in the statement of the corollary and
$i \in I_{p,q}$ that $\cE_i \cL(\la)$ is zero unless there
exists $\mu \in \La^\circ_{p,q}$ with $\mu
\stackrel{i}{\rightarrow} \la$ in the crystal graph,
in which case $\cL(\mu)$ is the unique irreducible quotient of
$\cE_i \cL(\la)$.
The corollary follows from this on
applying the exact functor $\hom_G(\nT_{d}^{p,q},?)$ and using
Lemma~\ref{rfun} and the definition (\ref{irreeps}).
\end{proof}

\iffalse
\begin{Remark}\rm
There is also an analogue of Lemma~\ref{rfun} for $\cF_i$:
we have that
$
\hom_G(\nT_{d}^{p,q},?) \circ \cF_i
\cong \cF_i \circ \hom_G(\nT_{d-1}^{p,q},?)$ as functors from $\sF$
to $\Rep{R^{p,q}_{d}}$. We don't need this here so omit the proof.
\end{Remark}
\fi

\begin{Remark}\label{sometimes}\rm
It is sometimes necessary to understand the homomorphism
$\theta_{d;i}$ from (\ref{thetai}) from a diagrammatic point of view. Using Theorem~\ref{iso2},
we can easily write down
$\theta_{d;i}$ on the Khovanov-Lauda-Rouquier
generators: it is the map
\begin{equation}
\theta_{d;i}:e(\bi) \mapsto e(\bi+i),\qquad
y_r e(\bi) \mapsto y_r e(\bi+i),\qquad
\psi_s e(\bi) \mapsto \psi_s e(\bi+i)
\end{equation}
for $\bi \in (I_{p,q})^{d}$,
$1 \leq r \leq d$ and $1 \leq s < d$,
where $\bi+i$ denotes the $(d+1)$-tuple $(i_1,\dots,i_{d},i)$.
It is harder to see $\theta_{d;i}$
in terms of the bases of
oriented stretched circle diagrams, but this is
worked out in detail in \cite[Corollary 6.12]{BS3}.
The basic idea to compute $\theta_{d;i}(|\bu^*[\bga^*] \wr \bt[\bde]|)$
is to insert two extra levels
chosen from

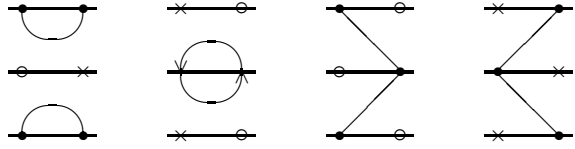
\begin{equation*}
% scriptstyles at (1,x) x = -68,-45,-22,1,24,47,70,93,...
% \times is -0.2
% \circ is +1
% \up is +0.3 and down 2.6
% \down is +0.3 and up 2.1
\begin{picture}(45,33)
\thicklines\put(-84,5){\line(1,0){33}}\thinlines
\put(-84,29){\line(1,0){33}}
\put(-67.5,29){\oval(23,23)[b]}
\put(-80.9,27.1){{$\scriptstyle\bullet$}}
\put(-57.9,27.1){{$\scriptstyle\bullet$}}
\put(-84,-19){\line(1,0){33}}
\put(-67.5,-19){\oval(23,23)[t]}
\put(-81.8,2.4){$\circ$}
\put(-59.2,3.1){$\scriptstyle\times$}
\put(-81.1,-20.9){{$\scriptstyle\bullet$}}
\put(-58.1,-20.9){{$\scriptstyle\bullet$}}

\put(-24,29){\line(1,0){33}}
\thicklines\put(-24,5){\line(1,0){33}}\thinlines
\put(-24,-19){\line(1,0){33}}
\put(-7.5,5){\oval(23,23)[b]}
\put(-7.5,5){\oval(23,23)[t]}
\put(-22.3,-20.9){$\scriptstyle\times$}
\put(1.2,-21.6){$\circ$}
\put(-22.3,27.1){$\scriptstyle\times$}
\put(1.2,26.4){$\circ$}
\put(-21.8,4.7){{$\scriptstyle\down$}}
\put(1.2,.8){{$\scriptstyle\up$}}

\thicklines\put(36,5){\line(1,0){33}}\thinlines
\put(36,29){\line(1,0){33}}
\put(61.2,26.4){$\circ$}
\put(64,5.8){\line(-1,1){22.9}}
\put(38.9,27.1){{$\scriptstyle\bullet$}}
\put(36,-19){\line(1,0){33}}
\put(61.2,-21.6){$\circ$}
\put(38.2,2.4){$\circ$}
\put(41.5,-18.2){\line(1,1){22.9}}
\put(38.9,-20.9){{$\scriptstyle\bullet$}}
\put(61.9,3.1){{$\scriptstyle\bullet$}}

\thicklines\put(96,5){\line(1,0){33}}\thinlines
\put(96,-19){\line(1,0){33}}
\put(120.7,3.1){$\scriptstyle\times$}
\put(97.7,-20.9){$\scriptstyle\times$}
\put(124.5,-19.2){\line(-1,1){22.9}}
\put(121.9,-20.9){{$\scriptstyle\bullet$}}
\put(98.9,3.1){{$\scriptstyle\bullet$}}
\put(96,29){\line(1,0){33}}
\put(97.7,27.1){$\scriptstyle\times$}
\put(101.1,5.8){\line(1,1){22.9}}
\put(121.9,27.1){{$\scriptstyle\bullet$}}
\end{picture}
\end{equation*}
\vspace{3mm}

\noindent
into the middle of the matching $\bu^* \bt$,
where we display only the strip between the $i$th
and $(i+1)$th vertices, the diagrams being trivial outside that strip.
In the first configuration here, this process involves making one
application of the generalised surgery procedure.
The construction is made precise in the paragraph after
\cite[(6.34)]{BS3}.
\end{Remark}

\section{Morita equivalence with generalised Khovanov algebras}

Next we construct an explicit Morita
equivalence between $R^{p,q} := \bigoplus_{d \geq 0}
R_d^{p,q}$ and a certain
generalised Khovanov algebra $K^{p,q}$.
Using this, we
replace the tensor space $T^{p,q} := \bigoplus_{d \geq 0} T^{p,q}_d$ from the
level two
Schur-Weyl duality
with a new space $P^{p,q}$
whose endomorphism algebra is $K^{p,q}$.
Exploiting the fact that $K^{p,q}$ is a basic algebra, we show that
the space $P^{p,q}$ has exactly
the same indecomposable summands as $T^{p,q}$ (up to isomorphism), but that they each
appear with multiplicity one.

\phantomsubsection{Generalised Khovanov algebras}
Given $p\leq q$, let
$K^{p,q}$ denote the subring
$e^{p,q} K e^{p,q}$ of $K$,
where
$e^{p,q}$ is the (non-central) idempotent
\begin{equation}\label{booky}
e^{p,q} := \sum_{\la \in \La_{p,q}^\circ}e_\la \in K.
\end{equation}

\begin{Lemma}\label{itsall}
The algebra $K$ is the union of the subalgebras $K^{p,q}$ for all $p \leq q$.
\end{Lemma}

\begin{proof}
This follows from (\ref{locun}) and the observation that, for any $\la,\mu \in \La$, we can find integers
$p \leq q$ such that both $\la$ and $\mu$ belong to $\La_{p,q}^\circ$.
\end{proof}

\begin{Remark}\rm
In terms of the diagram basis from (\ref{hbase}), $K^{p,q}$ has basis
\begin{equation*}
\left\{(a\la b)\:\Bigg|\:
\begin{array}{l}
\text{for all oriented circle diagrams $a \la b$ with $\la \in
  \La_{p,q}$ such that}\\
\text{cups and caps pass only through vertices in the interval $I_{p,q}^+$}
\end{array}
\right\}.
\end{equation*}
All the diagrams in the diagram basis of $K^{p,q}$
consist simply of straight lines oriented $\up$
outside of the interval $I_{p,q}^+$; these play no role when computing the multiplication.
So we can just ignore all of the diagram outside this strip without changing the algebra structure.
This shows that the algebra $K^{p,q}$
is a direct sum of the generalised Khovanov algebras
from \cite[$\S$6]{BS1} associated to the weights obtained from
$\La_{p,q}$ by erasing vertices outside the interval $I_{p,q}^+$.
\end{Remark}

\phantomsubsection{\boldmath Representations of $K^{p,q}$}
To understand the representation theory of the algebra $K^{p,q}$, we exploit the exact functor
\begin{equation}
e^{p,q}: \Rep{K} \rightarrow \Rep{K^{p,q}}
\end{equation}
arising by
left multiplication by the idempotent $e^{p,q}$; cf.
\cite[(6.13)]{BS1}.
It is easy to see that $e^{p,q} L(\la) \neq \{0\}$ if and only if $\la
\in \La_{p,q}^\circ$. Hence,
letting
$L^{p,q}(\la) := e^{p,q} L(\la)$
for $\la \in \La_{p,q}^\circ$,
the modules
\begin{equation}\label{theirrs}
\{L^{p,q}(\la)\:|\:\la \in \La_{p,q}^\circ\}
\end{equation}
give a complete set of pairwise inequivalent irreducible
$K^{p,q}$-modules.

Recalling also the $(K,K)$-bimodules $\widetilde{F}_i$ and $\widetilde{E}_i$
from (\ref{bims}), we get
$(K^{p,q},K^{p,q})$-bimodules
\begin{equation}
\widetilde{F}^{p,q}_i := e^{p,q} \widetilde{F}_i e^{p,q},
\qquad
\widetilde{E}^{p,q}_i :=
e^{p,q} \widetilde{E}_i
e^{p,q}
\end{equation}
for any $i \in I_{p,q}$. Let $F_i:= \widetilde{F}^{p,q}_i
\otimes_{K^{p,q}} ?$ and
$E_i := \widetilde{E}^{p,q}_i \otimes_{K^{p,q}} ?$ be the endofunctors
of $\Rep{K^{p,q}}$ defined by tensoring with these bimodules.

\begin{Lemma}\label{rfun2}
For any $i \in I_{p,q}$, there are isomorphisms
$F_i \circ e^{p,q} \cong e^{p,q} \circ F_i$
and
$E_i \circ e^{p,q} \cong e^{p,q} \circ E_i$
of functors from $\Rep{K}$ to $\Rep{K^{p,q}}$.
\end{Lemma}

\begin{proof}
We just explain the proof for $E_i$, since the argument for $F_i$ is
the same.
Suppose first that $P$ is any projective right $K$-module
that is isomorphic to a direct sum of summands of $e^{p,q} K$.
Then the natural multiplication map
$$
P e^{p,q} \otimes_{e^{p,q} K e^{p,q}} e^{p,q} K
\rightarrow P
$$
is an isomorphism of right $K$-modules.
This follows because it is obviously true if $P = e^{p,q} K$.
In the next paragraph, we show that $P = e^{p,q} \widetilde E_i$
satisfies the hypothesis that it is isomorphic to a direct sum of summands of
$e^{p,q} K$ as a right $K$-module.
Hence, we deduce that the multiplication map
\begin{equation}\label{nearly}
\widetilde E^{p,q}_i \otimes_{K^{p,q}} e^{p,q} K
\stackrel{\sim}{\rightarrow}
e^{p,q} \widetilde E_i
\end{equation}
is a $(K^{p,q},K)$-bimodule isomorphism. The desired
isomorphism $E_i \circ e^{p,q} \cong e^{p,q} \circ E_i$ follows at once,
since $E_i \circ e^{p,q}$ is
the functor defined by tensoring with the bimodule on the left hand
side
and $e^{p,q} \circ E_i$ is the functor defined by tensoring with the
bimodule on the right hand side of (\ref{nearly}).

It remains to show that $e^{p,q} \widetilde E_i$ is isomorphic to a
direct sum of summands of $e^{p,q} K$.
Equivalently, twisting with the obvious anti-automorphism $*$ that
reflects diagrams in a horizontal axis, we show that
$\widetilde F_i e^{p,q}$ is isomorphic to a direct sum of summands of
$K e^{p,q}$. The indecomposable summands of $K e^{p,q}$ are all of the
form
$P(\mu)$ for $\mu \in \La_{p,q}^\circ$, so using the definition
(\ref{booky}) this follows if we can show for any $\la \in \La_{p,q}^\circ$
that all indecomposable
summands of
$\widetilde F_i e_\la$
are of the form $P(\mu)$ for $\mu \in \La_{p,q}^\circ$.
As $\widetilde F_i e_\la \cong F_i P(\la)$, this follows easily from
\cite[Theorem 4.2]{BS2}, using also the assumption that $i \in I_{p,q}$.
\end{proof}

\begin{Corollary}\label{branching2}
Let $\la, \mu$ and $i$ be as in the statement of
Corollary~\ref{branching}.
Then $L^{p,q}(\la)$ is the unique irreducible representation of
$K^{p,q}$ with the property that
$E_i L^{p,q}(\la)$ has a quotient isomorphic to $L^{p,q}(\mu)$.
\end{Corollary}

\begin{proof}
This follows from Lemmas~\ref{ga} and \ref{rfun2} by the same argument
used to prove Corollary~\ref{branching}.
\end{proof}

\iffalse
\begin{Lemma}
We have that $K = \bigcup_{p,q} K^{p,q}$,
$\widetilde F_i = \bigcup_{p,q} \widetilde F_i^{p,q}$ and
$\widetilde E_i = \bigcup_{p,q} \widetilde E_i^{p,q}$,
where the
unions are taken over all pairs of integers $p\leq q$.
\end{Lemma}

\begin{proof}
This follows in all cases by considering diagram bases.
\end{proof}
\fi

\phantomsubsection{Morita bimodules}
Recall the $G$-module $\nT_d^{p,q}$ from (\ref{littlet}).
In view of
Theorem \ref{endo}, we can {identify}
its endomorphism algebra with
the algebra $R_d^{p,q}$.
Actually
it is convenient now to work with all $d$ simultaneously,
setting
\begin{align}
\nT^{p,q} &:= \bigoplus_{d \geq 0} \nT_d^{p,q},\\
R^{p,q} &:= \bigoplus_{d \geq 0} R_d^{p,q}
\equiv \End_G(\nT^{p,q})^{\op}.\label{splodge}
\end{align}
Note by Corollary~\ref{zer} and Lemma~\ref{eqd} that
$\nT^{p,q}$ and $R^{p,q}$ are both finite dimensional.

We next want to explain how the algebra $R^{p,q}$ is
Morita equivalent to the basic algebra $K^{p,q}$, by writing
down an explicit pair of bimodules
$A^{p,q}$ and $B^{p,q}$ that induce the Morita equivalence.
To do this, recall the notions of oriented
upper- and lower-stretched circle diagrams from
\cite[(6.17)]{BS3}.
They are the consistently oriented diagrams obtained by
gluing a cup diagram below an oriented stretched cap diagram,
or gluing a cap diagram above an oriented stretched cup diagram,
respectively.
Let $A^{p,q}$ and $B^{p,q}$ be the vector spaces
with bases\label{Ab}
\begin{align*}
\bigg\{(a\:\bt[\bga]|
\:&\bigg|\:
\begin{array}{l}
\text{for all oriented
upper-stretched circle diagrams $a\:\bt[\bga]$
of}\\
\text{height $d \geq 0$ such that
$\ga_0 = \la_{p,q}$ and $\ga_d \in \La_{p,q}$}
\end{array}
\bigg\},\\
\bigg\{|\bu^*[\bde^*]\:b)
\:&\bigg|\:
\begin{array}{l}
\text{for all oriented
lower-stretched circle diagrams $\bu^*[\bde^*]\:b$
of}\\
\text{height $d \geq 0$ such that
$\delta_0 = \la_{p,q}$ and $\delta_d \in \La_{p,q}$}
\end{array}
\bigg\},
\end{align*}
respectively.
We make $A^{p,q}$ into a
$(K^{p,q}, R^{p,q})$-bimodule as follows.
\begin{itemize}
\item
The left action of a basis vector $(a \la b) \in
K^{p,q}$
on $(c\:\bt[\bga]| \in A^{p,q}$ is by zero
unless $\la \sim \ga_d$ (where $\bga = \ga_d \cdots \ga_0$)
and $b=c^*$. Assuming these conditions hold, the product is computed by
drawing
$a \la b$ underneath $c\:\bt[\bga]$, then iterating the
generalised surgery procedure to smooth out the symmetric
middle section of the diagram.
\item
The right action of a basis vector
$|\bs^*[\btau^*]\wr \br[\bsigma]| \in R^{p,q}$
on $(a\:\bt[\bga]| \in A^{p,q}$ is by
zero unless $\bt = \bs$ and all mirror image pairs of internal circles
in $\bs^*[\btau^*]$ and $\bt[\bga]$ are oriented so that one is clockwise,
the other anti-clockwise.
Assuming these conditions hold, the product is computed by
drawing $a\:\bt[\bga]$ underneath
$\bs^*[\btau^*]\wr\br[\bsigma]$,
erasing all internal circles and number lines in
$\bt[\bga]$ and $\bs^*[\btau^*]$, then iterating the generalised
surgery procedure in the middle section once again.
\end{itemize}
Similarly we make $B^{p,q}$ into an
$(R^{p,q}, K^{p,q})$-bimodule.
We refer the reader to \cite[$\S$6]{BS3}
for detailed proofs (in an entirely analogous setting)
that these bimodules are well defined.

\begin{Theorem}
There are isomorphisms
\begin{align*}
\mu:A^{p,q} \otimes_{R^{p,q}} B^{p,q}\stackrel{\sim}{\rightarrow}
K^{p,q},\qquad
&\nu:B^{p,q} \otimes_{K^{p,q}} A^{p,q}\stackrel{\sim}{\rightarrow}
R^{p,q}
\end{align*}
of $(K^{p,q}, K^{p,q})$-bimodules
and of
$(R^{p,q}, R^{p,q})$-bimodules, respectively.
\end{Theorem}

\begin{proof}
This is a consequence of \cite[Theorem 6.2]{BS3} and \cite[Remark
6.7]{BS3}. These references give a somewhat indirect construction of the
desired isomorphisms $\mu$ and $\nu$. The same maps can also be constructed much
more directly by mimicking the definitions of
multiplication in the algebras $R^{p,q}$ and $K^{p,q}$, respectively.
\end{proof}

\begin{Corollary}[Morita equivalence]\label{sef}
The bimodule $B^{p,q}$ is a projective generator for $\Rep{R^{p,q}}$.
Also there is an algebra isomorphism
$K^{p,q} \stackrel{\sim}{\rightarrow}
\End_{R^{p,q}}(B^{p,q})^{\op}$
induced by the right action of $K^{p,q}$
on $B^{p,q}$.
Hence the functors
\begin{align*}
\hom_{R^{p,q}}(B^{p,q},?)&:\Rep{R^{p,q}} \rightarrow \Rep{K^{p,q}},\\
B^{p,q}\otimes_{K^{p,q}}?&:\Rep{K^{p,q}}\rightarrow \Rep{R^{p,q}}
\end{align*}
are quasi-inverse equivalences of categories.
\end{Corollary}

\begin{proof}
This follows immediately from the theorem by the usual arguments
of the Morita theory; see e.g. \cite[(3.5) Theorem]{Bass}.
\end{proof}

\phantomsubsection{\boldmath More about $i$-restriction and $i$-induction}
We will view the $i$-restriction and $i$-induction functors
$\cE_i$ and $\cF_i$ from (\ref{anotherde})--(\ref{iind})
now as endofunctors of $\Rep{R^{p,q}}$.
Summing the maps $\theta_{d;i}$ from (\ref{thetai}) over all $d \geq 0$,
we get a unital algebra homomorphism
\begin{equation}
\theta_i:R^{p,q} \rightarrow
1^{p,q}_i
R^{p,q} 1^{p,q}_i
\qquad\text{where}\qquad
1^{p,q}_i
:= \sum_{d \geq 0} 1^{p,q}_{d;i}
\end{equation}
(which makes sense as the sum has only finitely many non-zero terms).
Then $\cE_i$ is the functor defined by multiplying by the
idempotent $1^{p,q}_i$, viewing the result as an $R^{p,q}$-module
via $\theta_i$. The $i$-induction functor $\cF_i
=R^{p,q} 1^{p,q}_i \otimes_{R^{p,q}} ?$
is left adjoint to $\cE_i$;
here, we are viewing $R^{p,q} 1^{p,q}_i$ as a right $R^{p,q}$-module
via $\theta_i$.

\begin{Lemma}\label{it}
There is an isomorphism
$$
s':
B^{p,q}\otimes_{K^{p,q}}? \circ E_i
\stackrel{\sim}{\rightarrow}
\cE_i \circ B^{p,q}\otimes_{K^{p,q}}?
$$
of functors from $\Rep{K^{p,q}}$ to $\Rep{R^{p,q}}$.
\end{Lemma}

\begin{proof}
By the definitions of the various functors, it suffices to construct
an $(R^{p,q},K^{p,q})$-bimodule isomorphism
$$
B^{p,q} \otimes_{K^{p,q}} \widetilde E^{p,q}_i
\stackrel{\sim}{\rightarrow}
1^{p,q}_i B^{p,q},
$$
where $1^{p,q}_i B^{p,q}$ is viewed as a left $R^{p,q}$-module via the
homomorphism $\theta_i$.
There is an obvious multiplication map
defined on a tensor product of basis vectors of the form
$|\bu^*[\bde^*]\,b) \otimes (c \la t \mu d)$ so that it is zero
unless $c = b^*$ and $\la \sim \delta_d$ (where $\bde = \delta_d
\cdots \delta_0$), in which case it is the sum of basis vectors
obtained by applying the generalised surgery procedure to the
$bc$-part of the diagram obtained by
putting $\bu^*[\bde^*]\,b$ underneath $c \la t \mu d$.
The fact that this multiplication map
is an isomorphism of right $K^{p,q}$-modules is
a consequence of \cite[Theorem 3.5]{BS2}.
It remains to show that it is a left $R^{p,q}$-module
homomorphism.
Using the diagrammatic description of
the map $\theta_i$ from Remark~\ref{sometimes}, this reduces to checking
a statement which, on applying the anti-automorphism $*$, is equivalent to
the identity (6.38) established in the proof of \cite[Theorem 6.11]{BS3}.
\end{proof}

\begin{Corollary}\label{s}
There is an isomorphism
$$
s:
E_i \circ \hom_{R^{p,q}}(B^{p,q}, ?)
\stackrel{\sim}{\rightarrow}
\hom_{R^{p,q}}(B^{p,q}, ?)\circ \cE_i
$$
of functors from $\Rep{R^{p,q}}$ to $\Rep{K^{p,q}}$.
\end{Corollary}

\begin{proof}
In view of Corollary~\ref{sef},
the natural transformations arising from the canonical adjunction between tensor and hom
give isomorphisms of functors
\begin{align*}
\eta:\Id_{\Rep{K^{p,q}}} &\stackrel{\sim}{\rightarrow} \hom_{R^{p,q}}(B^{p,q},?) \circ
B^{p,q}\otimes_{K^{p,q}}?,\\
\eps:B^{p,q}\otimes_{K^{p,q}}? \circ
\hom_{R^{p.q}}(B^{p,q},?)&\stackrel{\sim}{\rightarrow} \Id_{\Rep{R^{p,q}}}.
\end{align*}
Now take the isomorphism from Lemma~\ref{it}, compose on the left and the right with the functor
$\hom_{R^{p,q}}(B^{p,q},?)$, then use the isomorphisms $\eta$ and $\eps$ to cancel the resulting pairs of quasi-inverse
functors.
\end{proof}

\iffalse
Applying Lemma~\ref{it}, we deduce that
the following composition of natural transformations
is an isomorphism:
\begin{align*}
E_i \circ \hom_{R^{p,q}}(B^{p,q}, ?)
&\stackrel{\eta\bid \bid}{\longrightarrow}
\hom_{R^{p,q}}(B^{p,q},?) \circ B^{p,q} \otimes_{K^{p,q}} ?
\circ E_i \circ \hom_{R^{p,q}}(B^{p,q},?)\\
&\stackrel{\bid s' \bid}{\longrightarrow}
\hom_{R^{p,q}}(B^{p,q},?) \circ \cE_i \circ B^{p,q}\otimes_{K^{p,q}}?
\circ \hom_{R^{p,q}}(B^{p,q},?)\\
&\stackrel{\bid \bid\eps}{\longrightarrow}
\hom_{R^{p,q}}(B^{p,q},?) \circ \cE_i
\end{align*}
This gives the desired isomorphism $s$.
\fi

\phantomsubsection{Identification of irreducible representations}
Now we can identify the labelling of the irreducible
$R^{p,q}$-modules
from Lemma~\ref{cycloclass} with the labelling of the irreducible $K^{p,q}$-modules from (\ref{theirrs}).

\begin{Lemma}\label{facup}
For $\la \in \La^\circ_{p,q}$,
we have that
$\hom_{R^{p,q}}(B^{p,q}, \cD^{p,q}(\la)) \cong L^{p,q}(\la)$
as $K^{p,q}$-modules.
\end{Lemma}

\begin{proof}
We first show that $L := \hom_{R^{p,q}}(B^{p,q}, \cD^{p,q}(\la_{p,q})) \cong
L^{p,q}(\la_{p,q})$.
It is obvious that $\cE_i \cD^{p,q}(\la_{p,q}) = \{0\}$
for all $i \in I_{p,q}$.
So by Corollary~\ref{s} we get that $E_i L = \{0\}$ for all $i \in
I_{p,q}$.
Combined with Corollary~\ref{branching2}, this implies that
$L \cong L^{p,q}(\mu)$ for some $\mu \in \La^\circ_{p,q}$ with $\deg(\mu) =
\deg(\la_{p,q})$,
hence $L \cong L^{p,q}(\la_{p,q})$ as $\la_{p,q}$ is the only such weight $\mu$.

Now take $\la \in \La^\circ_{p,q}$ different from $\la_{p,q}$,
so that $\deg(\la) > \deg(\la_{p,q})$. We again need to show that
$L := \hom_{R^{p,q}}(B^{p,q}, \cD^{p,q}(\la)) \cong L^{p,q}(\la)$.
Let $\mu$ and $i$ be as in
Corollary~\ref{branching}, so $\cD^{p,q}(\mu)$ is a quotient of $\cE_i
\cD^{p,q}(\la)$. We may assume by induction that $\hom_{R^{p,q}}(B^{p,q},
\cD^{p,q}(\mu)) \cong L^{p,q}(\mu)$.
Applying Corollary~\ref{s} again, we deduce that
$L^{p,q}(\mu)$ is a quotient of $E_i L$. So we get that $L \cong L^{p,q}(\la)$
by Corollary~\ref{branching2}.
\end{proof}

\phantomsubsection{Multiplicity-free version of level two Schur-Weyl duality}
Continue with  $p \leq q$.
Let
\begin{equation}\label{sp}
\nP^{p,q} := \nT^{p,q} \otimes_{R^{p,q}} B^{p,q}.
\end{equation}
This is a $(G,K^{p,q})$-bimodule, i.e. it is both a $G$-module and
right $K^{p,q}$-module so that the right action of $K^{p,q}$ is by $G$-module endomorphisms.

\begin{Theorem}\label{mfswd}
The homomorphism
$K^{p,q} \stackrel{\sim}{\rightarrow} \End_G(\nP^{p,q})^{\op}$
induced by the right action of $K^{p,q}$ on $\nP^{p,q}$ is an isomorphism.
Moreover:
\begin{itemize}
\item[(i)] There is an isomorphism
$$
\zeta:\hom_G(\nP^{p,q}, ?) \stackrel{\sim}{\rightarrow} \hom_{R^{p,q}}(B^{p,q}, ?) \circ
\hom_G(\nT^{p,q}, ?)
$$
of functors from $\sF$ to $\Rep{K^{p,q}}$.
\item[(ii)] There is an isomorphism
$$
t:E_i \circ \hom_G(P^{p,q}, ?) \stackrel{\sim}{\rightarrow}
\hom_G(P^{p,q}, ?) \circ \cE_i
$$
of functors from $\sF$ to $\Rep{K^{p,q}}$.
\item[(iii)] We have that $\hom_G(P^{p,q}, \cL(\la)) \cong L^{p,q}(\la)$ for
  each $\la \in \La^\circ_{p,q}$.
\item[(iv)] As a $G$-module, $\nP^{p,q}$ decomposes as $\bigoplus_{\la
    \in \La^\circ_{p,q}} \nP^{p,q} e_\la$ with $\nP^{p,q} e_\la \cong \cP(\la)$ for each $\la \in \La^\circ_{p,q}$.
\end{itemize}
\end{Theorem}

\begin{proof}
We have natural isomorphisms
\begin{align*}
\End_G(\nP^{p,q})^{\op}
&=
\hom_G(\nT^{p,q}\otimes_{R^{p,q}} B^{p,q}, \nT^{p,q}
\otimes_{R^{p,q}} B^{p,q})\\
&\cong
\hom_{R^{p,q}}(
B^{p,q}, \hom_G(\nT^{p,q},\nT^{p,q}
\otimes_{R^{p,q}} B^{p,q}))\\
&\cong
\hom_{R^{p,q}}(
B^{p,q}, \hom_G(\nT^{p,q},\nT^{p,q})
\otimes_{R^{p,q}} B^{p,q})\\
&\cong
\hom_{R^{p,q}}(
B^{p,q}, R^{p,q}
\otimes_{R^{p,q}} B^{p,q})
\cong \End_{R^{p,q}}(B^{p,q})^{\op}
\cong K^{p,q},
\end{align*}
using Corollary~\ref{sef}.
This proves the first statement in the theorem.

Then for (i), we use the natural isomorphisms
$$
\hom_G(\nP^{p,q}, M)
= \hom_G(\nT^{p,q} \otimes_{R^{p,q}} B^{p,q},
M)
\cong \hom_{R^{p,q}}(B^{p,q}, \hom_G(\nT^{p,q},M)).
$$
For (ii), we combine (i), Corollary~\ref{s} and Lemma~\ref{rfun};
the isomorphism $t$ is given explicitly by the natural transformation
$\zeta^{-1} \bid \circ \bid r \circ s \bid \circ \bid \zeta$.
For (iii), use Lemma~\ref{facup} and the definition (\ref{irreeps}).

Finally, consider (iv). The fact that $P^{p,q} = \bigoplus_{\la \in \La_{p,q}^\circ} P^{p,q} e_\la$
follows as the idempotents $\{e_\la\:|\:\la \in
\La^\circ_{p,q}\}$ sum to the identity in $K^{p,q}$.
Note as $B^{p,q}$ is projective as a left $R^{p,q}$-module,
it is a summand of a direct sum of copies of $R^{p,q}$ as a left
module. Hence as a $G$-module
$\nP^{p,q}$ is a summand of a direct sum
of copies of $\nT^{p,q}$. Applying Corollary~\ref{summands},
we deduce that the indecomposable summands of $\nP^{p,q}$
as a $G$-module are all of the form $\cP(\la)$ for various
$\la \in \La_{p,q}^\circ$. Moreover, for any $\la,\mu \in
\La^\circ_{p,q}$, we have that
$$
\dim \hom_G(\nP^{p,q} e_\la, \cL(\mu)) =
\dim e_\la \hom_G(\nP^{p,q}, \cL(\mu))
= \dim e_\la L^{p,q}(\mu) = \delta_{\la,\mu},
$$
using (iii) and the definition of $L^{p,q}(\mu)$.
This completes the proof.
\end{proof}

\section{Direct limits}

In this section we complete the proof of Theorem~\ref{main} by taking a limit
as $p \rightarrow -\infty$ and $q \rightarrow \infty$.

\phantomsubsection{Various embeddings}
In this subsection we fix $p' \leq p \leq q \leq q'$
such that {\em either} $p'=p-1$ and $q'=q$ {\em or}
$p'=p$ and $q'=q+1$.
By definition, the algebra $K^{p,q}$ is equal to the subring
$e^{p,q} K^{p',q'} e^{p,q}$ of $K^{p',q'}$.
So $P^{p',q'} e^{p,q}$
is a $(G, K^{p,q})$-bimodule.
The goal is to construct an isomorphism
$\pi_{p,q}^{p',q'}:P^{p,q} \stackrel{\sim}{\rightarrow} P^{p',q'}
e^{p,q}$.

Throughout the subsection, we set
\begin{equation}\label{bi}
\bi := \left\{\begin{array}{ll}
(p',p'-1,\dots,p'-m+1)&\text{if $p'=p-1$,}\\
 (q',q'+1,\dots,q'+n-1)&\text{if $q'=q+1$.}
\end{array}
\right.
\end{equation}
We have that $\bi \in (I_{p',q'})^c$ where $c := m$ if
$p'=p-1$ and $c := n$ if $q'=q+1$.
Introduce the idempotent
\begin{equation}\label{zi}
\zi_{\bi}
:= \sum_{d \geq 0} \zi_{\bi;d} \in R^{p',q'}
\qquad
\text{where}
\qquad
\zi_{\bi;d} = \sum_{\bj \in (I_{p,q})^d} e(\bi + \bj) \in
R^{p',q'}_{c+d},
\end{equation}
writing $\bi+\bj$ for the sequence $(i_1,\dots,i_c,j_1,\dots,j_d)$.
The following lemma explains how to identify $R^{p,q}$
with the subring $\zi_{\bi} R^{p',q'} \zi_{\bi}$ of $R^{p',q'}$.

\begin{Lemma}\label{map1}
Let $\bt = t_c \cdots t_1$ be the composite matching and
$\bGa = \Ga_c \cdots \Ga_0$ be the
block sequence associated to the $(p',q')$-admissible sequence $\bi$
from (\ref{bi}).
Let $\bga = \ga_c\cdots\ga_0$ be the unique sequence of weights
with $\ga_r \in \Ga_r$ for each $r$;
in particular, $\ga_0 = \la_{p',q'}$ and $\ga_c = \la_{p,q}$.
Then there is a unital algebra isomorphism
$$
\rho_{p,q}^{p',q'}:R^{p,q} \stackrel{\sim}{\rightarrow} \zi_{\bi} R^{p',q'} \zi_\bi
$$
defined on the basis of oriented stretched circle diagrams by setting
$$
\rho_{p,q}^{p',q'}(|\bs^*[\btau^*] \wr \br[\bsigma]|) := |
\bt^*[\bga^*]\wr
\bs^*[\btau^*]\wr \br[\bsigma] \wr \bt[\bga]|,
$$
i.e. we glue
$\ga_0 t_1^* \ga_1 \cdots \ga_{c-1} t_c^*$
onto the bottom and
$t_c \ga_{c-1} \cdots \ga_1 t_1 \ga_0$
onto the top
of the given diagram
$\bs^*[\btau^*] \wr \br[\bsigma]$.
Moreover, writing
$\rho^{p',q'}_{p,q} = \sum_{d \geq 0} \rho_d$ for isomorphisms
$\rho_d:R^{p,q}_d \stackrel{\sim}{\rightarrow} \zi_{\bi;d} R^{p',q'}_{c+d}
\zi_{\bi;d}$, the following two properties hold.
\begin{itemize}
\item[(i)]
On the Khovanov-Lauda-Rouquier
generators of $R^{p,q}_d$, we have that
$$
\rho_{d} (e(\bj)) = e(\bi+\bj),
\quad
\rho_{d} (\psi_r) = \zi_{\bi;d} \psi_{c+r},
\quad
\rho_{d} (y_s) = \zi_{\bi;d} y_{c+s},
$$
for $\bj \in (I_{p,q})^d, 1 \leq r < d$ and $1 \leq s \leq d$.
\item[(ii)]
On the Hecke generators of $R^{p,q}_d$, we have that
$$
\rho_{d} (s_r) = \zi_{\bi;d} s_{c+r},
\quad
\rho_{d} (x_s) = \zi_{\bi;d} x_{c+s},
$$
for $1 \leq r < d$ and $1 \leq s \leq d$.
\end{itemize}
\end{Lemma}

\begin{proof}
The existence of the isomorphism $\rho^{p',q'}_{p,q}$ is a consequence of the diagrammatic description of
the algebras $R^{p,q}$ and $\zi_\bi R^{p',q'} \zi_\bi$. One first checks by
inspecting bases that
the given linear map is a vector space isomorphism, then that it
preserves multiplication. The latter is obvious because we
have just added some extra line segments all oriented $\up$ at the top
and bottom of the diagram.

To check (i), it follows from (\ref{idmsps}) and the
definitions just before (\ref{psir}) and (\ref{yr})
that
$\rho_{d} (e(\bj)) = e(\bi+\bj)$,
$\rho_{d} (\bar\psi_r) = \zi_{\bi;d} \bar\psi_{c+r}$ and
$\rho_{d} (\bar y_s) = \zi_{\bi;d} \bar y_{c+s}$.
It remains to show that the signs (\ref{signs}) involved in passing from $\bar\psi$ to $\psi$
and from $\bar y$ to $y$ match up correctly, which amounts to
the observation that $\sigma_{p,q}^r(\bj) =
\sigma_{p',q'}^{c+r}(\bi+\bj)$ for $\bj \in (I_{p,q})^d$ and $1 \leq r \leq d$.
This follows from the identity
$$
\min(p,j_r)+\min(q,j_r) - p =
\min(p',j_r)+\min(q',j_r)-p'-\delta_{i_1,j_r}-\cdots-\delta_{i_c,j_r},
$$
which we leave as an exercise for the reader.
Then (ii) follows from (i) and Theorem~\ref{iso2}.
\end{proof}

We can make a very similar construction at the level of the bimodule
$B^{p,q}$. In the following lemma, we view $\zi_{\bi}
B^{p',q'} e^{p,q}$
as an $(R^{p,q}, K^{p,q})$-bimodule, where the left $R^{p,q}$-module
structure is defined via the isomorphism from Lemma~\ref{map1}.

\begin{Lemma}\label{map22}
Let $\bt$ and $\bga$ be as in Lemma~\ref{map1}.
There is an isomorphism of $(R^{p,q}, K^{p,q})$-bimodules
$$
\beta_{p,q}^{p',q'}:B^{p,q} \stackrel{\sim}{\rightarrow} \zi_{\bi}
B^{p',q'} e^{p,q}
$$
defined on the basis of oriented upper-stretched circle diagrams by setting
$$
\beta_{p,q}^{p',q'}(|\bu^*[\bde^*]\, b))
:=
|\bt^*[\bga^*]\wr \bu^*[\bde^*]\, b).
$$
Moreover,
the map
$$
\kappa_{p,q}^{p',q'}:R^{p',q'} \zi_{\bi} \otimes_{R^{p,q}} B^{p,q}
{\rightarrow}
B^{p',q'} e^{p,q},
\quad x \otimes b \mapsto x \beta^{p',q'}_{p,q}(b)
$$
is an isomorphism of $(R^{p',q'}, K^{p,q})$-bimodules.
\end{Lemma}

\begin{proof}
The fact that $\beta_{p,q}^{p',q'}$ is an isomorphism of vector spaces follows
by considering the explicit
diagram bases, and it is obviously a bimodule
homomorphism.
To deduce the final part of the lemma, it remains to show that the natural multiplication map
$$
R^{p',q'} \zi_{\bi} \otimes_{\zi_{\bi} R^{p',q'}\zi_{\bi}} \zi_{\bi} B^{p',q'}
e^{p,q}
\rightarrow B^{p',q'} e^{p,q}
$$
is an isomorphism.
For this, we argue like in the proof of Lemma~\ref{rfun2}, starting
from the trivial
observation that the multiplication map
$$R^{p',q'} \zi_{\bi} \otimes_{\zi_{\bi} R^{p',q'}\zi_{\bi}} \zi_{\bi} R^{p',q'}
\zi_{\bi}
\rightarrow R^{p',q'} \zi_{\bi}$$
is an isomorphism. Thus, we are reduced to showing
that all the
indecomposable summands of $B^{p',q'} e^{p,q}$
are also summands of $R^{p',q'} \zi_{\bi}$ as left
$R^{p',q'}$-modules.
By
Corollary~\ref{sef} and Lemma~\ref{facup}, we know that the indecomposable summands of $B^{p',q'} e^{p,q}$
are the projective covers of the irreducible
$R^{p',q'}$-modules
$\{\cD^{p',q'}(\la)\:|\:\la \in \La^\circ_{p,q}\}$.
Since $R^{p',q'}
\zi_{\bi}$ is projective, it just remains to check that
$$
\hom_{R^{p',q'}}(R^{p',q'} \zi_{\bi}, \cD^{p',q'}(\la))
=
\zi_{\bi}\cD^{p',q'}(\la) \neq \{0\}
$$ for $\la \in \La_{p,q}^\circ$.
By Lemma~\ref{cgt}, we can find $d \geq 0$ and a tuple $\bj \in (I_{p,q})^d$
such that $\la_{p,q}
\stackrel{j_1}{\rightarrow}\cdots\stackrel{j_d}{\rightarrow} \la$ is a
path in the crystal graph.
As $\la_{p',q'} \stackrel{i_1}{\rightarrow}
\cdots \stackrel{i_c}{\rightarrow}
\la_{p,q}$ is a path in the crystal graph too, we get by repeated
application of Corollary~\ref{branching}
that $\cE_{i_1} \cdots  \cE_{i_c} \cE_{j_1}
\cdots \cE_{j_d} \cD^{p',q'}(\la) \neq \{0\}$.
By the definition of the $i$-restriction functors, this means that
$e(\bi+\bj) \cD^{p',q'}(\la) \neq \{0\}$.
Since $\zi_{\bi} e(\bi+\bj) = e(\bi+\bj)$, this implies that
$\zi_{\bi} \cD^{p',q'}(\la) \neq \{0\}$ too.
\end{proof}

Next we explain how to identify $T^{p,q}$ with
$T^{p',q'}\zi_{\bi}$.

\begin{Lemma}\label{cooee}
There exists a (unique up to scalars) $G$-module isomorphism
$$
\tau^{p',q'}_{p,q} : T^{p,q}
\stackrel{\sim}{\rightarrow}
T^{p',q'} \zi_{\bi}
$$
such that
$\tau^{p',q'}_{p,q} = \sum_{d \geq 0} \tau_{d}$ for isomorphisms
$\tau_{d}:T^{p,q}_d \stackrel{\sim}{\rightarrow}
T^{p',q'}_{c+d} \zi_{\bi;d}$ with
$\tau_{d+1} = \sum_{k \in I_{p,q}} \cF_k(\tau_{d})$ for each $d \geq 0$.
Moreover, $\tau^{p',q'}_{p,q}$ is a homomorphism of right
$R^{p,q}$-modules, i.e. it is a $(G, R^{p,q})$-bimodule isomorphism,
where we are
viewing $T^{p',q'} \zi_{\bi}$
as a right $R^{p,q}$-module via the isomorphism
from Lemma~\ref{map1}.
\end{Lemma}

\begin{proof}
We first construct the map $\tau_{0}$.
Recall that $T^{p,q}_0 = \cV(\la_{p,q})$ and
$T^{p',q'}_c \zi_{\bi;0} = (\cV(\la_{p',q'}) \otimes V^{\otimes
  c}) e(\bi) = \cF_{\bi} \cV(\la_{p',q'})$.
By Lemma~\ref{ga2} we have that
$\cF_{\bi} \cV(\la_{p',q'}) \cong \cV(\la_{p,q})$; only
the analogues of the statements from (i) and (iii) of Lemma~\ref{ga}
are needed to see this.
So we can pick a $G$-module isomorphism
$$
\tau_{0}:T^{p,q}_0 = \cV(\la_{p,q}) \stackrel{\sim}{\rightarrow}
\cF_{\bi} \cV(\la_{p',q'}) = T^{p',q'}_c \zi_{\bi;0}.
$$
This map is unique up to a scalar.

Now we inductively {define} the higher
$\tau_{d}$'s. Note as
$T^{p,q}_d = \bigoplus_{\bj \in (I_{p,q})^d} \cF_{\bj} \cV(\la_{p,q})$
that $T^{p,q}_{d+1} = \bigoplus_{k \in I_{p,q}} \cF_k T^{p,q}_d$.
Similarly $T^{p',q'}_{c+d+1} \zi_{\bi;d+1}
= \bigoplus_{k \in I_{p,q}} \cF_k (T^{p',q'}_{c+d} \zi_{\bi;d})$.
So given a $G$-module isomorphism
$\tau_{d}:T^{p,q}_d \stackrel{\sim}{\rightarrow}
T^{p',q'}_{c+d} \zi_{\bi;d}$ for some $d \geq 0$, we get a
$G$-module isomorphism
$\tau_{d+1}:T^{p,q}_{d+1} \stackrel{\sim}{\rightarrow}
T^{p',q'}_{c+d+1} \zi_{\bi;d+1}$ on applying the functor
$\bigoplus_{k \in I_{p,q}} \cF_k$.
Starting from the map $\tau_0$ from the previous paragraph,
we obtain isomorphisms $\tau_d$ for every $d \geq 0$ in this way.
Then we set $\tau_{p,q}^{p',q'} := \sum_{d \geq 0} \tau_d$,
to get the desired $G$-module isomorphism $T^{p,q} \stackrel{\sim}{\rightarrow}
T^{p',q'}\xi_\bi$.

It remains to check that each
$\tau_{d}$ is a homomorphism of right $R^{p,q}_d$-modules, viewing $T^{p',q'}_{c+d} \zi_{\bi;d}$ as a
right $R^{p,q}_d$-module via the isomorphism $\rho_{d}$ from
Lemma~\ref{map1}.
Because $\tau_{0}$ is a $G$-module homomorphism,
the map
$$
\cV(\la_{p,q})\otimes V^{\otimes d}
\rightarrow \cV(\la_{p',q'}) \otimes V^{\otimes (c+d)},\:
u \otimes v \mapsto \tau_{0}(u) \otimes v\ \text{($u \in
\cV(\la_{p,q}), v \in V^{\otimes d}$)}
$$
intertwines the action of $x_s \in H_d$ with $x_{c+s} \in H_{c+d}$.
It obviously intertwines the action of each $s_r\in H_d$ with $s_{c+r} \in H_{c+d}$.
From this and the definition of $\tau_d$, we
deduce that $\tau_{d}$ intertwines the
actions of $s_r, x_s \in H_d$ on $T^{p,q}_d$ with the actions of
$s_{c+r}, x_{c+s} \in H_{c+d}$ on $T^{p',q'}_{c+d} \zi_{\bi;d}$.
So we are done by the description of $\rho_{d}$ from
Lemma~\ref{map1}(ii).
\end{proof}

Recall finally the spaces
$P^{p,q} = T^{p,q} \otimes_{R^{p,q}} B^{p,q}$ from (\ref{sp}).

\begin{Theorem}\label{maps99}
There is a unique (up to scalars) $(G, K^{p,q})$-bimodule isomorphism
$$
\pi_{p,q}^{p',q'}:P^{p,q}\stackrel{\sim}{\rightarrow} P^{p',q'} e^{p,q}
$$
such that
$\pi_{p,q}^{p',q'}(v \otimes b) = \tau_{p,q}^{p',q'}(v) \otimes \beta_{p,q}^{p',q'}(b)$
for $v \in T^{p,q}, b \in B^{p,q}$ and some choice of the isomorphism
$\tau_{p,q}^{p',q'}$ from Lemma~\ref{cooee}.
\end{Theorem}

\begin{proof}
Recalling the isomorphism $\kappa_{p,q}^{p',q'}$ from Lemma~\ref{map22}, we {define}
$\pi_{p,q}^{p',q'}$ to be the composition of the following
$(G, K^{p,q})$-bimodule isomorphisms:
\begin{align*}
T^{p,q} \otimes_{R^{p,q}} B^{p,q}
&\stackrel{\tau^{p',q'}_{p,q} \otimes \id}{\longrightarrow}
T^{p',q'} \zi_{\bi} \otimes_{R^{p,q}} B^{p,q}
\equiv
T^{p',q'} \otimes_{R^{p',q'}}
R^{p',q'}\zi_{\bi}
\otimes_{R^{p,q}} B^{p,q}\\
&\stackrel{\id \otimes \kappa^{p',q'}_{p,q}}{\longrightarrow}
T^{p',q'} \otimes_{R^{p',q'}} B^{p',q'}e^{p,q}.
\end{align*}
It remains to observe that
$\pi^{p',q'}_{p,q}(v \otimes b) =\tau^{p',q'}_{p,q}(v) \otimes
\beta^{p',q'}_{p,q}(b)$, which follows from the
definition of $\kappa^{p',q'}_{p,q}$.
\end{proof}

\phantomsubsection{Compatibility of embeddings}
Now we explain how to glue the isomorphisms $\pi_{p,q}^{p,q+1}$ and
$\pi_{p,q}^{p-1,q}$
from the previous subsection together in a consistent way
to obtain a compatible system of isomorphisms
$\pi_{p,q}^{p',q'}:P^{p,q} \stackrel{\sim}{\rightarrow} P^{p',q'}
e^{p,q}$
for every $p' \leq p \leq q \leq q'$.
The following lemma is the key ingredient making this possible.

\begin{Lemma}\label{difficult}
Let $p \leq q$ be fixed.
Given a choice of three out of the four maps
$$
\{\pi_{p,q}^{p,q+1}, \pi_{p,q}^{p-1,q},
\pi_{p,q+1}^{p-1,q+1},\pi_{p-1,q}^{p-1,q+1}\}
$$
from Theorem~\ref{maps99}, there is a unique way to choose the fourth
one so that
$\pi_{p,q+1}^{p-1,q+1}\circ \pi_{p,q}^{p,q+1}
= \pi_{p-1,q}^{p-1,q+1} \circ \pi_{p,q}^{p-1,q}$.
\end{Lemma}

\begin{proof}
We show equivalently given a choice of all four maps that there is a
(necessarily unique) scalar $z \in \C$ such that
$$
\pi_{p,q+1}^{p-1,q+1}\circ \pi_{p,q}^{p,q+1}
= z \pi_{p-1,q}^{p-1,q+1} \circ \pi_{p,q}^{p-1,q}.
$$
To see this, let
$\bh := (p-1,p-2,\dots,p-m)$
and $\bi := (q+1,q+2,\dots,q+n)$.
Let
\begin{align*}
\psi &:=
(\psi_m \psi_{m+1} \cdots \psi_{m+n-1}) \cdots (\psi_2 \psi_3 \cdots
\psi_{n+1})
(\psi_1 \psi_2 \cdots \psi_n),\\
\psi' &:=
(\psi_{n} \cdots \psi_2 \psi_1)
(\psi_{n+1} \cdots \psi_3 \psi_2)
\cdots
(\psi_{m+n-1} \cdots \psi_{m+1} \psi_m).
\end{align*}
It is easy to see from the defining relations between the
Khovanov-Lauda-Rouquier generators from
\cite[(6.8)--(6.16)]{BS3} that
$\psi \xi_{\bi+\bh} = \xi_{\bh+\bi} \psi$,
$\xi_{\bi+\bh} \psi' =
\psi' \xi_{\bh+\bi}$,
and
$\psi' \psi \xi_{\bi+\bh} = \xi_{\bi+\bh}$ in
$R^{p-1,q+1}$.

Now we claim that there exists a scalar $z \in \C$ such that the
following two diagrams commute:
\begin{equation}\label{diag1}
\begin{CD}
&&B^{p,q}&
\vspace{-2mm}
\\
&\hspace{-14mm}^{\beta_{p,q}^{p-1,q}}\!\!\!\swarrow&&\searrow^{\beta_{p,q}^{p,q+1}}\hspace{-12mm}&\\
\zi_{\bh} B^{p-1,q} e^{p,q}&&&&\xi_{\bi} B^{p,q+1} e^{p,q}\\
@V\beta_{p-1,q}^{p-1,q+1}VV&&@VV\beta_{p,q+1}^{p-1,q+1}V\\
\zi_{\bi+\bh} B^{p-1,q+1} e^{p,q}
&@>\sim>L_\psi>
&\zi_{\bh+\bi}
B^{p-1,q+1} e^{p,q},\\
\end{CD}
\end{equation}

\vspace{2mm}

\begin{equation}\label{diag2}
\begin{CD}
&&T^{p,q}&
\vspace{-2mm}
\\
&\hspace{-14mm}^{\tau_{p,q}^{p-1,q}}\!\!\!\swarrow&&\searrow^{\tau_{p,q}^{p,q+1}}\hspace{-12mm}&\\
T^{p-1,q}\zi_{\bh}&&&&T^{p,q+1}\xi_{\bi}\\
@V\tau_{p-1,q}^{p-1,q+1}VV&&@VV\tau_{p,q+1}^{p-1,q+1}V\\
T^{p-1,q+1}\zi_{\bi+\bh}
&@>\sim>R_{z\psi'}>
&
T^{p-1,q+1} \zi_{\bh+\bi},\\
\end{CD}
\end{equation}
where $L_\psi(b) := \psi b$ and $R_{z \psi'} (v) := z v \psi'$.
Given the claim and recalling Lemma~\ref{maps99}, we get for any $v \otimes b \in
T^{p,q}\otimes_{R^{p,q}} B^{p,q}$ that
\begin{align*}
\pi^{p-1,q+1}_{p,q+1}(\pi^{p,q+1}_{p,q}(v \otimes b))
&=
\tau_{p,q+1}^{p-1,q+1}(\tau_{p,q}^{p,q+1}(v))
\otimes \beta_{p,q+1}^{p-1,q+1}(\beta_{p,q}^{p,q+1}(b))\\
&=
z \tau_{p-1,q}^{p-1,q+1}(\tau_{p,q}^{p-1,q}(v)) \psi'
\otimes \psi\beta_{p-1,q}^{p-1,q+1}(\beta_{p,q}^{p-1,q}(b))\\
&=
z \tau_{p-1,q}^{p-1,q+1}(\tau_{p,q}^{p-1,q}(v))
\otimes \psi'
\psi\beta_{p-1,q}^{p-1,q+1}(\beta_{p,q}^{p-1,q}(b))\\
&=
z \tau_{p-1,q}^{p-1,q+1}(\tau_{p,q}^{p-1,q}(v))
\otimes
\beta_{p-1,q}^{p-1,q+1}(\beta_{p,q}^{p-1,q}(b))\\
&=z\pi^{p-1,q+1}_{p-1,q}(\pi^{p-1,q}_{p,q}(v \otimes b)).
\end{align*}
So the lemma follows from the claim.

To prove the claim, consider first the diagram (\ref{diag1}).
The point for this is that all the $\psi_r$'s in the element $\psi$
are acting successively on the left on $\xi_{\bi+\bh} B^{p-1,q+1}
e^{p,q}$ as a sequence
of height moves
in the sense of \cite[$\S$5]{BS3}.
Combined with the diagrammatic definition from Lemma~\ref{cooee}
this is enough to see that (\ref{diag1}) commutes.
Next consider the diagram (\ref{diag2}). Here one first reduces using
the definition of the higher $\tau_d$'s in Lemma~\ref{cooee}
to checking just that the diagram commutes on restriction to
$T^{p,q}_0 = \cV(\la_{p,q})$. In that case, both of
$T^{p-1,q+1}_{m+n}\xi_{\bi+\bh}$ and
$T^{p-1,q+1}_{m+n}\xi_{\bh+\bi}$ are isomorphic to $\cV(\la_{p,q})$,
and the map defined by right multiplication by $\psi'$ is a non-zero
isomorphism.
So the diagram must commute up to a scalar
as $\End_{G}(\cV(\la_{p,q}))$ is one dimensional.
\end{proof}

\begin{Theorem}\label{braverman}
We can
choose $(G, K^{p,q})$-bimodule isomorphisms
$$
\pi_{p,q}^{p',q'}:P^{p,q} \stackrel{\sim}{\rightarrow} P^{p',q'}
e^{p,q}
$$
for all $p' \leq p \leq q \leq q'$,
in such a way that
$\pi_{p,q}^{p'',q''} = \pi_{p',q'}^{p'',q''}\circ \pi_{p,q}^{p',q'}$
whenever $p''\leq p' \leq p \leq q \leq q' \leq q''$.
\end{Theorem}

\begin{proof}
First of all we make arbitrary choices for the maps
$\pi_{p,q}^{p,q+1}$ from Theorem~\ref{maps99} for all $p \leq q$. Also we make arbitrary
choices for the maps $\pi_{p,p}^{p-1,p}$ from Theorem~\ref{maps99} for all $p$.
Then we repeatedly apply Lemma~\ref{difficult}
, proceeding by induction on $(q-p)$, to get
maps $\pi_{p,q}^{p-1,q}$ so that the following local relation holds
$$
\pi_{p,q+1}^{p-1,q+1}\circ \pi_{p,q}^{p,q+1}
= \pi_{p-1,q}^{p-1,q+1} \circ \pi_{p,q}^{p-1,q}
$$
for all $p \leq q$.
Finally we define the maps $\pi_{p,q}^{p',q'}$ in general by
setting
$\pi_{p,q}^{p',q'} :=
\pi^{p',q'}_{p'+1,q'}\circ\cdots\circ\pi^{p-1,q'}_{p,q'}\pi^{p,q'}_{p,q'-1}\circ\cdots\circ\pi_{p,q}^{p,q+1}$.
The equality
$\pi_{p,q}^{p'',q''} = \pi_{p',q'}^{p'',q''}\circ
\pi_{p,q}^{p',q'}$
follows from this definition and the local relation.
\end{proof}

\phantomsubsection{\boldmath Proof of the main theorem}
Consider the directed set $\{(p,q)\:|\:p \leq q\}$ where
$(p,q) \rightarrow (p',q')$ if $p' \leq p \leq q \leq q'$.
By Theorem~\ref{braverman}, it is possible to choose
a direct system $\{\pi_{p,q}^{p',q'}:P^{p,q}\rightarrow P^{p',q'}e^{p,q}\}$ of $(G, K^{p,q})$-bimodule isomorphisms
for every $(p,q) \rightarrow (p',q')$.
Let
\begin{equation}\label{PPP}
P := \varinjlim P^{p,q}
\end{equation}
be the corresponding direct limit taken in the category of all
$G$-modules, and denote the canonical inclusion
of each $P^{p,q}$ into $P$ by $\phi^{p,q}$.
We make $P$ into a locally unital right $K$-module as follows.
Take $x \in K$ and $v \in P$.
Recalling Lemma~\ref{itsall}, we can
choose $p \leq q$ so that $x  = e^{p,q}
x e^{p,q}$ and
$v = \phi^{p,q}(v^{p,q})$
for some $v^{p,q} \in P^{p,q}$.
Then set $vx := \phi^{p,q}(v^{p,q} x)$.

\begin{Remark}\rm
Note that $P$ is
independent of the particular choice of the maps
$\{\pi^{p',q'}_{p,q}\}$ in the sense that if $\bar P = \varinjlim
P^{p,q}$ is another such direct limit taken with respect to
maps $\{\bar \pi_{p,q}^{p',q'}\}$, then there is a unique bimodule isomorphism
$P \stackrel{\sim}{\rightarrow} \bar P$ such that
$\phi^{p,q}(v) \mapsto \bar \phi^{p,q}(v)$ for all $v \in P^{p,q}$ and $p \leq q$.
\end{Remark}

Roughly speaking, the following lemma shows that $P$ is a minimal projective
generator for the category $\sF$ (except that as $P$ is not
finite dimensional it is not actually an object in the category).

\begin{Lemma}\label{l1}
As a $G$-module,
we have that $P = \bigoplus_{\la \in \La} P e_\la$ with $P e_\la \cong
\cP(\la)$ for each $\la \in \La$.
\end{Lemma}

\begin{proof}
The first part of the lemma is immediate because $P$ is a locally
unital right $K$-module.
To show that $P e_\la \cong \cP(\la)$, we have by
the above definitions that
$P e_\la = \varinjlim (P^{p,q} e_\la)$ where the direct limit is taken
over all $p \leq q$ with $\la \in \La^\circ_{p,q}$ (so that $e_\la \in
K^{p,q}$).
Each $P^{p,q} e_\la$ is isomorphic to $\cP(\la)$ by Theorem~\ref{mfswd}(iv).
Hence the direct limit is isomorphic to $\cP(\la)$ too.
\end{proof}

Now we want to identify the algebra $K$ with the endomorphism algebra of
$P$. A little care is needed here as $P$ is an infinite direct sum.
So for any $G$-module $M$, we let
\begin{equation}
\hom^{\fin}_G(P, M)
:= \bigoplus_{\la \in \La} \hom_G(P e_\la,M) \subseteq \hom_G(P,M),
\end{equation}
which is the {\em locally finite part} of $\hom_G(P,M)$.
Note if $M$ is finite dimensional that $\hom^{\fin}_G(P,M) = \hom_G(P,M)$.
In particular, we denote $\hom^{\fin}_G(P,P)$ by $\End_G^{\fin}(P)$ and
write
$\End_G^{\fin}(P)^{\op}$ for the opposite algebra, which acts naturally
on the right on $P$ by $G$-module endomorphisms.

\begin{Lemma}\label{l2}
The right action of $K$ on $P$ defined above induces an algebra isomorphism
$K \stackrel{\sim}{\rightarrow} \End_G^{\fin}(P)^{\op}$.
\end{Lemma}

\begin{proof}
We need to show that right multiplication induces a vector space isomorphism
$e_\la K \stackrel{\sim}{\rightarrow} \hom_G(P e_\la, P)$
for each $\la \in \La$.
By definition, the right hand space is
$$
\hom_G(P e_\la, \bigcup P e^{p,q}) =
\bigcup \hom_G(P e_\la, P e^{p,q})
$$
where we can take the union
just over $p \leq q$ with $\la \in \La_{p,q}^{\circ}$.
As $P e_\la = \phi^{p,q} (P^{p,q} e_\la)$ and $P e^{p,q} = \phi^{p,q}(P^{p,q})$
for all such $p \leq q$, the first statement from Theorem~\ref{mfswd} implies that right multiplication
induces an isomorphism
$e_\la K e^{p,q}
\stackrel{\sim}{\rightarrow} \hom_G(P e_\la, P e^{p,q})$.
Taking the union and recalling Lemma~\ref{itsall},
we deduce that we do get an isomorphism
$e_\la K \stackrel{\sim}{\rightarrow} \bigcup \hom_G(P e_\la, P e^{p,q})$.
\end{proof}

Finally we record the following variation on a basic fact.

\begin{Lemma}\label{l3}
Let $B$ be a $G$-module that is also a locally unital right
$K$-module, such that the action of $K$ on $B$ is by $G$-module
endomorphisms. Let $M$ be any finite dimensional left $K$-module
and assume that $B \otimes_K M$ is finite dimensional.
Then there is a natural $G$-module isomorphism
$$
\hom^{\fin}_G(P, B) \otimes_K M \rightarrow \hom_G(P, B \otimes_K M)
$$
sending $f \otimes m$ to the homomorphism
$v \mapsto f(v) \otimes m$.
\end{Lemma}

\begin{proof}
It suffices to show that $\hom_G(Pe_\la, B) \otimes_K M \cong \hom_G(Pe_\la, B \otimes_K M)$
for each $\la \in \La$,
which is well known.
\end{proof}

Now we can prove the main result of the article, which is essentially
Theorem~\ref{main} from the introduction with the functor $\mathbb{E}$
there constructed explicitly.
The proof is a rather standard consequence of
the last three lemmas, but we include some details since
we are in a slightly unusual locally finite setting.

\begin{Theorem}\label{reallymain}
The functors $$
\hom_G(P, ?):\sF \rightarrow
\Rep{K},\qquad
P \otimes_K ?: \Rep{K} \rightarrow \sF
$$
are quasi-inverse equivalences of categories.
Moreover $P \otimes_K P(\la) \cong \cP(\la)$ for each $\la \in \La$.
\end{Theorem}

\begin{proof}
Note using Lemma~\ref{l1} that both the functors map finite dimensional modules to finite
dimensional modules, so the first statement makes sense.
Lemmas~\ref{l3} and \ref{l2} yield a natural isomorphism
$$
\hom_G(P, P \otimes_K M) \stackrel{\sim}{\rightarrow} \hom^{\fin}_G(P, P)
\otimes_K M  \cong K \otimes_K M \equiv M
$$
for any $M \in \Rep{K}$. Thus $\hom_G(P, ?) \circ P \otimes_K ? \cong
\Id_{\Rep{K}}$.
Conversely, to show that $P \otimes_K ? \circ \hom_G(P, ?) \cong
\Id_{\sF}$, we have a natural homomorphism
$$
P \otimes_K \hom_G(P, N) \rightarrow N, \qquad
v \otimes f \mapsto f(v)
$$
for every $N \in \sF$.
Because of Lemma~\ref{l1} this map is surjective. To show that it is
injective too, denote its kernel by $U$. Applying the exact functor
$\hom_G(P, ?)$, we get a short exact sequence
$$
0 \rightarrow \hom_G(P, U)
\rightarrow \hom_G(P, P \otimes_K \hom_G(P, N)) \rightarrow \hom_G(P,
N) \rightarrow 0.
$$
By the fact established just before, the middle space here is
isomorphic to $\hom_G(P, N)$, so the right hand map is an isomorphism.
Hence $\hom_G(P, U) = \{0\}$, which implies that $U = \{0\}$.
So our natural transformation is an isomorphism, and we have
established the equivalence of categories. Moreover,
$$
P \otimes_K P(\la) = P \otimes_K K e_\la \equiv P e_\la \cong \cP(\la)
$$
by Lemma~\ref{l1}.
\end{proof}

Theorem~\ref{main} from the introduction is a consequence of
Theorem~\ref{reallymain}, taking
$\mathbb E := \hom_G(P, ?)$. We have already proved that
$\mathbb E \cP(\la) \cong P(\la)$, which immediately implies that
$\mathbb E \cL(\la) \cong L(\la)$.
The fact that $\mathbb E \cV(\la) \cong V(\la)$ follows because both the
categories $\sF$ and $\Rep{K}$ are highest weight categories in which the modules
$\{\cV(\la)\}$ and $\{V(\la)\}$ give the standard modules; see
Theorem~\ref{hwc} for the former and
\cite[Theorem 5.3]{BS1} for the latter fact.

\phantomsubsection{Identification of special projective functors}
Finally we discuss briefly how to relate the special projective
functors on the two sides of our equivalence of categories.

\begin{Theorem}
For each $i \in I$, we have that
$$
E_i \cong
\hom_G(P, ?) \circ \cE_i \circ P \otimes_K ?,\qquad
F_i \cong
\hom_G(P, ?) \circ \cF_i \circ P \otimes_K ?
$$
as endofunctors of $\Rep{K}$.
\end{Theorem}

\begin{proof}
Since $F_i$ is left adjoint to $E_i$ and $\cF_i$ is left adjoint to
$\cE_i$, the second isomorphism is a consequence of the first, by unicity of adjoints. To
prove the first, we note using Lemma~\ref{l3} that there  are natural
isomorphisms
$$
\hom_G(P, \cE_i(P \otimes_K M))
\cong \hom_G(P, (\cE_i P) \otimes_K M)
\cong \hom_G^{\fin}(P, \cE_i P) \otimes_K M
$$
for any $M \in \Rep{K}$.
Hence it suffices to show that
$\hom_G^{\fin}(P, \cE_i P) \cong \widetilde E_i$
as $(K,K)$-bimodules. For this, we just sketch how to construct the appropriate
map, leaving details to the reader.
Take any $\la \in \La$ and any $p \leq q$ so that we actually have
$\la \in \La_{p,q}^{\circ}$.
When applied to the module $P^{p,q}$, the natural isomorphism from
Theorem~\ref{mfswd}(ii)
produces a $(K^{p,q}, K^{p,q})$-bimodule isomorphism
$$
\eps^{p,q}:\widetilde E_i^{p,q}
\stackrel{\sim}{\rightarrow} \hom_G(P^{p,q}, \cE_i P^{p,q}).
$$
Restricting this to $e_\la \widetilde E_i^{p,q} = e_\la
\widetilde E_i e^{p,q}$ and using $\phi^{p,q}$ to identify $P^{p,q}$
with $P e^{p,q}$, we get from this a vector space isomorphism
$$
\eps^{p,q}:e_\la \widetilde E_i e^{p,q} \stackrel{\sim}{\rightarrow}
e_\la \hom_G(P e^{p,q}, \cE_i P e^{p,q})
= \hom_G(P e_\la, \cE_i P e^{p,q}).
$$
Now one checks for $p' \leq p \leq q \leq q'$ that
$\eps^{p,q}(v) = \eps^{p',q'}(v)$ for all $v \in e_\la
\widetilde E_i e^{p,q}$; it
suffices to do this in the cases
$(p',q') = (p-1,q)$ or $(p,q+1)$.
Hence it makes sense to take the union over all $p \leq q$ to get an
isomorphism
$$
\eps:e_\la \widetilde E_i \stackrel{\sim}{\rightarrow} \hom_G(P
e_\la, \cE_i P).
$$
Taking the direct sum of these maps over all $\la \in \La$ gives
finally the desired map $\widetilde E_i\stackrel{\sim}{\rightarrow}
\hom^{\fin}_G(P, \cE_i P)$.
\end{proof}

\section*{Index of notation}

\small
\begin{tabular}{llr}
\hspace{-5mm}$G = GL(m|n)$&General linear supergroup &\pageref{supermat}\\
\hspace{-5mm}$V, V^*$&Natural representation of $G$ and its dual&\pageref{nv}\\
\hspace{-5mm}$B, T$&Standard Borel subgroup and maximal torus of $G$&\pageref{borels}\\
\hspace{-5mm}$\sF = \sF(m|n)$&Half of the category of finite
dimensional $G$-modules\phantom{pace}&\pageref{fcat}\\
\hspace{-5mm}$X^+(T)$&Dominant weights&\pageref{XT}\\
\hspace{-5mm}$\cL(\la), \cV(\la), \cP(\la)$&Irreducibles, standards
and PIMs for $G$ for $\la \in X^+(T)$&\pageref{reps}\\
\hspace{-5mm}$\cE_i, \cF_i$&Special projective functors for $G$&\pageref{spf}\\
\hspace{-5mm}$K = K(m|n)$&Generalised Khovanov algebra&\pageref{locun}\\
\hspace{-5mm}$\La = \La(m|n)$
&Diagrammatic
weights in bijection with $X^+(T)$&\hspace{6.5mm}\pageref{dwts}\\
\hspace{-5mm}$L(\la), V(\la), P(\la)$&Irreducibles, standards and PIMs
for $K$ for $\la \in \La$&\pageref{form1}\\
\hspace{-5mm}$E_i, F_i$&Special projective functors for $K$&\pageref{bims}\\
\hspace{-5mm}$\la_{p,q}$&Ground-state weight&\pageref{laab}\\
\hspace{-5mm}$H_d^{p,q}$&Cyclotomic Hecke algebra
which acts on $\cV(\la_{p,q})\otimes V^{\otimes d}$
&\pageref{cyclo}\\
\hspace{-5mm}$I_{p,q}^+$&
Index set for $pq$-strip&\pageref{intervals}\\
\hspace{-5mm}$\La_{p,q}, \La_{p,q}^\circ$&Weights, weights of maximal
defect in $pq$-strip&\pageref{lamn1}\\
\hspace{-5mm}$1_d^{p,q}$&Central idempotent in $H_d^{p,q}$
corresponding to $pq$-strip&\pageref{theidemp}\\
\hspace{-5mm}$T_d^{p,q}$&Tensor space
$(\cV(\la_{p,q})\otimes V^{\otimes d}) 1_d^{p,q}$&\pageref{lab}\\
\hspace{-5mm}$R^{p,q}_d$ &Cyclotomic KLR algebra $
\cong  1_d^{p,q} H_d^{p,q}\cong
\End_G(T^{p,q}_d)^{\op}$
&\pageref{tpr}\\
\hspace{-5mm}$T^{p,q}, R^{p,q}$&
Direct sums $\bigoplus_{d \geq 0} T^{p,q}_d$ and $\bigoplus_{d \geq 0} R^{p,q}_d$&\pageref{splodge}\\
\hspace{-5mm}$K^{p,q}$&Subring of $K$ that is Morita equivalent to
$R^{p,q}$&\pageref{booky}\\
\hspace{-5mm}$A^{p,q}, B^{p,q}$&Morita bimodules
&\pageref{Ab}\\
\hspace{-5mm}$P^{p,q}$&
Multiplicity-free projective module
$T^{p,q} \otimes_{R^{p,q}} B^{p,q}$&\pageref{sp}\\
\hspace{-5mm}$P = \varinjlim P^{p,q}$&Canonical minimal projective generator for $\sF$&\pageref{PPP}
\end{tabular}

\end{document}